\documentclass{article}
\usepackage{amsfonts}
\usepackage{mathrsfs}
\usepackage[dvips]{graphicx}

\usepackage{amsmath}
\usepackage{latexsym,amsfonts,amssymb}

\newcommand{\bfv}{\mbox{\boldmath $v$}}
\newcommand{\bfn}{\mbox{\boldmath $n$}}

\newcommand{\bfb}{\mbox{\boldmath $b$}}
\newcommand{\bfc}{\mbox{\boldmath $c$}}

\begin{document}

\title{Quasi-differentiable Banach manifold and phase-diagram of invariant parabolic differential equation in
  such manifold\thanks{This work is supported by China National
  Natural Science Foundation under the grant number 11571381}}
\author{Shangbin Cui\\[0.2cm]
  {\small School of Mathematics, Sun Yat-Sen University, Guangzhou, Guangdong 510275,}\\
  {\small People's Republic of China. E-mail:\,cuishb@mail.sysu.edu.cn}}
\date{}
 \maketitle

\begin{abstract}
  The purpose of this paper is twofold. First we study a class of Banach manifolds which are not differentiable in traditional
  sense but they are quasi-differentiable in the sense that a such Banach manifold has an embedded submanifold such that all
  points in that submanifold are differentiable and tangent spaces at those points can be defined. It follows that differential
  calculus can be performed in that submanifold and, consequently, differential equations in a such Banach manifold can be
  considered. Next we study the structure of phase diagram near center manifold of a parabolic differential equation in Banach
  manifold which is invariant or quasi-invariant under a finite number of mutually quasi-commutative Lie group actions. We prove
  that under certain conditions, near the center manifold $\mathcal{M}_c$ the underline manifold is a homogeneous fibre bundle
  over $\mathcal{M}_c$, with fibres being stable manifolds of the differential equation. As an application, asymptotic behavior
  of the solution of a two-free-surface Hele-Shaw problem is also studied.
\medskip

\medskip

\end{abstract}

\section{Introduction}
\setcounter{equation}{0}

\hskip 2em
  Classical linearized stability theorem is an important fundamental result in the theory of ordinary differential equations.
  It states that for a differential equation $x'=F(x)$ in ${\mathbb{R}}^n$ with an isolated stationary point or equilibrium
  $x^*$, i.e., $F(x^*)=0$, where $F\in C^1(O,{\mathbb{R}}^n)$ for some open subset $O$ of ${\mathbb{R}}^n$ and $x^*\in O$,
  if $s:=\displaystyle\max_{1\leqslant j\leqslant n}{\rm Re}\lambda_j<0$, where $\lambda_j$'s are all eigenvalues of $F'(x^*)$,
  then $x^*$ is asymptotically stable, whereas if $s>0$ then $x^*$ is unstable.

  The above theorem has been successfully extended to nonlinear parabolic differential equations in Banach spaces during the
  last two decades of the last century, cf., e.g., Poitier-Ferry \cite{Poi}, Lunardi \cite{Lun1}, Drangeid \cite{Dra}, and
  Da Prato and Lunardi \cite{PraLun}; see Chapter 5 of \cite{Henry} and Chapter 9 of \cite{Lun2} for
  expositions on this topic. The extended theorem states as follows: Let $X$ be a Banach space and $X_0$ an embedded Banach
  subspace of $X$ (i.e., as linear spaces $X_0$ is a subspace of $X$ and the norm of $X$ can be controlled by that of $X_0$
  when restricted to $X_0$; see Section 3). Let $O$ be an open subset of $X_0$ and $F\in C^{2-0}(O,X)$. Consider the autonomous
  differential equation
\begin{equation}
  x'=F(x)
\end{equation}
%---(1.1)---
  in $X$. We say this equation is of {\em parabolic type} at a point $x_0\in O$ if the operator $F'(x_0)\in L(X_0,X)$ is a
  sectorial operator when regarded as an unbounded linear operator in $X$ with domain $X_0$, and the graph norm
  of $X_0$ is equivalent to its own norm $\|\cdot\|_{X_0}$. If (1.1) is of parabolic type at every point in $O$ then it is
  called of parabolic type in $O$. Assume that the differential equation (1.1) is of parabolic type in $O$ and it has an
  isolated stationary point $x^*\in O$, i.e., $F(x^*)=0$. Let $s$ be the spectrum bound of $F'(x^*)$, i.e.,
  $s=\sup\{{\rm Re}\lambda:\lambda\in\sigma(F'(x^*))\}$, where as usual $\sigma(\cdot)$ denotes the spectrum of a linear
  operator. Then the following assertion holds: If $s<0$ then $x^*$ is asymptotically stable, whereas if $s>0$ then $x^*$ is
  unstable.

  In application, however, we often encounter the critical case $s=0$ which is usually caused by unisolation of the
  stationary point $x^*$, i.e., $x^*$ is not isolated but is contained in a manifold made up of stationary points. Analysis
  to such critical case has been in the scope of researchers for over fourty years. Many authors including Hale \cite{Hale1},
  Hausrath \cite{Haus}, Carr \cite{Carr}, Chow and Lu \cite{ChowLu1}, Da Prato and Lunardi \cite{Dra}, Bates and Jones
  \cite{BatJon}, Mielke \cite{Mielke}, Iooss and Vanderbauwhede \cite{IooVan} and et al made a lot of contribution on this
  topic. Investigation shows that in the case $s=0$, if in addition to the condition $\sup\{{\rm Re}\lambda:\lambda\in
  \sigma(F'(x^*))\backslash \{0\}\}<0$ some other conditions is satisfied, then the equation has a center manifold
  made by stationary points which attracts all nearing flows.

  Clearly, phase diagram of a differential equation which is invariant under some Lie group action has certain special
  structure. In our previous work \cite{Cui2}, we studied this problem for parabolic differential equation in Banach space
  which is invariant under a local Lie group action. The result of \cite{Cui2} shows that local phase diagram of such
  differential equation at a neighborhood of the stationary point usually has a nice structure induced by the Lie group
  action. We note that limited by its purpose for application to free boundary problems, the condition of ``local Lie group
  action'' can not be replaced with ``Lie group action'' in \cite{Cui2}. However, localness of the Lie group action
  makes the condition very inconvenient to verify. To remove the condition ``local'', we must appeal to the concept of Banach
  manifold. Unfortunately, some important Banach manifolds such as the Banach manifold made by bounded $C^{m+\mu}$-domains
  in $\mathbb{R}^n$, where $m$ is a positive integer and $0\leqslant\mu\leqslant 1$, are not differentiable, i.e., they do
  not have differentiable structure, so that it is impossible to study differential equations in such Banach manifolds. It
  follows that it is impossible to reform the main result of \cite{Cui2} into a nicely presented result in traditional Banach
  manifold frame such that the reformed result is applicable to free boundary problems.

  The purpose of this paper is twofold. First we study a class of Banach manifolds which are not differentiable in traditional
  sense but they are quasi-differentiable in the sense that a such Banach manifold has an embedded submanifold such that all
  points in that submanifold are differentiable and tangent spaces at those points can be defined. It follows that differential
  calculus can be performed in that submanifold and, consequently, differential equations in a such Banach manifold can be
  considered. This analysis is motivated by and aims at applications to the Banach manifold of bounded $C^{m+\mu}$-domains
  in $\mathbb{R}^n$, which is not differentiable in traditional sense but is quasi-differentiable in the sense studied in
  the present work. We note that this Banach manifold has great importance in the study of free boundary problems; the
  reason that it has not been widely used as it should be is due to the obstacle of nondifferentiability. We also note that the
  Fr\'{e}chet manifold of bounded smooth domains in $\mathbb{R}^n$ was already used long ago in the study of free boundary
  problems, cf. \cite{Ham} for instance.

  The next goal of this paper is to extend the result of \cite{Cui2} from Banach space to Banach manifold and, furthermore,
  to extend the situation of one Lie group action into more general situation of finite many mutually quasi-commutative Lie
  group actions. As an application of our abstract result obtained in this paper, we shall also make a rigorous analysis to
  asymptotic behavior of solutions of a two-free-surface Hele-Shew problem. In what follows we briefly state our main results
  in this aspect. Let us first introduce some basic concepts and notations.

  Let $\mathfrak{M}$ be a Banach manifold and $\mathfrak{M}_0$ a $C^k$-embedded Banach submanifold of $\mathfrak{M}$,
  $k\geqslant 2$, built on the Banach spaces $X$ and $X_0$, respectively, where $X_0$ is a densely embedded Banach subspace
  of $X$; see Section 2 for details of these concepts. Definitions of the concepts of quasi-differentiable Banach manifold,
  $C^1$-kernel, inner $C^2$-kernel and shell etc. are also given in Section 2.

   Let $\mathscr{F}$ be a vector field in $\mathfrak{M}$ with domain $\mathfrak{M}_0$. Consider the following differential
   equation in $\mathfrak{M}$:
\begin{equation}
  \eta'=\mathscr{F}(\eta).
\end{equation}
%---(1.2)---
  We say this equation is of {\em parabolic type} in $\mathfrak{M}_0$ if for any $\eta_0\in\mathfrak{M}_0$ there exists a
  $\mathfrak{M}_0$-regular local chart of $\mathfrak{M}$ at $\eta_0$ such that its representation in that local chart is of
  parabolic type; see Definition 3.1 for details. We say the vector field $\mathscr{F}$ is {\em Fredholm} at a point $\eta_0
  \in\mathfrak{M}_0$ if the representation $F$ of $\mathscr{F}$ in a $\mathfrak{M}_0$-regular local chart is Fredholm. We
  say the derivative of $\mathscr{F}$ at $\eta_0\in\mathfrak{M}_0$ is a {\em standard Fredholm operator} if the derivative
  of its representation $F$ at the corresponding point $x_0$ of $\eta_0$ in $X_0$ is a standard Fredholm operator, i.e.,
  the following relations hold:
$$
  \dim{\rm Ker}F'(x_0)<\infty,\quad \mbox{\rm Range}F'(x_0)\;\mbox{is closed in $X$},\quad
  \mbox{\rm codim}{\rm Range}F'(x_0)=\dim{\rm Ker}F'(x_0), \quad \mbox{and}
$$
$$
  X={\rm Ker}F'(x_0)\oplus{\rm Range}F'(x_0).
$$
  The notations ${\rm Ker}\mathscr{F}'(\eta)$ and $\sigma(\mathscr{F}'(\eta^*))$ can be similarly defined with the aid
  of the representation $F$ of $\mathscr{F}$; see Sections 3 and 4 for details.

  Let $(G,p)$ be a $\mathfrak{M}_0$-regular Lie group action to $\mathfrak{M}$, i.e., $p\in C(G\times\mathfrak{M},\mathfrak{M})$,
  the map $a\mapsto p(a,\cdot)$ is a continuous group isomorphism from $G$ onto the transformation group of $\mathfrak{M}$, and
  $p(G\times\mathfrak{M}_0)\subseteq\mathfrak{M}_0$. We say $\mathscr{F}$ is {\em quasi-invariant} under the Lie group action
  $(G,p)$ if there exists a positive-valued function $\theta$ defined in $G$ such that the following condition is satisfied:
\begin{equation}
  \mathscr{F}(p(a,\eta))=\theta(a)\partial_{\eta}p(a,\eta)\mathscr{F}(\eta), \quad \forall a\in G, \;\; \forall\eta\in\mathfrak{M}_0.
\end{equation}
%---(1.3)---
  In this case we also say $\mathscr{F}$ is {\em $\theta$-quasi-invariant} and call $\theta$ {\em quasi-invariance factor}, and
  also say the differential equation $(1.2)$ is $\theta$-quasi-invariant under the Lie group action $(G,p)$. If in particular
  $\theta(a)=1$, $\forall a\in G$, then we simply say the vector field $\mathscr{F}$ and the differential equation $(1.2)$ are
  {\em invariant} under the Lie group action $(G,p)$.

  Let $(G_i,p_i)$, $i=1,2,\cdots,N$, be a finite number of $\mathfrak{M}_0$-regular Lie group actions to $\mathfrak{M}$. We say
  these Lie group actions are {\em mutually quasi-commutative}
  if for any $1\leqslant i,j\leqslant N$ with $i\neq j$ there exists corresponding smooth function $f_{ij}:G_i\times G_j\to G_i$
  such that
\begin{equation}
   p_j(b,p_i(a,\eta))=p_i(f_{ij}(a,b),p_j(b,\eta)), \quad \forall\eta\in\mathfrak{M},\;\; \forall a\in G_i,\;\; \forall b\in G_j,
\end{equation}
%---(1.4)---
  and for every fixed $b\in G_j$, the map $a\mapsto f_{ij}(a,b)$ is bijective with a smooth inverse. Moreover, we say these Lie
  group actions are {\em fully ranked} if by denoting $g:G\times\mathfrak{M}\to\mathfrak{M}$, where
  $G=G_1\times G_2\times\cdots\times G_N$, to be the function
\begin{equation}
   g(a,\eta)=p_1(a_1,p_2(a_2,\cdots,p_N(a_N,\eta)\cdots))
\end{equation}
%---(1.5)---
  for $\eta\in\mathfrak{M}$ and $a=(a_1,a_2,\cdots,a_N)\in G$, then
\begin{equation}
   {\rm rank}\,\partial_ag(a,\eta)=n:=n_1+n_2+\cdots+n_N, \quad \forall\eta\in\mathfrak{M}, \;\; \forall a\in G,
\end{equation}
%---(1.6)---
  where $n_i=\dim G_i$, $i=1,2,\cdots,N$.

  The first main result of this paper can be briefly stated as follows:
\medskip

  {\bf Theorem 1.1}\ \ {\em Let $(\mathfrak{M},\mathscr{A})$ be a quasi-differentiable Banach manifold with a $C^1$-kernel
  $\mathfrak{M}_0$, an inner $C^2$-kernel $\mathfrak{M}_1$ and a shell $\widetilde{\mathfrak{M}}$ with the property that
  $\mathfrak{M}_0$ is a $C^2$-embedded Banach submanifold of $\widetilde{\mathfrak{M}}$. Let $(G_i,p_i)$, $i=1,2,\cdots,N$,
  be a finite number of mutually commutative, $\mathfrak{M}_0$-regular and fully ranked Lie group actions to $\mathfrak{M}$,
  or more precisely, the conditions $(L1)\sim (L7)$ in Section 3 are satisfied. Consider the initial value problem of the
  equation $(1.2)$, where $\mathscr{F}$ is a vector field in $\mathfrak{M}$ with domain $\mathcal{O}\subseteq\mathfrak{M}_0$.
  Let $\eta^*\in\mathcal{O}\cap\mathfrak{M}_1$ be a zero of $\mathscr{F}$. Assume that the following conditions are satisfied:
\begin{enumerate}
\item[]$(G_1)$\ \ The equation $(1.2)$ is of $C^2$-class parabolic type.\vspace*{-0.2cm}
\item[]$(G_2)$\ \ $\mathscr{F}$ is quasi-invariant under all Lie group actions $(G_i,p_i)$, $i=1,2,\cdots,N$. \vspace*{-0.2cm}
\item[]$(G_3)$\ \ $\mathscr{F}'(\eta^*)$ is a standard Fredholm operator. \vspace*{-0.2cm}
\item[]$(G_4)$\ \ ${\rm dim}\,{\rm Ker}\mathscr{F}'(\eta)=n$, where $n$ is as in $(1.6)$.
\item[]$(G_5)$\ \ $\sup\{{\rm Re}\lambda:\lambda\in\sigma(\mathscr{F}'(\eta^*)))\backslash\{0\}\}<0$.
\end{enumerate}
  Then we have the following assertions:

  $(1)$\ \ The set $\mathcal{M}_c=\{g(a,\eta^*):a\in G\}$ is a $n$-dimensional
  submanifold of $\mathfrak{M}_0$.

  $(2)$\ \ There is a neighborhood $\mathcal{O}$ of $\mathcal{M}_c$ in $\mathfrak{M}_0$ such that for any $\eta_0\in
  \mathcal{O}$, the solution of $(1.2)$ with initial data $\eta(0)=\eta_0$ is global, i.e., $\eta\in C([0,\infty),
  \mathfrak{M}_0)\cap C^1((0,\infty),\mathfrak{M}_0)$.

  $(3)$\ \ There exists a submanifold $\mathcal{M}_s\subseteq\mathcal{O}$ of $\mathfrak{M}_0$ of codimension $n$ passing
  $\eta^*$ such that for any $\eta_0\in\mathcal{M}_s$, the solution of $(1.2)$ with initial data $\eta(0)=\eta_0$ possesses
  the following property:
$$\displaystyle\lim_{t\to\infty}\eta(t)=\eta^*.$$
  Conversely, if a solution of $(1.2)$ lying in $\mathcal{O}$ satisfies this property then its initial data $\eta_0\in
  \mathcal{M}_s$.

  $(4)$\ \ For any $\eta_0\in\mathcal{O}$ there exist unique $a\in G$ and $\xi_0\in\mathcal{M}_s$ such that
  $\eta_0=g(a,\xi_0)$ and for the solution $\eta=\eta(t)$ of $(1.2)$ with initial data $\eta(0)=\eta_0$,}
$$
   \lim_{t\to\infty}\eta(t)=g(a,\eta^*).
$$

  $\mathcal{M}_c$ is called {\em center manifold} of the equation (1.2), and $\mathcal{M}_s$ is called the {\em stable manifold}
  of (1.2) corresponding to the stationary point $\eta^*$. The above theorem shows that the center manifold $\mathcal{M}_c$ of
  (1.2) is exactly the combined trajectory of the Lie group actions $(G_i,p_i)$, $i=1,2,\cdots,N$, to a stationary point
  $\eta^*$, and the phase diagram of the equation (1.2) in the neighborhood $\mathcal{O}$ of $\mathcal{M}_c$ has a nice structure:
  All flows in $\mathcal{O}$ converges to the center manifold $\mathcal{M}_c$ in such a way that by collecting flows converging
  to the stationary point $\eta^*$ together to make up the stable manifold $\mathcal{M}_s$ corresponding to this point, all
  other flows are images of flows in $\mathcal{M}_s$ under the combined actions of Lie groups $(G_i,p_i)$, $i=1,2,\cdots,N$. It
  follows that $\mathcal{O}$ is a homogeneous fibre bundle over the center manifold $\mathcal{M}_c$ when regarding stable manifolds
  as fibres, with each fibre being the image of one fibre under the combined actions of Lie groups $(G_i,p_i)$, $i=1,2,\cdots,N$;
  see Figure 1.

\vspace*{2cm}
\centerline{Put figure1 here $\qquad\qquad\qquad\qquad$ Put figure2 here}
\vspace*{2cm}
\centerline{Figure 1 $\qquad\qquad\qquad\qquad\qquad$ Figure 2}
%\begin{figure}
%\begin{center}
%\includegraphics[width=10.5cm,height=6.5cm]{figure1.jpg}
%\end{center}
%\end{figure}

  Phase diagram depicted in Theorems 1.1 can be regarded as a generalization of the phase diagram of the following planer system:
$$
  x'=0, \quad y'=-y.
$$
  Let $\mathfrak{M}=\mathfrak{M}_0=X=X_0={\mathbb{R}}^2$, and define $F:{\mathbb{R}}^2\to T({\mathbb{R}}^2)={\mathbb{R}}^2$
  to be the 2-vector function $F(x,y)=(0,-y)$, $\forall (x,y)\in{\mathbb{R}}^2$. Then the above system can be rewritten as
  $u'=F(u)$ with $u=(x,y)$. Let $G={\mathbb{R}}$ be the usual one-dimensional additive Lie group. We introduce an action $p$ of
  $G$ to ${\mathbb{R}}^2$ as follows:
$$
  p(z,(x,y))=(x+z,y), \quad \forall (x,y)\in{\mathbb{R}}^2, \;\; \forall z\in G.
$$
  Since $\partial_up(z,u)=id$ and $F(p(z,u))=F(u)$, $\forall u\in{\mathbb{R}}^2$, $\forall z\in G$, we see this equation
  is invariant under the group action $(G,p)$. Hence Theorem 1.1 (with $N=1$) applies to it. The phase diagram of this equation
  is as in Figure 2. Figure 1 is clearly a sophistic generalization of Figure 2.
\medskip

  {\em Remark}.\ \ Recall that for a Banach space $X$, a mapping $F:X\to X$ and a linear group action $(G,\cdot)$ on $X$, i.e.,
  $\cdot$ is a continuous group homorphism from the Lie group $G$ to the group $L_{is}(X)$ of all self-isomorphisms of
  $X$, such that for any $a\in G$, $[x\mapsto a\cdot x]\in L_{is}(X)$, the mapping $F$ is called {\em equivariant} with respect
  to this Lie group action if it satisfies the following condition (cf., e.g., \cite{EscPro}):
$$
  F(a\cdot x)=a\cdot F(x), \quad \forall a\in G, \;\; \forall x\in X.
$$
  Linearity of this action implies that if we redenote it as $p$, i.e., $p(a,x)=a\cdot x$ for $(a,x)\in G\times X$, then
  $D_xp(a,x)y=a\cdot y$ for all $a\in G$ and $x,y\in X$. Hence, the notion of equivariance in existing literatures such as
  \cite{EscPro} is a special situation of the notion of invariance here, or in another word, the concept
  of invariance defined here is an extension of the concept of equivariance for linear Lie group actions to general possibly
  nonlinear Lie group actions. Note that in \cite{Cui3}, the phrase ``quasi-invariance'' rather than ``invariance'' as here
  is used.

  As an application of Theorem 1.1, let us consider the following two-free-surface Hele-Shaw problem::
\begin{equation}
\left\{
\begin{array}{ll}
   \Delta u(x,t)=0, &\quad x\in D(t), \;\; t>0,\\
   u(x,t)=\gamma\kappa(x,t),   &\quad x\in S(t), \;\; t>0,\\
   u(x,t)=\mu\kappa(x,t),   &\quad x\in\Gamma(t), \;\; t>0,\\
   V_n(x,t)=-\partial_{\mbox{\footnotesize $\bfn$}}u(x,t), &\quad x\in S(t), \;\; t>0,\\
   \partial_{\mbox{\footnotesize $\bfn$}}u(x,t)=\bfc\cdot\bfn,   &\quad x\in\Gamma(t), \;\; t>0,\\
   \displaystyle\;\frac{1}{|S(t)|}\oint_{S(t)}x{\rm d}\sigma=\!\!\!&\!\!\!\displaystyle\frac{1}{|\Gamma(t)|}\oint_{\Gamma(t)}x{\rm d}\sigma,
   \quad t>0,\\
   S(0)=S_0. &
\end{array}
\right.
\end{equation}
%---(1.7)---
  Here $D(t)$ is an unknown annular domain in ${\mathbb{R}}^n$ varying in time $t$, with an outer surface $S(t)$ and
  an inner surface $\Gamma(t)$, $u=u(x,t)$ is an unknown function defined for $x\in\overline{D(t)}$ and $t\geqslant 0$,
  $\kappa(\cdot,t)$ is the mean curvature of $\partial D(t)=S(t)\cup\Gamma(t)$, $V_n$ is the normal velocity of the outer
  surface $S(t)$, $\bfn$ denotes the outward-pointing normal field of $\partial D(t)$, $\gamma,\mu$ are positive constants,
  $\gamma>\mu$, $\bfc$ is an unknown constant real $n$-vector to be determined together with $D(t)$ and $u(x,t)$, and $S_0$
  is a given initial hypersurface. We note that $\Delta$ is the Laplacian in $n$ variables, and $|S(t)|$,
  $|\Gamma(t)|$ denote surface measures of the surfaces $S(t)$ and $\Gamma(t)$, respectively. We take the convention that
  for a convex closed surface its mean curvature takes nonnegative values. For simplicity we only consider the case that
  $S_0$ is a $\dot{C}^{m+\mu}$-hypersurface homeomorphic and sufficiently close to a sphere with a given radius $R$ centered
  at the origin.
  The condition $(1.7)_6$ is imposed to ensure the surfaces $S(t)$ and $\Gamma(t)$ have a common center; it is imposed to
  ensure uniqueness of the solution. Indeed, if this condition is removed, then for a given sphere surface $S_0=\partial B(0,R)$
  ($R>0$), the above problem has infinitely many solutions: For an arbitrary closed sphere $\overline{B(x_0,K)}\subseteq B(0,R)$
  with radius $K>0$ satisfying the following condition
\begin{equation}
   \frac{\gamma}{R}=\frac{\mu}{K},
\end{equation}
%---(1.8)---
  by putting $S(t)\equiv\partial B(0,R)$, $\Gamma(t)\equiv\partial B(x_0,K)$
  and $u(x,t)\equiv{\gamma}/{R}$ for $x\in \overline{B(0,R)}\backslash B(x_0,K)$ and $t\geqslant 0$, we see that $(u,S,\Gamma)$
  is a solution if we fix $\bfc=0$. Unknown $n$-vector $\bfc$ is introduced to balance the condition $(1.7)_6$, because without
  this unknown $n$-vector (1.7) is an overdetermined system. Physical explanation of this condition is as follows: We know
  that in physics $u$ represents pressure of a fluid in a porous medium. Let $\bfv$ be velocity field of the fluid. The equation
  $(1.7)_1$ is a composition of the following two equations:
$$
  \mbox{\rm div}\,\bfv=0, \qquad \bfv=-\nabla u.
$$
  The first equation means that the fluid is incompressible, and the second equation is the Darcy law. The equation $(1.7)_5$
  can be rewritten as follows:
$$
  \bfn\cdot(\bfv+\bfc)=0 \quad \mbox{on}\;\; \Gamma(t),
$$
  which means that on the inner surface the normal velocity of the fluid might not be vanishing; but by shifting the velocity
  with a constant vector, the normal component vanishes.

  The problem (1.7) is an extension of the classical Hele-Shew problem which has been extensively studied by many authors
  during the past fifty years, cf. \cite{Chen, EscSim, FriRei} and the references cited therein. Note that classical Hele-Shew
  problem has only one-free boundary, whereas the above problem has two free-surfaces: Both $S(t)$ and $\Gamma(t)$ need to be
  determined together with the solution $u(x,t)$ of the partial differential equation. Similar to the classical Hele-Shew
  problem, the above problem has also infinitely many stationary solutions depending on $n\!+\!1$ parameters, given by
$$
  u(x,t)\equiv u_s(x)=\frac{\gamma}{R}, \quad S(t)\equiv S_s=\partial B(x^*,R), \quad \Gamma(t)\equiv\Gamma_s=\partial B(x^*,K),
$$
  where $x^*$ is an arbitrary point in $\mathbb{R}^n$, $R$ is an arbitrary positive number, and $K$ is the number determined
  by the relation (1.8). By using Theorem 1.1, we shall prove that this stationary solution is asymptotically stable module
  translation and scaling, i.e., for each stationary solution $(u_s,S_s,\Gamma_s)$ there exists corresponding $\delta>0$ such
  that for any initial hypersurface $S_0$ in the $\delta$-neighborhood of $S_s$, the solution of (1.7) exists globally and, as
  $t\to\infty$, converges to a stationary solution $(u_s',S_s',\Gamma_s')$ which may not be $(u_s,S_s,\Gamma_s)$ but is obtained
  from it through translation and scaling. To give a precise statement of this result, let us first introduce some basic concepts
  and notations.

  Given a nonnegative integer $m$, a real number $\mu\in [0,1]$, a bounded open set $\Omega\subseteq{\mathbb{R}}^n$
  and a sufficiently smooth (e.g., smooth up to order $k$ for some integer $k\geqslant m+\mu$) closed hypersurface
  $S\subseteq{\mathbb{R}}^n$, the notations $C^{m+\mu}(\overline{\Omega})$ and $C^{m+\mu}(S)$ denote the usual $m\!+\!\mu$-th
  order H\"{o}lder spaces on $\overline{\Omega}$ and $S$, respectively, and the notation $C^{m+\mu}(\overline{\Omega},{\mathbb{R}}^n)$
  denotes the usual $m\!+\!\mu$-th order $n$-vector H\"{o}lder space on $\overline{\Omega}$, with the cases $\mu=0,1$
  understood in conventional sense. We use the notation $\dot{C}^{m+\mu}(\overline{\Omega})$ to denote the closure of
  $C^{\infty}(\overline{\Omega})$ in $C^{m+\mu}(\overline{\Omega})$, and similarly for the notations $\dot{C}^{m+\mu}(S)$
  and $\dot{C}^{m+\mu}(\overline{\Omega},{\mathbb{R}}^n)$. The last three spaces are called $m\!+\!\mu$-th order little
  H\"{o}lder spaces. A significant difference between little H\"{o}lder spaces and H\"{o}lder spaces is that for
  nonnegative integers $k,m$ and real numbers $\mu,\nu\in [0,1]$, if $k+\nu>m+\mu$ then $\dot{C}^{k+\nu}(\overline{\Omega})$
  (resp. $\dot{C}^{k+\nu}(S)$, $\dot{C}^{k+\nu}(\overline{\Omega},{\mathbb{R}}^n)$) is dense in
  $\dot{C}^{m+\mu}(\overline{\Omega})$ (resp. $\dot{C}^{m+\mu}(S)$, $\dot{C}^{m+\mu}(\overline{\Omega},{\mathbb{R}}^n)$),
  but $C^{k+\nu}(\overline{\Omega})$ (resp. $C^{k+\nu}(S)$, $C^{k+\nu}(\overline{\Omega},{\mathbb{R}}^n)$) is not dense
  in $C^{m+\mu}(\overline{\Omega})$ (resp. $C^{m+\mu}(S)$, $C^{m+\mu}(\overline{\Omega},{\mathbb{R}}^n)$).

  Recall (cf. \cite{Cui3}) that an open set $\Omega\subseteq{\mathbb{R}}^n$ is said to be a simple $C^{m+\mu}$-domain if
  $\Omega$ is $C^{m+\mu}$-diffeomorphic to the open unit sphere $B(0,1)$ in ${\mathbb{R}}^n$, i.e., there exists a bijective
  mapping $\Phi:\overline{B(0,1)}\to\overline{\Omega}$ satisfying the following properties:
$$
  \Phi\in C^{m+\mu}(\overline{B(0,1)},{\mathbb{R}}^n) \quad \mbox{and} \quad \Phi^{-1}\in C^{m+\mu}(\overline{\Omega},{\mathbb{R}}^n).
$$
  We use the notations $\mathfrak{D}^{m+\mu}({\mathbb{R}}^n)$ and $\mathfrak{S}^{m+\mu}({\mathbb{R}}^n)$ to denote the sets
  of all simple $C^{m+\mu}$-domains and their boundaries in ${\mathbb{R}}^n$, respectively. If instead of $C^{m+\mu}$ the
  notation $\dot{C}^{m+\mu}$ is used in the above relations, then the notations $\dot{\mathfrak{D}}^{m+\mu}({\mathbb{R}}^n)$
  and $\dot{\mathfrak{S}}^{m+\mu}({\mathbb{R}}^n)$ are used correspondingly. From the discussion in \cite{Cui3} we know that
  $\mathfrak{D}^{m+\mu}({\mathbb{R}}^n)$ and $\dot{\mathfrak{D}}^{m+\mu}({\mathbb{R}}^n)$ are Banach manifolds built on the
  Banach spaces $C^{m+\mu}(\mathbb{S}^{n-1})$ and $\dot{C}^{m+\mu}(\mathbb{S}^{n-1})$, respectively, and each point $\Omega
  \in\mathfrak{D}^{m+1+\mu}({\mathbb{R}}^n)$ (resp. $\dot{\mathfrak{D}}^{m+1+\mu}({\mathbb{R}}^n)$) is a differentiable point
  in $\mathfrak{D}^{m+\mu}({\mathbb{R}}^n)$ (resp. $\dot{\mathfrak{D}}^{m+\mu}({\mathbb{R}}^n)$). Since
  $\mathfrak{S}^{m+\mu}({\mathbb{R}}^n)$ (resp. $\dot{\mathfrak{S}}^{m+\mu}({\mathbb{R}}^n)$) can be identified with
  $\mathfrak{D}^{m+\mu}({\mathbb{R}}^n)$ (resp. $\dot{\mathfrak{D}}^{m+\mu}({\mathbb{R}}^n)$), it follows that
  $\mathfrak{S}^{m+\mu}({\mathbb{R}}^n)$ and $\dot{\mathfrak{S}}^{m+\mu}({\mathbb{R}}^n)$ are also Banach manifolds built on the
  Banach spaces $C^{m+\mu}(\mathbb{S}^{n-1})$ and $\dot{C}^{m+\mu}(\mathbb{S}^{n-1})$, respectively, and all the results
  obtained in \cite{Cui3} for $\mathfrak{D}^{m+\mu}({\mathbb{R}}^n)$ (resp. $\dot{\mathfrak{D}}^{m+\mu}({\mathbb{R}}^n)$)
  work for $\mathfrak{S}^{m+\mu}({\mathbb{R}}^n)$ (resp. $\dot{\mathfrak{S}}^{m+\mu}({\mathbb{R}}^n)$).

  Let $m$ be a positive integer $\geqslant 2$ and $0<\mu<1$. Let $\mathfrak{M}:=\dot{\mathfrak{S}}^{m+\mu}({\mathbb{R}}^n)$
  and $\mathfrak{M}_0:=\dot{\mathfrak{S}}^{m+3+\mu}({\mathbb{R}}^n)$. We know that $\mathfrak{M}$ and $\mathfrak{M}_0$ are
  Banach manifolds built on the Banach spaces $\dot{C}^{m+\mu}({\mathbb{S}}^{n-1})$ and $\dot{C}^{m+3+\mu}({\mathbb{S}}^{n-1})$,
  respectively, and $\dot{C}^{m+3+\mu}({\mathbb{S}}^{n-1})$ is dense in $\dot{C}^{m+\mu}({\mathbb{S}}^{n-1})$. We introduce a
  vector field $\mathscr{F}$ in $\mathfrak{M}$ with domain $\mathfrak{M}_0$ in the following way. Given $S\in\mathfrak{M}_0$,
  let $\Omega\in\dot{\mathfrak{D}}^{m+3+\mu}({\mathbb{R}}^n)$ be the simple domain such that $\partial\Omega=S$. Consider the
  following problem: Find a hypersurface $\Gamma\subseteq\Omega$, a real $n$-vector $\bfc$ and a function $u$ defined in the
  domain $D$ enclosed by $S$ and $\Gamma$ (i.e., $\partial D=S\cup\Gamma$) such that the following equations are satisfied:
\begin{equation}
\left\{
\begin{array}{ll}
   \Delta u(x)=0, &\quad x\in D,\\
   u(x)=\gamma\kappa(x),   &\quad x\in S,\\
   u(x)=\mu\kappa(x),   &\quad x\in\Gamma,\\
   \partial_{\mbox{\footnotesize$\bfn$}}u(x)=\bfc\cdot\bfn,   &\quad x\in\Gamma,\\
   \displaystyle\frac{1}{|S|}\oint_{S}x{\rm d}\sigma
   =&\displaystyle\frac{1}{|\Gamma|}\oint_{\Gamma}x{\rm d}\sigma.
\end{array}
\right.
\end{equation}
%---(1.9)---
  We shall prove this problem has an unique solution; see Section 5. We now define
\begin{equation}
   \mathscr{F}(S)=-\partial_{\mbox{\footnotesize$\bfn$}}u|_{S}\in\dot{C}^{m+\mu}(S)
   =\mathcal{T}_{S}(\mathfrak{M}).
\end{equation}
%---(1.10)---
  It follows that the problem (1.7) reduces into the following initial value problem of a differential equation in the
  Banach manifold $\mathfrak{M}$:
\begin{equation}
\left\{
\begin{array}{ll}
   S'(t)=\mathscr{F}(S(t)), &\quad  t>0,\\
   S(0)=S_0. &
\end{array}
\right.
\end{equation}
%---(1.11)---
  See Section 5 for details.

  Let $G_{tl}={\mathbb{R}}^n$ be the additive group of $n$-vectors. Given $z\in{\mathbb{R}}^n$ and $S\in
  \dot{\mathfrak{S}}^{m+\mu}({\mathbb{R}}^n)$, let
$$
  p(z,S)=S+z=\{x+z:\,x\in S\}.
$$
  It is clear that $p(z,S)\in\dot{\mathfrak{S}}^{m+\mu}({\mathbb{R}}^n)$, $\forall S\in\dot{\mathfrak{S}}^{m+\mu}({\mathbb{R}}^n)$,
  $\forall z\in{\mathbb{R}}^n$. It can be easily seen that $(G_{tl},p)$ is a Lie group action on
  $\dot{\mathfrak{S}}^{m+\mu}({\mathbb{R}}^n)$. By Lemma 4.1 of \cite{Cui3} we know that the action $p(z,S)$ is differentiable
  at every point $S\in\dot{\mathfrak{S}}^{m+1+\mu}({\mathbb{R}}^n)$, and $\mbox{rank}D_{z}p(z,S)=n$, $\forall z\in G_{tl}$,
  $\forall S\in\dot{\mathfrak{D}}^{m+1+\mu}({\mathbb{R}}^n)$.

  Next let $G_{\!dl}={\mathbb{R}}_+=(0,\infty)$ be the multiplicative group of all positive numbers. Given $\lambda\in G_{\!dl}$
  and $S\in\dot{\mathfrak{S}}^{m+\mu}({\mathbb{R}}^n)$, let
$$
  q(\lambda,S)=\lambda S=\{\lambda x:\,x\in S\}.
$$
  Clearly $q(\lambda,S)\in\dot{\mathfrak{S}}^{m+\mu}({\mathbb{R}}^n)$, $\forall S\in\dot{\mathfrak{S}}^{m+\mu}({\mathbb{R}}^n)$,
  $\forall\lambda\in G_{\!dl}$, and $(G_{\!dl},q)$ is also a Lie group action on $\dot{\mathfrak{S}}^{m+\mu}({\mathbb{R}}^n)$.
  By Lemma 4.2 of \cite{Cui3} we know that the action $q(\lambda,\Omega)$ is differentiable at every point $S\in
  \dot{\mathfrak{S}}^{m+1+\mu}({\mathbb{R}}^n)$, and $\mbox{rank}D_{z}q(\lambda,S)=1$, $\forall\lambda\in G_{dl}$,
  $\forall S\in\dot{\mathfrak{S}}^{m+1+\mu}({\mathbb{R}}^n)$).

  The group actions $(G_{tl},p)$ and $(G_{\!dl},q)$ to $\dot{\mathfrak{S}}^{m+\mu}({\mathbb{R}}^n)$ are not mutually commutative.
  However, they are quasi-commutative in the sense that they satisfy the following relation:
$$
  q(\lambda,p(z,S))=p(\lambda z, q(\lambda,S)), \quad \forall S\in\dot{\mathfrak{S}}^{m+\mu}({\mathbb{R}}^n),\;\;
  \forall z\in G_{tl},\;\; \forall\lambda\in G_{\!dl}.
$$
  Besides, denoting
$$
  g(z,\lambda,S)=p(z, q(\lambda,S)), \quad \forall S\in\dot{\mathfrak{S}}^{m+\mu}({\mathbb{R}}^n),\;\;
  \forall z\in G_{tl},\;\; \forall\lambda\in G_{\!dl},
$$
  we easily see that the following relation holds:
$$
  {\rm rank}\,\partial_{(z,\lambda)}g(z,\lambda,S)=n+1, \quad
  \forall S\in\dot{\mathfrak{D}}^{m+1+\mu}({\mathbb{R}}^n),\;\;
  \forall z\in G_{tl},\;\; \forall\lambda\in G_{\!dl}.
$$
  It can be easily shown that the vector field $\mathscr{F}$ is invariant under the translation group action $(G_{tl},p)$, and
  quasi-invariant under the dilation group action $(G_{dl},q)$ with quasi-invariant factor $\theta(\lambda)=\lambda^{-3}$,
  $\lambda>0$, i.e., the following relations hold:
\begin{equation}
  \mathscr{F}(p(z,S))=\partial_{S}p(z,S)\mathscr{F}(S), \quad \forall z\in G_{tl},\;\; \forall S\in\mathfrak{M}_0,
\end{equation}
%---(5.15)---
\begin{equation}
  \mathscr{F}(q(\lambda,S))=\lambda^{-3}\partial_{S}q(\lambda,S)\mathscr{F}(S), \quad \forall\lambda\in G_{dl},\;\;
  \forall S\in\mathfrak{M}_0;
\end{equation}
%---(5.16)---
  see Section 5 for details. Based on these facts and an analysis of the spectrum of the derivative of $\mathscr{F}$ at its
  stationary point, it follows by applying Theorem 1.1 that the following result holds:
\medskip

  {\bf Theorem 1.2}\ \ {\em Let $\mathcal{M}_c$ be the $(n\!+\!1)$-dimensional submanifold of $\mathfrak{M}_0$ consisting of
  all surface spheres in ${\mathbb{R}}^n$. We have the following assertions:

  $(1)$\ \ There is a neighborhood $\mathcal{O}$ of $\mathcal{M}_c$ in $\mathfrak{M}_0$ such that for any $S_0\in\mathcal{O}$,
  the initial value problem $(1.11)$ has a unique solution $S\in C([0,\infty),\mathfrak{M}_0)\cap C^1((0,\infty),\mathfrak{M}_0)$.
  Correspondingly, the free-boundary problem $(1.7)$ has a unique solution $(S,\Gamma,u)$ with the property that $S,\Gamma\in
  C([0,\infty),\mathfrak{M}_0)\cap C^1((0,\infty),\mathfrak{M}_0)$.

  $(2)$\ \ There exists a submanifold $\mathcal{M}_s$ of $\mathfrak{M}_0$ of codimension $n\!+\!1$ passing $S_s=\partial B(0,1)$
  such that for any $S_0\in\mathcal{M}_s$, the solution of the problem $(1.11)$ satisfies $\displaystyle\lim_{t\to\infty}
  S(t)=S_s$ and, conversely, if the solution of $(1.11)$ satisfies this property then $S_0\in\mathcal{M}_s$.

  $(3)$\ \ For any $S_0\in\mathcal{O}$ there exist unique $x_0\in {\mathbb{R}}^n$, $R>0$ and $T_0\in\mathcal{M}_s$ such that
  $S_0=x_0+RT_0$ and, for the solution $S=S(t)$ of $(1.11)$, we have
$$
   \lim_{t\to\infty}S(t)=\partial B(x_0,R).
$$
  Correspondingly, the second component of the solution $(S,\Gamma,u)$ of the free-boundary problem $(1.7)$ has the following
  property:
$$
   \lim_{t\to\infty}\Gamma(t)=\partial B(x_0,K),
$$
  where $K$ is the number determined by the relation $(1.8)$. Moreover, convergence rate of the above limit relations is of the
  form $C\mbox{\rm e}^{-\nu t}$ for some positive constants $C$ and $\nu$ depending on $\gamma$ and $R$.}
\medskip

  {\em Remark}.\ \ In the above theorem properties of $u(x,t)$ are not stated. This is because statement of such properties
  is very complex. Actually, having known properties of $S(t)$ and $\Gamma(t)$, properties of $u(x,t)$ easily follow from
  well-known theory of elliptic boundary value problems.
\medskip

  Organization of the rest part is as follows. In Section 2 we introduce some basic concepts concerning
  quasi-differential structure for a class of nondifferentiable Banach manifolds. In Section 3 we discuss some basic
  concepts and results concerning differential equations in quasi-differential Banach manifolds. In Section 4 we give the
  proof of Theorem 1.1. In the last section we give the proof of Theorem 1.2.

\section{Basic concepts on differential calculus in quasi-differentiable Banach manifold}
\setcounter{equation}{0}

\hskip 2em
  In this section we introduce some basic concepts concerning differential calculus in a class of nondifferentiable Banach
  manifolds. We shall show that for this class of Banach manifolds, it is possible to introduce the concept of differentiable
  point and tangent space at a such point. This enables us to partially extend the technique of differential calculus for
  traditional differentiable Banach manifolds to this class of nondifferentiable Banach manifolds and, consequently,
  differential equations in such nondifferentiable Banach manifolds can be considered. All the discussion made in this section
  is modeled by and toward to application to the Banach manifold $\dot{\mathfrak{D}}^{m+\mu}(\mathbb{R}^n)$ ($m\in\mathbb{N}$,
  $0\leqslant\mu\leqslant 1$) of simple $\dot{C}^{m+\mu}$-domains in $\mathbb{R}^n$ studied in the reference \cite{Cui3}, and
  the Banach manifold $\dot{\mathfrak{S}}^{m+\mu}(\mathbb{R}^n)$ ($m\in\mathbb{N}$, $0\leqslant\mu\leqslant 1$), which
  are merely topological Banach manifolds and do not possess differentiable structure in traditional sense, but in the study
  of many evolutionary free boundary problems we often need to consider differential equations in these Banach manifolds or their
  tangent bundles.

  Let $X$ and $X_0$ be two Banach spaces. We say $X_0$ is an {\em embedded Banach subspace} of $X$ if as linear spaces
  $X_0$ is a subspace of $X$, and the restriction of the norm of $X$ to $X_0$ is majorized by that of $X_0$, i.e., there exists
  a constant $C>0$ such that
\begin{equation}
  \|x\|_{X}\leqslant C\|x\|_{X_0}, \qquad \forall x\in X_0.
\end{equation}
%---(2.1)---
  We say $X_0$ is {\em densely embedded} in $X$ if $X_0$ is dense in $X$.

  Let $X$ and $Y$ be two Banach spaces. Recall that for a positive integer $k$, the notation $L^k(X,Y)$ denotes the Banach
  space of all bounded $k$-linear mappings from $X\times X\times\cdots\times X$ ($k$ times) to $Y$. Recall that a $k$-linear
  mapping $A$ from $X\times X\times\cdots\times X$ ($k$ times) to $Y$ is said to be {\em bounded} if there exists a positive
  constant $C$ such that the following relation holds for all $(x_1,x_2,\cdots,x_k)\in\underbrace{X\times X\times\cdots\times X}_{k}$:
$$
  \|A(x_1,x_2,\cdots,x_k)\|_Y\leqslant C\|x_1\|_X\|x_2\|_X\cdots,\|x_k\|_X.
$$
  The infimum of all such constant $C$ is called the norm of $A$ and is denoted as $\|A\|_{L^k(X,Y)}$. Note that
$$
  L^k(X,Y)\approxeq L(X,L(X,\cdots,L(X,Y)\cdots))\;\; \mbox{($k$-times $L$ and $X$)}.
$$
  Recall that an element $A\in L^k(X,Y)$ is said to be {\em symmetric} if for any permutation $i_1,i_2,\cdots,i_k$ of
  $1,2,\cdots,k$ there holds
$$
  A(x_1,x_2,\cdots,x_k)=A(x_{i_1},x_{i_2},\cdots,x_{i_k}), \quad \forall (x_1,x_2,\cdots,x_k)\in
  \underbrace{X\times X\times\cdots\times X}_{k}.
$$
  The notation $L_s^k(X,Y)$ denotes the Banach subspace of $L^k(X,Y)$ consisting of all bounded symmetric $k$-linear mappings
  from $X\times X\times\cdots\times X$ ($k$ times) to $Y$.

  Let $X$ and $X_0$ be two Banach spaces such that $X_0$ is a densely embedded Banach subspace of $X$. Let $U_0$ be an open subset
  of $X_0$. Let $Y$ be another Banach space. For a map $F:U_0\to Y$ and a point $x_0\in U_0$, we say $x_0$ is a {\em $d^k$-point}
  of $F$, where $k$ is a positive integer, if for each $1\leqslant j\leqslant k$ there exists a corresponding operator $A_j
  \in L_s^j(X,Y)$ such that the following relation holds:
$$
  \lim_{\|x-x_0\|_{X_0}\to 0}\frac{\Big\|F(x)-F(x_0)
  -\displaystyle\sum_{j=1}^k\frac{1}{j!}A_j(x-x_0,x-x_0,\cdots,x-x_0)\Big\|_{Y}}{\|x-x_0\|_{X_0}^k}=0.
$$
  $A_j$ is called the $j$-th order differential or $j$-th order Fr\'{e}chet derivative of $F$ at $x_0$ and is denoted as
  $D^jF(x_0)=A_j$ or $F^{(j)}(x_0)=A_j$, $j=1,2,\cdots,k$. In particular, for $j=1,2$ the operators $F^{(1)}(x_0)$ and
  $F^{(2)}(x_0)$ are also denoted as $F'(x_0)$ and $F''(x_0)$, respectively. Note that since $L^j(X,Y)\subseteq
  L^j(X_0,Y)$, it follows that $A_j\in L_s^j(X,Y)$ implies $A_j\in L_s^j(X_0,Y)$. From this fact it can be easily seen
  that if $x_0$ is a $d^k$-point of $F$ then it is also a $d^j$-point of $F$ for any $1\leqslant j\leqslant k\!-\!1$.
\medskip

  {\em Remark}.\ Note that if $X_0$ is not dense in $X$, then the operators $A_1$, $A_2$, $\cdots$, $A_k$ might not be uniquely
  determined by $F$ and $x_0$. The reason is that the above relation does not use values of these operators outside $X_0$,
  so that it is possible to change values of them in any subspace of $X$ which is complementary to $\bar{X}_0$ without
  changing the above relation if $X_0$ is not dense in $X$. If, however, $X_0$ is dense in $X$, then clearly $A_1$, $A_2$,
  $\cdots$, $A_k$ are uniquely determined by $F$ and $x_0$.
\medskip

  Let $X,X_0,Y$ and $U_0$ be as above. Given a positive integer $k$, we use the notation $\mathfrak{C}^k(U_0;X,Y)$ to
  denote the set of all mappings $F:U_0\to Y$ satisfing the following two conditions:
\begin{description}
\item[] (1)\ \ All points in $U_0$ are $d^k$-points of $F$;
\item[] (2)\ \ $[x\mapsto F^{(j)}(x)]\in C(U_0,L^j(X,Y))$, $j=1,2,\cdots,k$, where $U_0$ uses the topology of $X_0$.
\end{description}
  It is clear that $\mathfrak{C}^k(U_0;X,Y)\subseteq C^k(U_0,Y)$, where $C^k(U_0,Y)$ denotes the set of all $k$-th order
  continuously differentiable mappings $F:U_0\subseteq X_0\to Y$, and an element $F\in C^k(U_0,Y)$ belongs to
  $\mathfrak{C}^k(U_0;X,Y)$ if and only if for each $1\leqslant j\leqslant k$ and any $x\in U_0$, $D^jF(x)\in L_s^j(X_0,Y)$\footnotemark
\footnotetext[1]{Here $D^jF(x)$ denotes the $j$-th order Fr\'{e}chet derivative of $F$ at $x$ in usual sense, which is
inductively defined to be the Fr\'{e}chet derivative of $D^{j-1}F$ at $x$, and $D^0F=F$.}
   can be extended into an operator belonging to $L_s^j(X,Y)$, and after extension $D^jF\in C(U_0,L^j(X,Y))$, where $U_0$ uses
   the topology of $X_0$. Hence the condition $F\in\mathfrak{C}^k(U_0;X,Y)$ is stronger than the condition $F\in C^k(U_0,Y)$.

\medskip
  {\em Example}\ \ As before for $m\in\mathbb{N}$($=$the set of all positive integers) and $0\leqslant\mu\leqslant 1$ we use
  the notation $\dot{C}^{m+\mu}$ to denote $m+\mu$-th order little H\"{o}lder functions or spaces (for $\mu=1$, $\dot{C}^{m+\mu}$
  refers to $\dot{C}^{m+1-0}$). Let $X=\dot{C}^{m+\mu}(\mathbb{S}^{n-1})$ and $X_0=\dot{C}^{m+k+\mu}(\mathbb{S}^{n-1})$, where
  $k\in\mathbb{N}$. Let $f\in C^{\infty}(\mathbb{S}^{n-1}\times\mathbb{R},\mathbb{S}^{n-1})$. Then the map $F:\rho\mapsto
  [x\mapsto\rho(f(x,\rho(x)))]$ belongs to $\mathfrak{C}^k(X_0;X,X)$.

\medskip
  {\bf Definition 2.1}\ \ {\em Let $\mathfrak{M}$ and $\mathfrak{M}_0$ be two topological $($i.e., they need not have
  differentiable structure$)$ Banach manifolds built on Banach spaces $X$ and $X_0$, respectively. Let $\mathscr{A}$ be
  a family of local charts of $\mathfrak{M}$. Let $k$ be a positive integer. We say $\mathfrak{M}_0$ is a \mbox{\small\boldmath
  $C^k$-$embedded\;Banach\; submanifold$} of $\mathfrak{M}$ with respect to $\mathscr{A}$ if the following conditions are satisfied:
\begin{enumerate}
\item[]$(D1)$\ \ $X_0$ is a densely embedded Banach subspace of $X$.\vspace*{-0.2cm}
\item[]$(D2)$\ \ $\mathfrak{M}_0$ is an embedded topological subspace of $\mathfrak{M}$, i.e., $\mathfrak{M}_0\subseteq
  \mathfrak{M}$, and for any open subset $U$ of $\mathfrak{M}$, $U\cap\mathfrak{M}_0$ is an open subset of $\mathfrak{M}_0$.
\vspace*{-0.2cm}
\item[]$(D3)$\ \ For any $\eta\in\mathfrak{M}$ there exists a local chart $(\mathcal{U},\varphi)\in\mathscr{A}$ such
  that $\eta\in\mathcal{U}$.
\item[]$(D4)$\ \ For any $\eta\in\mathfrak{M}_0$ and any local chart $(\mathcal{U},\varphi)\in\mathscr{A}$ such that
  $\eta\in\mathcal{U}$, by letting $\mathcal{U}_0=\mathcal{U}\cap\mathfrak{M}_0$, $(\mathcal{U}_0,\varphi|_{\mathcal{U}_0})$
  is a local chart of $\mathfrak{M}_0$ at $\eta$.  \vspace*{-0.2cm}
\item[]$(D5)$\ \ For any $\eta\in\mathfrak{M}_0$ and any $(\mathcal{U},\varphi),(\mathcal{V},\psi)\in\mathscr{A}$ such that
  $\eta\in\mathcal{U}\cap\mathcal{V}$, by letting $\mathcal{U}_0=\mathcal{U}\cap\mathfrak{M}_0$ and $\mathcal{V}_0=\mathcal{V}\cap
  \mathfrak{M}_0$, the following relations hold:
$$
  \psi\circ\varphi^{-1}\in\mathfrak{C}^k(\varphi(\mathcal{U}_0\cap\mathcal{V}_0);X,X) \quad \mbox{and} \quad
  \varphi\circ\psi^{-1}\in\mathfrak{C}^k(\psi(\mathcal{U}_0\cap\mathcal{V}_0);X,X).
$$
\end{enumerate}
  We call any local chart $(\mathcal{U},\varphi)\in\mathscr{A}$ such that $\eta\in\mathcal{U}$ a \mbox{\small\boldmath
  $(C^k,\mathfrak{M}_0)$-$regular\;local\; chart$} of $\mathfrak{M}$ at $\eta$, and call the family $\mathscr{A}$ a
  \mbox{\small\boldmath $(C^k,\mathfrak{M}_0)$-$regular\; local\; chart\;system$} of $\mathfrak{M}$. For the special
  case $k=1$ we simply call $(C^k,\mathfrak{M}_0)$-regular as \mbox{\small\boldmath $\mathfrak{M}_0$-$regular$}.}
\medskip

  {\em Remark}.\ From the conditions (D1)$\sim$(D4) we easily see that $\mathfrak{M}_0$ is dense in $\mathfrak{M}$. As usual,
  we call the Banach spaces $X$ and $X_0$ the {\em base spaces} of $\mathfrak{M}$ and $\mathfrak{M}_0$, respectively.
\medskip

  {\em Example}.\ Let $m\in\mathbb{N}$ and $0\leqslant\mu\leqslant1$. As in Section 1 we use the notations
  $\dot{\mathfrak{D}}^{m+\mu}(\mathbb{R}^n)$ and $\dot{\mathfrak{S}}^{m+\mu}(\mathbb{R}^n)$ to denote the Banach manifolds
  of all simple $\dot{C}^{m+\mu}$-domains in $\mathbb{R}^n$ and their boundaries, respectively, where $\dot{C}^{m+\mu}$
  represents $m+\mu$-th order little H\"{o}lder space. Let $\mathscr{A}$ be the set of all regular local charts of
  $\dot{\mathfrak{D}}^{m+\mu}(\mathbb{R}^n)$ (recall that a local chart is called to be a regular local chart if its base
  hypersurface is smooth). Then for any integer $k\geqslant 1$, $\dot{\mathfrak{D}}^{m+k+\mu}(\mathbb{R}^n)$ is a
  $C^k$-embedded Banach submanifold of $\dot{\mathfrak{D}}^{m+\mu}(\mathbb{R}^n)$. The base spaces of
  $\dot{\mathfrak{D}}^{m+\mu}(\mathbb{R}^n)$ and $\dot{\mathfrak{D}}^{m+k+\mu}(\mathbb{R}^n)$ are respectively
  $\dot{C}^{m+\mu}(\mathbb{S}^{n-1})$ and $\dot{C}^{m+k+\mu}(\mathbb{S}^{n-1})$. Similarly, let $\mathscr{A}'$ be the set
  of all local charts of $\dot{\mathfrak{S}}^{m+\mu}(\mathbb{R}^n)$ with smooth base hypersurfaces (a such local chart is
  also called a regular local chart of $\dot{\mathfrak{S}}^{m+\mu}(\mathbb{R}^n)$). Then for any integer $k\geqslant 1$,
  $\dot{\mathfrak{S}}^{m+k+\mu}(\mathbb{R}^n)$ is a $C^k$-embedded Banach submanifold of $\dot{\mathfrak{S}}^{m+\mu}(\mathbb{R}^n)$.
  The base spaces of $\dot{\mathfrak{S}}^{m+\mu}(\mathbb{R}^n)$ and $\dot{\mathfrak{S}}^{m+k+\mu}(\mathbb{R}^n)$ are also
  $\dot{C}^{m+\mu}(\mathbb{S}^{n-1})$ and $\dot{C}^{m+k+\mu}(\mathbb{S}^{n-1})$, respectively.
\medskip

  {\bf Definition 2.2}\ \ {\em Let $\mathfrak{M}$ be a Banach manifold, $\mathscr{A}$ a family of local charts of $\mathfrak{M}$,
  and $\mathfrak{M}_0$ a $C^1$-embedded Banach submanifold of $\mathfrak{M}$ with respect to $\mathscr{A}$. We have the following
  notions:

  $(1)$\ We say $(\mathfrak{M},\mathfrak{M}_0,\mathscr{A})$ is \mbox{\small\boldmath $inward\;spreadable$} if there exists a
  Banach manifold $\mathfrak{M}_1\subseteq\mathfrak{M}_0$ such that $\mathfrak{M}_1$ is a $C^1$-embedded Banach submanifold
  of $\mathfrak{M}_0$ with respect to the restriction of $\mathscr{A}$ to $\mathfrak{M}_0$. In this case, a local chart in
  $\mathscr{A}$ is called a $(\mathfrak{M}_0,\mathfrak{M}_1)$-regular local chart.

  $(2)$\ We say $(\mathfrak{M},\mathfrak{M}_0,\mathscr{A})$ is \mbox{\small\boldmath $outward\;spreadable$} if there exists a
  Banach manifold $\widetilde{\mathfrak{M}}\supseteq\mathfrak{M}$ and a family $\widetilde{\mathscr{A}}$ of local charts of
  $\widetilde{\mathfrak{M}}$, such that $\mathscr{A}$ is the restriction of $\widetilde{\mathscr{A}}$ to $\mathfrak{M}$ and
  $\mathfrak{M}$ is a $C^1$-embedded Banach submanifold of $\widetilde{\mathfrak{M}}$ with respect to $\widetilde{\mathscr{A}}$.

  $(3)$\ If $(\mathfrak{M},\mathfrak{M}_0,\mathscr{A})$ is both inward spreadable and outward spreadable then we call the pair
  $(\mathfrak{M},\mathscr{A})$ a \mbox{\small\boldmath $quasi$-$di\!f\!ferentiable$} Banach manifold with a \mbox{\small\boldmath
  $C^1$-$kernel$} $\mathfrak{M}_0$, or simply call it a quasi-differentiable Banach manifold without mentioning the
  $C^1$-kernel $\mathfrak{M}_0$. $\mathfrak{M}_1$ is called an \mbox{\small\boldmath $inner$ $C^1$-$kernel$} of
  $(\mathfrak{M},\mathscr{A})$, and $(\widetilde{\mathfrak{M}},\widetilde{\mathscr{A}})$ a \mbox{\small\boldmath $C^1$-$shell$}
  of $(\mathfrak{M},\mathscr{A})$. Later on we shall often omit mentioning the local chart families $\mathscr{A}$ and
  $\widetilde{\mathscr{A}}$.

  $(4)$\ If $(\mathfrak{M},\mathscr{A})$ is a quasi-differentiable Banach manifold and the $C^1$-kernel $\mathfrak{M}_0$ of
  $(\mathfrak{M},\mathscr{A})$ is a $C^k$-embedded Banach submanifold of $\mathfrak{M}$, where $k\in\mathbb{N}$, we call
  $\mathfrak{M}_0$ a \mbox{\small\boldmath $C^k$-$kernel$} of $(\mathfrak{M},\mathscr{A})$. Similarly, if the inner $C^1$-kernel
  $\mathfrak{M}_1$ of $(\mathfrak{M},\mathscr{A})$ is a $C^k$-embedded Banach submanifold of $\mathfrak{M}$, where $k\in\mathbb{N}$,
  we call $\mathfrak{M}_1$ an \mbox{\small\boldmath $inner$ $C^k$-$kernel$} of $(\mathfrak{M},\mathscr{A})$.}
\medskip

  {\em Example}.\ Let $m,k\in\mathbb{N}$ and $0\leqslant\mu\leqslant1$. Then $(\dot{\mathfrak{D}}^{m+\mu}(\mathbb{R}^n),
  \dot{\mathfrak{D}}^{m+k+\mu}(\mathbb{R}^n))$ is clearly inward spreadable, with $\dot{\mathfrak{D}}^{m+k+l+\mu}(\mathbb{R}^n)$
  for any $l\in\mathbb{N}$ being an inner $C^1$-kernel. If $m\geqslant 2$ then $(\dot{\mathfrak{D}}^{m+\mu}(\mathbb{R}^n),
  \dot{\mathfrak{D}}^{m+k+\mu}(\mathbb{R}^n))$ is outward spreadable, with $\dot{\mathfrak{D}}^{m-1+\mu}(\mathbb{R}^n)$
  being a shell. It follows that if $m\geqslant 2$ then $\dot{\mathfrak{D}}^{m+\mu}(\mathbb{R}^n)$ is a quasi-differentiable
  Banach manifold, with $\dot{\mathfrak{D}}^{m+1+\mu}(\mathbb{R}^n)$ being a kernel, $\dot{\mathfrak{D}}^{m+2+\mu}(\mathbb{R}^n)$
  a inner kernel, and $\dot{\mathfrak{D}}^{m-1+\mu}(\mathbb{R}^n)$ a shell. Similarly, $(\dot{\mathfrak{S}}^{m+\mu}(\mathbb{R}^n),
  \dot{\mathfrak{S}}^{m+k+\mu}(\mathbb{R}^n))$ is also inward spreadable, with $\dot{\mathfrak{S}}^{m+k+l+\mu}(\mathbb{R}^n)$
  for any $l\in\mathbb{N}$ being an inner $C^1$-kernel, and if $m\geqslant 2$ then $(\dot{\mathfrak{S}}^{m+\mu}(\mathbb{R}^n),
  \dot{\mathfrak{S}}^{m+k+\mu}(\mathbb{R}^n))$ is outward spreadable, with $\dot{\mathfrak{S}}^{m-1+\mu}(\mathbb{R}^n)$
  being a shell. Hence if $m\geqslant 2$ then $\dot{\mathfrak{S}}^{m+\mu}(\mathbb{R}^n)$ is also a quasi-differentiable
  Banach manifold, with $\dot{\mathfrak{S}}^{m+1+\mu}(\mathbb{R}^n)$ being a kernel, $\dot{\mathfrak{S}}^{m+2+\mu}(\mathbb{R}^n)$
  a inner kernel, and $\dot{\mathfrak{S}}^{m-1+\mu}(\mathbb{R}^n)$ a shell.
\medskip

  {\bf Lemma 2.3}\ \ {\em Let $\mathfrak{M}$ be a Banach manifold and $\mathfrak{M}_0$ a $C^1$-embedded Banach submanifold of
  $\mathfrak{M}$. Assume that $(\mathfrak{M},\mathfrak{M}_0)$ is inward spreadable and $\mathfrak{M}_1$ is an inner $C^1$-kernel.
  Then for any $\eta\in\mathfrak{M}_0$ and any three $(\mathfrak{M}_0,\mathfrak{M}_1)$-regular local charts $(\mathcal{U},\varphi)$,
  $(\mathcal{V},\psi)$ and $(\mathcal{X},\chi)$ of $\mathfrak{M}$ at $\eta$, the following relation holds:
\begin{equation}
  (\varphi\circ\chi^{-1})'(\chi(\eta))=(\varphi\circ\psi^{-1})'(\psi(\eta))(\psi\circ\chi^{-1})'(\chi(\eta)).
\end{equation}
%---(2.2)---
  In particular,}
\begin{equation}
  (\varphi\circ\psi^{-1})'(\psi(\eta))=(\psi\circ\varphi^{-1})'(\varphi(\eta))^{-1}.
\end{equation}
%---(2.3)---

  {\em Proof}.\ Let $X$, $X_0$ and $X_1$ be the base spaces of $\mathfrak{M}$, $\mathfrak{M}_0$ and $\mathfrak{M}_1$, respectively.
  We first prove that the relation (2.2) holds for any $\eta\in\mathfrak{M}_1$. For a such $\eta$, we denote
  $u_0=\varphi(\eta)$, $v_0=\psi(\eta)$ and $w_0=\chi(\eta)$. Then $u_0,v_0,w_0\in X_1$, and the conditions in this lemma ensure that
$$
  (\varphi\circ\psi^{-1})'(v_0),(\varphi\circ\chi^{-1})'(w_0),(\psi\circ\chi^{-1})'(w_0)\in L(X)\cap L(X_0),
$$
  and
$$
  [v\mapsto(\varphi\circ\psi^{-1})'(v)]\in C(\psi(\mathcal{U}_0\cap\mathcal{V}_0),L(X))\cap
  C(\psi(\mathcal{U}_1\cap\mathcal{V}_1),L(X)\cap L(X_0)),
$$
$$
  [w\mapsto(\varphi\circ\chi^{-1})'(w)]\in C(\chi(\mathcal{U}_0\cap\mathcal{W}_0),L(X))\cap
  C(\chi(\mathcal{U}_1\cap\mathcal{W}_1),L(X)\cap L(X_0)),
$$
$$
  [w\mapsto(\psi\circ\chi^{-1})'(w)]\in C(\chi(\mathcal{V}_0\cap\mathcal{W}_0),L(X))\cap
  C(\chi(\mathcal{V}_1\cap\mathcal{W}_1),L(X)\cap L(X_0)),
$$
  where $\mathcal{U}_0=\mathcal{U}\cap\mathfrak{M}_0$, $\mathcal{V}_0=\mathcal{V}\cap\mathfrak{M}_0$,
  $\mathcal{W}_0=\mathcal{W}\cap\mathfrak{M}_0$, $\mathcal{U}_1=\mathcal{U}\cap\mathfrak{M}_1$,
  $\mathcal{V}_1=\mathcal{U}\cap\mathfrak{M}_1$, $\mathcal{W}_1=\mathcal{W}\cap\mathfrak{M}_1$, with
  $\psi(\mathcal{U}_0\cap\mathcal{V}_0)$, $\chi(\mathcal{U}_0\cap\mathcal{V}_0)$ using the topology of $X_0$ and
  $\psi(\mathcal{U}_1\cap\mathcal{V}_1)$, $\chi(\mathcal{U}_1\cap\mathcal{V}_1)$ using the topology of $X_1$.
  For $u\in X_1$ sufficiently close to $u_0$, let $v=(\psi\circ\varphi^{-1})(u)$, $w=(\chi\circ\varphi^{-1})(u)
  =(\chi\circ\psi^{-1})(v)$. Then $v,w\in X_1$ and $v,w$ are sufficiently close to $v_0,w_0$, respectively. The conditions
  given in this lemma ensure that the following relations hold:
\begin{equation}
  u-u_0=(\varphi\circ\psi^{-1})'(v_0)(v-v_0)+o(\|v-v_0\|_{X_0}) \; (\,\mbox{in}\;X), \quad \mbox{as }\|v-v_0\|_{X_0}\to 0,
\end{equation}
%---(2.4)---
\begin{equation}
  u-u_0=(\varphi\circ\chi^{-1})'(w_0)(w-w_0)+o(\|w-w_0\|_{X_0}) \; (\,\mbox{in}\;X), \quad \mbox{as }\|w-w_0\|_{X_0}\to 0,
\end{equation}
%---(2.5)---
\begin{equation}
  v-v_0=(\psi\circ\chi^{-1})'(w_0)(w-w_0)+o(\|w-w_0\|_{X_1}) \; (\,\mbox{in}\;X_0), \quad \mbox{as }\|w-w_0\|_{X_1}\to 0.
\end{equation}
%---(2.6)---
  Substituting (2.6) into (2.4), we get
$$
  u-u_0=(\varphi\circ\psi^{-1})'(v_0)(\psi\circ\chi^{-1})'(w_0)(w-w_0)
  +o(\|w-w_0\|_{X_1}) \; (\,\mbox{in}\;X), \quad \mbox{as }\|w-w_0\|_{X_1}\to 0.
$$
  Comparing this relation with (2.5), we get
$$
  (\varphi\circ\chi^{-1})'(w_0)w=(\varphi\circ\psi^{-1})'(v_0)(\psi\circ\chi^{-1})'(w_0)w \quad \mbox{for}\; w\in X_1.
$$
  Since $(\varphi\circ\chi^{-1})'(w_0),(\varphi\circ\psi^{-1})'(v_0),(\psi\circ\chi^{-1})'(w_0)\in L(X)$, and clearly $X_1$ is
  dense in $X$, it follows that the above relation also holds for all $w\in X$. This proves that (2.2) holds for all $\eta\in
  \mathfrak{M}_1$. For a general $\eta\in\mathfrak{M}_0$, we use density of $\mathfrak{M}_1$ in $\mathfrak{M}_0$ and continuity
  of the three operators appearing in (2.2) with respect to $\eta$ in the topology of $\mathfrak{M}_0$. Hence (2.2) is proved.
  Having proved (2.2), by applying it to the special case $\chi=\varphi$ (note that in this case $w_0=u_0$), we see that
  $(\varphi\circ\psi^{-1})'(v_0)$ is a left-inverse of $(\psi\circ\varphi^{-1})'(u_0)$, and $(\psi\circ\varphi^{-1})'(u_0)$
  is a right-inverse of $(\varphi\circ\psi^{-1})'(v_0)$. Interchanging the roles of $\varphi$ and $\psi$, we see that
  $(\varphi\circ\psi^{-1})'(v_0)$ is also a right-inverse of $(\psi\circ\varphi^{-1})'(u_0)$. Hence (2.3) follows.
  $\quad\Box$
\medskip

  {\bf Definition 2.4}\ \ {\em Let $\mathfrak{M}$ be a Banach manifold built on the Banach space $X$ and $\mathfrak{M}_0$ a
  $C^1$-embedded Banach submanifold of $\mathfrak{M}$ built on the Banach space $X_0$, where $X_0$ is a densely embedded
  Banach subspace of $X$. Let $\eta\in\mathfrak{M}_0$. We define the following notions:

  $(1)$\ Let $\mathcal{O}$ be a neighborhood of $\eta$ in $\mathfrak{M}_0$ and $F:\mathcal{O}\to{\mathbb{R}}$ a real-valued
  function defined in $\mathcal{O}$. We say $F$ is \mbox{\small\boldmath $di\!f\!ferentiable$} at $\eta$ if there exists a
  $\mathfrak{M}_0$-regular local chart $(\mathcal{U},\varphi)$ of $\mathfrak{M}$ at $\eta$, such that the function
  $F\circ\varphi^{-1}:\varphi(\mathcal{O}\cap\mathcal{U})\to{\mathbb{R}}$ is differentiable at $\varphi(\eta)$ in the topology
  of $X_0$ and $(F\circ\varphi^{-1})'(u)\in X^*=L(X,\mathbb{R})$. We say $F$ is \mbox{\small\boldmath $continuously$
  $di\!f\!ferentiable$} at $\eta$ if there exist a neighborhood $\mathcal{O}'\subseteq\mathcal{O}$ of $\eta$ in $\mathfrak{M}_0$
  and a $\mathfrak{M}_0$-regular local chart $(\mathcal{U},\varphi)$ of $\mathfrak{M}$ at $\eta$, such that the function
  $F\circ\varphi^{-1}:\varphi(\mathcal{O}\cap\mathcal{U})\to{\mathbb{R}}$ is differentiable in $O=\varphi(\mathcal{O}'\cap
  \mathcal{U})$ in the topology of $X_0$ and $[u\mapsto(F\circ\varphi^{-1})'(u)]\in C(O,X^*)=C(O,L(X,\mathbb{R}))$, where $O$
  uses the topology of $X_0$. We denote by $\mathscr{D}^1_{\eta}$ the set of all real-valued functions $F$ defined in a
  neighborhood of $\eta$ in $\mathfrak{M}_0$ which are continuously differentiable at $\eta$.

  $(2)$\ Let $\mathcal{O}$ and $F$ be as in $(1)$. We say $F$ is \mbox{\small\boldmath $fully$ $di\!f\!ferentiable$} at $\eta$
  if for any $\mathfrak{M}_0$-regular local chart $(\mathcal{U},\varphi)$ of $\mathfrak{M}$ at $\eta$, the function
  $F\circ\varphi^{-1}:\varphi(\mathcal{O}\cap\mathcal{U})\to{\mathbb{R}}$ is differentiable at $\varphi(\eta)$ in the topology
  of $X_0$ and $(F\circ\varphi^{-1})'(u)\in X^*=L(X,\mathbb{R})$. We say $F$ is \mbox{\small\boldmath $fully$ $continuously$
  $di\!f\!ferentiable$} at $\eta$ if there exists a neighborhood $\mathcal{O}'\subseteq\mathcal{O}$ of $\eta$ in $\mathfrak{M}_0$
  such that for any $\mathfrak{M}_0$-regular local chart $(\mathcal{U},\varphi)$ of $\mathfrak{M}$ at $\eta$, the function
  $F\circ\varphi^{-1}:\varphi(\mathcal{O}\cap\mathcal{U})\to{\mathbb{R}}$ is differentiable in $O=\varphi(\mathcal{O}'\cap
  \mathcal{U})$ in the topology of $X_0$ and $[u\mapsto(F\circ\varphi^{-1})'(u)]\in C(O,X^*)=C(O,L(X,\mathbb{R}))$, where
  $O$ uses the topology of $X_0$. We denote by $\dot{\mathscr{D}}^1_{\eta}$ the set of all real-valued functions $F$ defined in
  a neighborhood of $\eta$ in $\mathfrak{M}_0$ which are fully continuously differentiable at $\eta$.

  $(3)$\ Let $\mathcal{O}$ be a neighborhood of $\eta$ in $\mathfrak{M}$ and $F:\mathcal{O}\to{\mathbb{R}}$ a real-valued
  function defined in $\mathcal{O}$. We say $F$ is \mbox{\small\boldmath $strongly$ $di\!f\!ferentiable$} at $\eta$ if there
  exists a $\mathfrak{M}_0$-regular local chart $(\mathcal{U},\varphi)$ of $\mathfrak{M}$ at $\eta$, such that the function
  $F\circ\varphi^{-1}:\varphi(\mathcal{O}\cap\mathcal{U})\to{\mathbb{R}}$ is differentiable at $\varphi(\eta)$ in the topology
  of $X$, so that $(F\circ\varphi^{-1})'(\varphi(\eta))\in X^*=L(X,\mathbb{R})$. We say $F$ is \mbox{\small\boldmath $strongly$
  $continuously$ $di\!f\!ferentiable$} at $\eta$ if there exist a neighborhood $\mathcal{O}'\subseteq\mathcal{O}$ of $\eta$ in
  $\mathfrak{M}$ and a $\mathfrak{M}_0$-regular local chart $(\mathcal{U},\varphi)$ of $\mathfrak{M}$ at $\eta$, such that the
  function $F\circ\varphi^{-1}: \varphi(\mathcal{O}\cap\mathcal{U})\to{\mathbb{R}}$ is differentiable in $O=
  \varphi(\mathcal{O}'\cap\mathcal{U})$ in the topology of $X$ and $[u\mapsto(F\circ\varphi^{-1})'(u)]\in C(O,X^*)=
  C(O,L(X,\mathbb{R}))$, where $O$ uses the topology of $X$. We denote by
  $\mathscr{D}^{1\!s}_{\eta}$ the set of all real-valued functions $F$ defined in a neighborhood of $\eta$ in $\mathfrak{M}$
  which are strongly continuously differentiable at $\eta$.

  $(4)$\ Let $\mathcal{O}$ and $F$ be as in $(3)$. We say $F$ is \mbox{\small\boldmath $fully$ $strongly$ $di\!f\!ferentiable$}
  at $\eta$ if for any $\mathfrak{M}_0$-regular local chart $(\mathcal{U},\varphi)$ of $\mathfrak{M}$ at $\eta$, the function
  $F\circ\varphi^{-1}:\varphi(\mathcal{O}\cap\mathcal{U})\to{\mathbb{R}}$ is differentiable at $\varphi(\eta)$ in the topology
  of $X$, so that $(F\circ\varphi^{-1})'(\varphi(\eta))\in X^*=L(X,\mathbb{R})$. We say $F$ is \mbox{\small\boldmath $fully$
  $strongly$ $continuously$ $di\!f\!ferentiable$} at $\eta$ if there exists a neighborhood $\mathcal{O}'\subseteq\mathcal{O}$
  of $\eta$ in $\mathfrak{M}$ such that for any $\mathfrak{M}_0$-regular local chart $(\mathcal{U},\varphi)$ of $\mathfrak{M}$
  at $\eta$, the function $F\circ\varphi^{-1}:\varphi(\mathcal{O}\cap\mathcal{U})\to{\mathbb{R}}$ is differentiable in
  $O=\varphi(\mathcal{O}'\cap\mathcal{U})$ in the topology of $X$ and $[u\mapsto(F\circ\varphi^{-1})'(u)]\in C(O,X^*)=
  C(O,L(X,\mathbb{R}))$, where $O$ uses the topology of $X$. We denote by $\dot{\mathscr{D}}^{1\!s}_{\eta}$ the set of all
  real-valued functions $F$ defined in a neighborhood of $\eta$ in $\mathfrak{M}$ which are fully strongly continuously
  differentiable at $\eta$.}
\medskip

  {\em Remark}.\ We make several remarks here. Firstly, we note that the condition ``$F\circ\varphi^{-1}$ is differentiable at
  $\varphi(\eta)$ in the topology of $X_0$'' implies that $(F\circ\varphi^{-1})'(\varphi(\eta))\in X_0^*=L(X_0,\mathbb{R})
  \supseteq X^*=L(X,\mathbb{R})$. The definitions (1) and (2) impose a stronger condition on $F\circ\varphi^{-1}$ than this:
  $(F\circ\varphi^{-1})'(\varphi(\eta))$ can be extended into a continuous linear functional in $X$, so that after
  extension $(F\circ\varphi^{-1})'(\varphi(\eta))\in X^*=L(X,\mathbb{R})$. Note that density of $X_0$ in $X$ ensures that the
  extension is unique. Secondly, we note that the following deduction relations hold:
$$
  \mbox{\em fully strongly differentiable}\;\Rightarrow\;\mbox{\em strongly differentiable}\;\Rightarrow\;
  \mbox{\em fully differentiable}\;\Rightarrow\;\mbox{\em differentiable}.
$$
  The first and the last deduction relations are obvious. The proof of the middle one is also easy (see the proof of the
  assertion (1) of Lemma 2.5 below). Thirdly, we note that if for any two neighborhoods $\mathcal{O}_j$ ($j=1,2$) of $\eta$
  in $\mathfrak{M}_0$ and two real-valued functions $F_j:\mathcal{O}_j\to {\mathbb{R}}$ ($j=1,2$), we define their linear
  combination $\alpha_1F_1+\alpha_2F_2$ for two real numbers $\alpha_1,\alpha_2$ to be the function $\alpha_1F_1+\alpha_2F_2:
  \mathcal{O}'\cap\mathcal{O}_2\to {\mathbb{R}}$ such that
$$
  (\alpha_1F_1+\alpha_2F_2)(\xi)=\alpha_1F_1(\xi)+\alpha_2F_2(\xi), \quad \forall\xi\in\mathcal{O}'\cap\mathcal{O}_2,
$$
  then clearly $\dot{\mathscr{D}}^1_{\eta}$ and $\dot{\mathscr{D}}^{1\!s}_{\eta}$ are linear spaces. Note that
  $\mathscr{D}^1_{\eta}$ and $\mathscr{D}^{1\!s}_{\eta}$ might not be linear spaces. Finally, assume that
  $(\mathfrak{M},\mathfrak{M}_0)$ is outward spreadable and $\widetilde{\mathfrak{M}}$ is a shell. Let $\widetilde{X}$ be
  the base Banach space of $\widetilde{\mathfrak{M}}$. For any $h\in\widetilde{X}^*$ and any $\mathfrak{M}_0$-regular local
  chart $(\mathcal{U},\varphi)$ of $\mathfrak{M}$ at $\eta$, it can be easily seen that $h\circ\varphi\in
  \dot{\mathscr{D}}^{1\!s}_{\eta}$. Hence $\dot{\mathscr{D}}^{1\!s}_{\eta}$ contains many functions.
\medskip

  {\bf Lemma 2.5}\ \ {\em Let $\mathfrak{M}$ and $\mathfrak{M}_0$ be two Banach manifolds such that $\mathfrak{M}_0$ is
  a $C^1$-embedded Banach submanifold of $\mathfrak{M}$. Assume that $(\mathfrak{M},\mathfrak{M}_0)$ is inward spreadable
  and $\mathfrak{M}_1$ is an inner $C^1$-kernel. Let $\mathcal{O}$ be an open subset of $\mathfrak{M}_0$ and
  $F:\mathcal{O}\to{\mathbb{R}}$ a real-valued function defined in $\mathcal{O}$. Then we have the following assertions:

  $(1)$\ If $\eta\in\mathcal{O}$ and $F$ is strongly differentiable at $\eta$, then for any two $\mathfrak{M}_0$-regular
  local charts $(\mathcal{U},\varphi)$ and $(\mathcal{V},\psi)$ of $\mathfrak{M}$ at $\eta$, the following relation holds:
\begin{equation}
  (F\circ\psi^{-1})'(\psi(\eta))=(F\circ\varphi^{-1})'(\varphi(\eta))(\varphi\circ\psi^{-1})'(\psi(\eta)).
\end{equation}
%---(2.7)---

  $(2)$\ If $\eta\in\mathcal{O}\cap\mathfrak{M}_1$ and $F$ is differentiable at $\eta$, then for any two
  $(\mathfrak{M}_0,\mathfrak{M}_1)$-regular local charts $(\mathcal{U},\varphi)$ and $(\mathcal{V},\psi)$ of $\mathfrak{M}$
  at $\eta$, the relation $(2.7)$ holds.

  $(3)$\ If $F$ is continuously differentiable at $\eta$, then for any $\eta\in\mathcal{O}$ and any two
  $(\mathfrak{M}_0,\mathfrak{M}_1)$-regular local charts $(\mathcal{U},\varphi)$ and $(\mathcal{V},\psi)$ of $\mathfrak{M}$ at
  $\eta$, the relation $(2.7)$ holds.}
\medskip

  {\em Proof}.\ (1)\ Let $X$ and $X_0$ be the base spaces of $\mathfrak{M}$ and $\mathfrak{M}_0$, respectively. Since $F$
  is strongly differentiable at $\eta$, there exists a $\mathfrak{M}_0$-regular local chart
  $(\mathcal{W},\chi)$ of $\mathfrak{M}$ at $\eta$, such that the function $F\circ\chi^{-1}:\chi(\mathcal{O}\cap
  \mathcal{W})\to{\mathbb{R}}$ is differentiable at $\chi(\eta)$ in the topology of $X$. For any other $\mathfrak{M}_0$-regular
  local chart $(\mathcal{V},\psi)$ of $\mathfrak{M}$ at $\eta$, we write
$$
  F\circ\psi^{-1}=(F\circ\chi^{-1})\circ(\chi\circ\psi^{-1}).
$$
  The condition that $\mathfrak{M}_0$ is a $C^1$-embedded Banach submanifold of $\mathfrak{M}$ implies that $\chi\circ\psi^{-1}$
  is differentiable as a mapping from $X_0$ to $X$. Hence the above relation shows that $F\circ\psi^{-1}$ is differentiable at
  $v_0=\psi(\eta)$ in the topology of $X_0$ (so that $F$ is fully differentiable at $\eta$ by arbitrariness of $(\mathcal{V},\psi)$),
  and
$$
  (F\circ\psi^{-1})'(\psi(\eta))v=(F\circ\chi^{-1})'(\chi(\eta))(\chi\circ\psi^{-1})'(\psi(\eta))v,
  \quad \forall v\in X_0.
$$
  Since $X_0$ is dense in $X$, $(F\circ\chi^{-1})'(\chi(\eta))\in L(X,\mathbb{R})$ and $(\chi\circ\psi^{-1})'(\psi(\eta))
  \in L(X)$, it follows that $(F\circ\psi^{-1})'(\psi(\eta))\in L(X,\mathbb{R})$ and the above relation also holds for all
  $v\in X$ (note that this does not imply that $F\circ\psi^{-1}$ is differentiable at $\psi(\eta)$ in the topology of $X$).
  Hence
$$
  (F\circ\psi^{-1})'(\psi(\eta))=(F\circ\chi^{-1})'(\chi(\eta))(\chi\circ\psi^{-1})'(\psi(\eta)).
$$
  Applying this assertion to any other $\mathfrak{M}_0$-regular local chart $(\mathcal{U},\varphi)$ of $\mathfrak{M}$
  at $\eta$, we also have
$$
  (F\circ\varphi^{-1})'(\varphi(\eta))=(F\circ\chi^{-1})'(\chi(\eta))(\chi\circ\varphi^{-1})'(\varphi(\eta)),
$$
  which implies that
$$
  (F\circ\chi^{-1})'(\chi(\eta))=(F\circ\varphi^{-1})'(\varphi(\eta))[(\chi\circ\varphi^{-1})'(\varphi(\eta))]^{-1}.
$$
  Hence
\begin{align}
  (F\circ\psi^{-1})'(\psi(\eta))=&(F\circ\varphi^{-1})'(\varphi(\eta))
  [(\chi\circ\varphi^{-1})'(\varphi(\eta))]^{-1}(\chi\circ\psi^{-1})'(\psi(\eta))
\nonumber\\
  =&(F\circ\varphi^{-1})'(\varphi(\eta))(\varphi\circ\psi^{-1})'(\psi(\eta)).
\nonumber
\end{align}
  In getting the last equality we used Lemma 2.3. This proves (2.7).

  (2)\ The proof of the assertion (2) is similar to that of the assertion (1) with some minor modification. We omit it here.

  (3)\ Let $\eta\in\mathcal{O}$. Let $\mathcal{O}'$ be as in Definition 2.4 (2). Choose a sequence $\eta_j\in\mathcal{O}'
  \cap\mathfrak{M}_1$ ($j=1,2,\cdots$) such that $\eta_j\to\eta$ in the topology of $\mathfrak{M}_0$. By the assertion (2),
  for every $\eta_j$ we have
$$
  (F\circ\psi^{-1})'(\psi(\eta_j))=(F\circ\varphi^{-1})'(\varphi(\eta_j))(\varphi\circ\psi^{-1})'(\psi(\eta_j)).
$$
  Since $[x\mapsto (F\circ\varphi^{-1})'(x)]\in C(\varphi(\mathcal{U}\cap\mathcal{O}'),X^*)$, $[x\mapsto
  (F\circ\psi^{-1})'(x)]\in C(\psi(\mathcal{U}\cap\mathcal{O}'),X^*)$, and $[x\mapsto (\varphi\circ\psi^{-1})'(x)\in
  C(\psi(\mathcal{U}\cap\mathcal{O}),L(X))$, by letting $j\to\infty$ in the above relation, we see that (2.7) follows.
  This proves the assertion (3). $\quad\Box$
\medskip

  {\bf Definition 2.6}\ \ {\em Let $\mathfrak{M}$ be a Banach manifold built on the Banach space $X$ and $\mathfrak{M}_0$ a
  $C^1$-embedded Banach submanifold of $\mathfrak{M}$ built on the Banach space $X_0$, where $X_0$ is an embedded Banach
  subspace of $X$. Let $\eta\in\mathfrak{M}_0$.

  $(1)$ Let $f:(-\varepsilon,\varepsilon)\to\mathfrak{M}$ $(\varepsilon>0)$ be a curve in $\mathfrak{M}$ passing $\eta$,
  i.e., $f(0)=\eta$. We say $f(t)$ is \mbox{\small\boldmath $di\!f\!ferentiable$} at $t=0$ if there exists a
  $\mathfrak{M}_0$-regular local chart $(\mathcal{U},\varphi)$ of $\mathfrak{M}$ at $\eta$, such that the function
  $t\mapsto\varphi(f(t))$ is differentiable at $t=0$ in the topology of $X$, and if $(\mathcal{V},\psi)$ is any other
  $\mathfrak{M}_0$-regular local chart of $\mathfrak{M}$ at $\eta$ such that the function $t\mapsto\psi(f(t))$ is
  differentiable at $t=0$, then the following relation holds:
\begin{equation}
  (\psi\circ f)'(0)=(\psi\circ\varphi^{-1})'(\varphi(\eta))(\varphi\circ f)'(0).
\end{equation}
%---(2.8)---
  Moreover, we define the \mbox{\small\boldmath $tangent$ $vector$} of this curve at $\eta$, or the \mbox{\small\boldmath
  $derivative$} $f'(0)$ of $f(t)$ at $t=0$, to be the mapping $f'(0):\dot{\mathscr{D}}^{1\!s}_{\eta}\to{\mathbb{R}}$
  defined by
\begin{equation}
  f'(0)F=(F\circ f)'(0), \quad \forall F\in\dot{\mathscr{D}}^{1\!s}_{\eta}.
\end{equation}
%---(2.9)---
  We say $f(t)$ is \mbox{\small\boldmath $continuously$ $differentiable$} at $t=0$ if by taking $\varepsilon>0$ smaller
  when necessary, there exists a $\mathfrak{M}_0$-regular local chart $(\mathcal{U},\varphi)$ of $\mathfrak{M}$ at $\eta$,
  such that the function $t\mapsto\varphi(f(t))$ is continuously differentiable for all $|t|<\varepsilon$ in the topology
  of $X$ and the relation $(2.8)$ holds for any $\mathfrak{M}_0$-regular local chart $(\mathcal{V},\psi)$ of $\mathfrak{M}$
  at $\eta$ such that the function $t\mapsto\psi(f(t))$ is continuously differentiable at $t=0$.

  $(2)$ Let $f:(-\varepsilon,\varepsilon)\to\mathfrak{M}$ $(\varepsilon>0)$ be a curve in $\mathfrak{M}$ passing $\eta$,
  i.e., $f(0)=\eta$. We say $f(t)$ is \mbox{\small\boldmath $fully$ $di\!f\!ferentiable$} at $t=0$ if for any
  $\mathfrak{M}_0$-regular local chart $(\mathcal{U},\varphi)$ of $\mathfrak{M}$ at $\eta$, the function $t\mapsto
  \varphi(f(t))$ is differentiable at $t=0$ in the topology of $X$, and for any two $\mathfrak{M}_0$-regular local charts
  $(\mathcal{U},\varphi)$ and $(\mathcal{V},\psi)$ of $\mathfrak{M}$ at $\eta$, the relation $(2.8)$ holds. In this case,
  the tangent vector $f'(0)$ can be re-defined as a mapping $f'(0):\mathscr{D}^{1\!s}_{\eta}\to{\mathbb{R}}$ as follows:
$$
  f'(0)F=(F\circ f)'(0), \quad \forall F\in\mathscr{D}^{1\!s}_{\eta}.
$$
  We say $f(t)$ is \mbox{\small\boldmath $fully$ $continuously$ $di\!f\!ferentiable$} at $t=0$ if by taking $\varepsilon>0$
  smaller when necessary, for any $\mathfrak{M}_0$-regular local chart $(\mathcal{U},\varphi)$ of $\mathfrak{M}$ at $\eta$
  the function $t\mapsto\varphi(f(t))$ is continuously differentiable for all $|t|<\varepsilon$, and for any two
  $\mathfrak{M}_0$-regular local charts $(\mathcal{U},\varphi)$ and $(\mathcal{V},\psi)$ of $\mathfrak{M}$ at $\eta$, the
  relation $(2.8)$ holds.

  $(3)$ Let $f:(-\varepsilon,\varepsilon)\to\mathfrak{M}_0$ $(\varepsilon>0)$ be a curve in $\mathfrak{M}_0$ passing $\eta$,
  i.e., $f(0)=\eta$. We say $f(t)$ is \mbox{\small\boldmath $strongly$ $di\!f\!ferentiable$} at $t=0$ if there exists a
  $\mathfrak{M}_0$-regular local chart $(\mathcal{U},\varphi)$ of $\mathfrak{M}$ at $\eta$ such that the function
  $t\mapsto\varphi(f(t))$ is differentiable at $t=0$ in the topology of $X_0$. In this case, the tangent vector $f'(0)$ can
  be re-defined as a mapping $f'(0):\dot{\mathscr{D}}^1_{\eta}\to{\mathbb{R}}$ as follows:
$$
  f'(0)F=(F\circ f)'(0), \quad \forall F\in\dot{\mathscr{D}}^1_{\eta}.
$$
  We say $f(t)$ is \mbox{\small\boldmath $strongly$ $continuously$ $di\!f\!ferentiable$} at $t=0$ if by taking $\varepsilon>0$
  smaller when necessary, there exists a $\mathfrak{M}_0$-regular local chart $(\mathcal{U},\varphi)$ of $\mathfrak{M}$ at
  $\eta$, such that the function $t\mapsto\varphi(f(t))$ is continuously differentiable for all $|t|<\varepsilon$ in the
  topology of $X_0$.

  $(4)$ Let $f:(-\varepsilon,\varepsilon)\to\mathfrak{M}_0$ $(\varepsilon>0)$ be a curve in $\mathfrak{M}_0$ passing $\eta$,
  i.e., $f(0)=\eta$. We say $f(t)$ is \mbox{\small\boldmath $fully$ $strongly$ $di\!f\!ferentiable$} at $t=0$ if for any
  $\mathfrak{M}_0$-regular local chart $(\mathcal{U},\varphi)$ of $\mathfrak{M}$ at $\eta$, the function $t\mapsto
  \varphi(f(t))$ is differentiable at $t=0$ in the topology of $X_0$. In this case, the tangent vector $f'(0)$ can be re-defined
  as a mapping $f'(0):\mathscr{D}^1_{\eta}\to{\mathbb{R}}$ as follows:
$$
  f'(0)F=(F\circ f)'(0), \quad \forall F\in\mathscr{D}^1_{\eta}.
$$
  We say $f(t)$ is \mbox{\small\boldmath $fully$ $strongly$ $continuously$ $di\!f\!ferentiable$} at $t=0$ if by taking
  $\varepsilon>0$ smaller when necessary, for any $\mathfrak{M}_0$-regular local chart $(\mathcal{U},\varphi)$ of
  $\mathfrak{M}$ at $\eta$, the function $t\mapsto\varphi(f(t))$ is continuously differentiable for all $|t|<\varepsilon$ in
  the topology of $X_0$.}

  {\em Remark}.\ We make four remarks here. (i) The following deduction relations are obvious:
$$
  \mbox{\em fully strongly differentiable}\;\Rightarrow\;\mbox{\em strongly differentiable}; \qquad
  \mbox{\em fully differentiable}\;\Rightarrow\;\mbox{\em differentiable}.
$$
  For a quasi-differentiable Banach manifold there also holds the following deduction relation:
$$
  \mbox{\em strongly differentiable}\;\Rightarrow\;\mbox{\em fully differentiable};
$$
  see Lemma 2.7 below. (ii) The definition (2.9) makes sense. Indeed, since $f(t)$ is differentiable at $t=0$, there exists
  a $\mathfrak{M}_0$-regular local chart $(\mathcal{U},\varphi)$ of $\mathfrak{M}$ at $\eta$, such that the function
  $t\mapsto\varphi(f(t))$ is differentiable at $t=0$ in the topology of $X$. Let $F\in\dot{\mathscr{D}}^{1s}_{\eta}$. Then
  $F\circ\varphi^{-1}$ is differentiable at $\varphi(\eta)$ in the topology of $X$. Since $F\circ f=(F\circ\varphi^{-1})
  \circ(\varphi\circ f)$, we see that $F\circ f$ is differentiable at $t=0$, and
$$
  (F\circ f)'(0)=(F\circ\varphi^{-1})'(\varphi(\eta))(\varphi\circ f)'(0).
$$
  If $(\mathcal{V},\psi)$ is another $\mathfrak{M}_0$-regular local chart of $\mathfrak{M}$ at $\eta$, such that the function
  $t\mapsto\psi(f(t))$ is differentiable at $t=0$ in the topology of $X$, then similarly we have
$$
  (F\circ f)'(0)=(F\circ\psi^{-1})'(\psi(\eta))(\psi\circ f)'(0).
$$
  By (2.8) and (2.7), the right-hand sides of the above two equations are equal. Hence (2.9) makes sense.
  Similarly, the equations in definitions (2), (3) and (4) make sense. (iii) It can be easily
  seen that if $f(t)$ is strongly differentiable at $t=0$ then $f'(0)$ is a linear functional in $\dot{\mathscr{D}}^1_{\eta}$.
  (iv) There are many strongly continuously differentiable functions. Indeed, for $\eta\in\mathfrak{M}_0$ and any $v\in X_0$,
  arbitrarily choose a $\mathfrak{M}_0$-regular local chart $(\mathcal{U},\varphi)$ of $\mathfrak{M}$ at $\eta$ and let
  $f(t)=\varphi^{-1}(\varphi(\eta)+tv)$ for $|t|<\varepsilon$, where $\varepsilon>0$ is sufficiently small so that
  $\varphi(\eta)+tv\in\varphi(\mathcal{U}\cap\mathfrak{M}_0)$ for all $|t|<\varepsilon$. Then $f(t)$ is strongly continuously
  differentiable at $t=0$. (v) Let $(\mathfrak{M},\mathfrak{M}_0)$ be inward spreadable with an inner $C^1$-kernel
  $\mathfrak{M}_1$. Let $X$, $X_0$ and $X_1$ be the base spaces of
  $\mathfrak{M}$, $\mathfrak{M}_0$ and $\mathfrak{M}_1$, respectively. Let $\eta\in\mathfrak{M}_1$. For an arbitrary $v\in X_1$
  and an arbitrary $(\mathfrak{M}_0,\mathfrak{M}_1)$-regular local chart $(\mathcal{U},\varphi)$ of $\mathfrak{M}$ at $\eta$,
  let $f(t)=\varphi^{-1}(\varphi(\eta)+tv)$ for $|t|<\varepsilon$, where $\varepsilon>0$ is sufficiently small so that
  $\varphi(\eta)+tv\in\varphi(\mathcal{U}\cap\mathfrak{M}_1)$ for all $|t|<\varepsilon$. Then it can be easily seen that $f(t)$
  is fully strongly continuously differentiable at $t=0$.
\medskip

  {\bf Lemma 2.7}\ \ {\em Let $\mathfrak{M}$ be a Banach manifold and $\mathfrak{M}_0$ a $C^1$-embedded Banach submanifold of
  $\mathfrak{M}$. Assume that $(\mathfrak{M},\mathfrak{M}_0)$ is inward spreadable. Let $\eta\in\mathfrak{M}_0$. Let
  $f:(-\varepsilon,\varepsilon)\to\mathfrak{M}_0$ $(\varepsilon>0)$ be a curve in $\mathfrak{M}_0$ passing $\eta$, i.e.
  $f(0)=\eta$, which is strongly differentiable at $t=0$. Then $f(t)$ is fully differentiable at $t=0$, i.e., for any
  $\mathfrak{M}_0$-regular local chart $(\mathcal{U},\varphi)$ of $\mathfrak{M}$ at $\eta$, the function $t\mapsto
  \varphi(f(t))$ is differentiable at $t=0$ in the topology of $X$, and for any two $\mathfrak{M}_0$-regular local charts
  $(\mathcal{U},\varphi)$ and $(\mathcal{V},\psi)$ of $\mathfrak{M}$ at $\eta$, the relation $(2.8)$ holds.}
\medskip

  {\em Proof}.\ Let $\mathfrak{M}_1$ be an inner $C^1$-kernel. Let $X$, $X_0$ and $X_1$ be the base spaces of $\mathfrak{M}$,
  $\mathfrak{M}_0$ and $\mathfrak{M}_1$, respectively. The condition that $f(t)$ is strongly differentiable at $t=0$ implies
  that there exists a $\mathfrak{M}_0$-regular local chart $(\mathcal{W},\chi)$ of $\mathfrak{M}$ at $\eta$ such that the
  function $t\mapsto(\chi\circ f)(t)$ is differentiable at $t=0$ in the topology of $X_0$. Let $(\mathcal{U},\varphi)$ be an
  arbitrary $\mathfrak{M}_0$-regular local chart of $\mathfrak{M}$ at $\eta$. Since $\varphi\circ f=(\varphi\circ \chi^{-1})
  \circ(\chi\circ f)$ and $\varphi\circ \chi^{-1}$ is differentiable at $\chi(\eta)$ as a map from $X_0$ to $X$, it follows
  that the functions $t\mapsto(\varphi\circ f)(t)$ is differentiable at $t=0$ in the topology of $X$. Moreover, for any two
  $\mathfrak{M}_0$-regular local charts $(\mathcal{U},\varphi)$ and $(\mathcal{V},\psi)$ of $\mathfrak{M}$ at $\eta$, we
  have
$$
  (\varphi\circ f)'(0)=(\varphi\circ\chi^{-1})'(\chi(\eta))(\chi\circ f)'(0),
$$
$$
  (\psi\circ f)'(0)=(\psi\circ\chi^{-1})'(\chi(\eta))(\chi\circ f)'(0).
$$
  Applying the relation (2.2), we see that $(\psi\circ\chi^{-1})'(\chi(\eta))=(\psi\circ\varphi^{-1})'(\varphi(\eta))
  (\varphi\circ\chi^{-1})'(\chi(\eta))$, so that
$$
\begin{array}{rl}
  (\psi\circ f)'(0)=&(\psi\circ\varphi^{-1})'(\varphi(\eta))(\varphi\circ\chi^{-1})'(\chi(\eta))(\chi\circ f)'(0)
\\
  =&(\psi\circ\varphi^{-1})'(\varphi(\eta))(\varphi\circ f)'(0).
\end{array}
$$
  This proves (2.8). $\quad\Box$
\medskip

  {\bf Lemma 2.8}\ \ {\em Let $\mathfrak{M}$ be a Banach manifold built on the Banach space $X$ and $\mathfrak{M}_0$ a
  $C^1$-embedded Banach submanifold of $\mathfrak{M}$. If $(\mathfrak{M},\mathfrak{M}_0)$ is outward spreadable, then for
  any $\eta\in\mathfrak{M}_0$, any $\mathfrak{M}_0$-regular local chart $(\mathcal{U},\varphi)$ of $\mathfrak{M}$ at
  $\eta$, and any function $f:(-\varepsilon,\varepsilon)\to\mathfrak{M}$ with the property $f(0)=\eta$, if the function
  $\varphi\circ f:(-\varepsilon,\varepsilon)\to X$ is differentiable at $t=0$ $($in this case we say $f$ is
  $\varphi$-differentiable at $t=0$ later on$)$, then $f$ is differentiable at $t=0$, i.e., for any $\mathfrak{M}_0$-regular
  local chart $(\mathcal{V},\psi)$ of $\mathfrak{M}$ at $\eta$ such that $\psi\circ f:(-\varepsilon,\varepsilon)\to X$
  is differentiable at $t=0$, the relation $(2.8)$ holds.}
\medskip

  {\em Proof}.\ The condition of this lemma ensures that there exists a Banach manifold $\widetilde{\mathfrak{M}}$ such that
  $\mathfrak{M}$ is a $C^1$-embedded Banach submanifold of $\widetilde{\mathfrak{M}}$ and $(\widetilde{\mathfrak{M}},
  \mathfrak{M})$ is inward spreadable, with $C^1$-kernel $\mathfrak{M}$ and inner $C^1$-kernel $\mathfrak{M}_0$. Assume that
  $\widetilde{\mathfrak{M}}$, $\mathfrak{M}$ and $\mathfrak{M}_0$ are built on the Banach spaces $\widetilde{X}$, $X$ and
  $X_0$, respectively. Let $\eta$, $\varphi$ and $f$ be as in the lemma. Let $(\widetilde{\mathcal{U}},\widetilde{\varphi})$
  be the corresponding ($\mathfrak{M},\mathfrak{M}_0$)-regular local chart of $\widetilde{\mathfrak{M}}$ at $\eta$, i.e.,
  $(\widetilde{\mathcal{U}},\widetilde{\varphi})$ is a local chart of $\widetilde{\mathfrak{M}}$ at $\eta$ such that
  $\mathcal{U}=\widetilde{\mathcal{U}}|_{\mathfrak{M}}$ and $\varphi=\widetilde{\varphi}|_{\mathcal{U}}$. The condition
  that $f$ is $\varphi$-differentiable at $t=0$ implies that it is strongly differentiable at $t=0$ as a function in
  $\widetilde{\mathfrak{M}}$, so that by Lemma 2.7 we see that for any $\mathfrak{M}_0$-regular local chart
  $(\mathcal{V},\psi)$ of $\mathfrak{M}$ at $\eta$ there holds
$$
  (\psi\circ f)'(0)=(\psi\circ\varphi^{-1})'(\varphi(\eta))(\varphi\circ f)'(0).
$$
  Here $(\psi\circ f)'(0),(\varphi\circ f)'(0)$ are regarded as elements in $\widetilde{X}$, and
  $(\psi\circ\varphi^{-1})'(\varphi(\eta))$ is regarded as an element in $L(\widetilde{X})$. The condition that
  $f$ is $\varphi$-differentiable at $t=0$ implies that $(\varphi\circ f)'(0)\in X$. If there exists another
  $\mathfrak{M}_0$-regular local chart $(\mathcal{V},\psi)$ of $\mathfrak{M}$ at $\eta$ such that $f$ is $\psi$-differentiable
  at $t=0$, then also $(\psi\circ f)'(0)\in X$. Since $\eta\in\mathfrak{M}_0$ implies that $(\psi\circ\varphi^{-1})'(\varphi(\eta))
  \in L(X)$, the above relation ensures that (2.8) holds. This proves the desired assertion. $\quad\Box$
\medskip

  {\bf Definition 2.9}\ \ {\em Let $\mathfrak{M}$ be a Banach manifold built on the Banach space $X$ and $\mathfrak{M}_0$ a
  $C^1$-embedded Banach submanifold of $\mathfrak{M}$ built on the Banach space $X_0$, where $X_0$ is a densely embedded
  Banach subspace of $X$. Assume that $(\mathfrak{M},\mathfrak{M}_0)$ is outward spreadable. Let $\eta\in\mathfrak{M}_0$.

  $(1)$\ \ We denote
$$
  \mathcal{T}_{\eta}(\mathfrak{M})=\{f'(0):f:(-\varepsilon,\varepsilon)\to\mathfrak{M},\; f(0)=\eta,\;
  \mbox{$f(t)$ is continuously differentiable at $t=0$}\}.
$$

  $(2)$\ \ If $(\mathcal{U},\varphi)$ is a $\mathfrak{M}_0$-regular local chart of $\mathfrak{M}$ at $\eta$, we define
  the derivative of $\varphi:\mathcal{U}\cap\mathfrak{M}_0\to X_0\subseteq X$ at $\eta$ to be a mapping $\varphi'(\eta):
  \mathcal{T}_{\eta}(\mathfrak{M})\to X$ such that
$$
  \varphi'(\eta)\upsilon=(\varphi\circ\psi^{-1})'(\psi(\eta))(\psi\circ f)'(0) \quad \mbox{for}\;\,
  \upsilon=f'(0)\in \mathcal{T}_{\eta}(\mathfrak{M}),
$$
  where $(\mathcal{V},\psi)$ is a $\mathfrak{M}_0$-regular local chart of $\mathfrak{M}$ at $\eta$ such that $\psi\circ f$ is
  continuously differentiable at $t=0$ in the topology of $X$.

  $(3)$\ \ Let $(\mathcal{U},\varphi)$ be as in $(2)$. We define a nonnegative-valued function $\|\cdot\|_{\varphi}$ on
  $\mathcal{T}_{\eta}(\mathfrak{M})$ as follows:}
$$
  \|\upsilon\|_{\varphi}=\|\varphi'(\eta)\upsilon\|_X, \quad \forall\upsilon\in\mathcal{T}_{\eta}(\mathfrak{M}),
$$

  {\bf Lemma 2.10}\ \ {\em Let $\mathfrak{M}$ be a quasi-differentiable Banach manifold built on the Banach space $X$.
  Let $\mathfrak{M}_0$ be the $C^1$-kernel. Then the following assertions hold:

  $(1)$\ \ The definition $(2)$ in Definition 2.9 makes sense, i.e., the value of $\varphi'(\eta)\upsilon$ does not
  depend on specific choice of the function $f$ such that $\upsilon=f'(0)$ and the $\mathfrak{M}_0$-regular local chart
  $(\mathcal{V},\psi)$ of $\mathfrak{M}$ at $\eta$ such that $\psi\circ f$ is continuously differentiable at $t=0$.

  $(2)$\ \ If $(\mathcal{U},\varphi)$ and $(\mathcal{V},\psi)$ are two different $\mathfrak{M}_0$-regular local charts of
  $\mathfrak{M}$ at $\eta$, then $\|\cdot\|_{\varphi}$ and $\|\cdot\|_{\psi}$ are equivalent when restricted to the set
$$
  \mathcal{T}_{\eta}^{00}(\mathfrak{M})=\{f'(0):f:(-\varepsilon,\varepsilon)\to\mathfrak{M}_0,\; f(0)=\eta,\;
  \mbox{$f(t)$ is fully strongly continuously differentiable at $t=0$}\}.
$$
  i.e., there exists constants $C_1,C_2>0$ such that
\begin{equation}
  \|\upsilon\|_{\psi}\leqslant C_1\|\upsilon\|_{\varphi}, \quad
  \|\upsilon\|_{\varphi}\leqslant C_2\|\upsilon\|_{\psi}, \quad
   \forall\upsilon\in\mathcal{T}_{\eta}^{00}(\mathfrak{M}).
\end{equation}
%---(2.10)---

  $(3)$\ \ For any $\mathfrak{M}_0$-regular local chart $(\mathcal{U},\varphi)$ of $\mathfrak{M}$ at $\eta$, the mapping
  $\varphi'(\eta):\mathcal{T}_{\eta}(\mathfrak{M})\to X$ is a bijection.}
\medskip

  {\em Proof}.\ (1)\ Let $f,g:(-\varepsilon,\varepsilon)\to\mathfrak{M}$ ($\varepsilon>0$) be two functions such that
  $f(0)=g(0)=\eta$, both continuously differentiable at $t=0$ and $f'(0)=g'(0)$. Let $(\mathcal{V},\psi)$ and $(\mathcal{W},\chi)$
  be two $\mathfrak{M}_0$-regular local charts of $\mathfrak{M}$ at $\eta$ such that $\psi\circ f$ is continuously differentiable
  at $t=0$ and $\chi\circ g$ is continuously differentiable at $t=0$. Since $f'(0)=g'(0)$, we have
$$
  (F\circ f)'(0)=(F\circ g)'(0)
$$
  for all $F\in\dot{\mathscr{D}}^{1\!s}_{\eta}$, which implies that
$$
  (F\circ\psi^{-1})'(\psi(\eta))(\psi\circ f)'(0)=(F\circ\chi^{-1})'(\chi(\eta))(\chi\circ g)'(0)
$$
  for all $F\in\dot{\mathscr{D}}^{1\!s}_{\eta}$.  By the assertion (1) of Lemma 2.5, this implies that for any $\mathfrak{M}_0$-regular
  local chart $(\mathcal{U},\varphi)$ of $\mathfrak{M}$ at $\eta$, there holds
\begin{equation}
  (F\circ\varphi^{-1})'(\varphi(\eta))(\varphi\circ\psi^{-1})'(\psi(\eta))(\psi\circ f)'(0)
  =(F\circ\varphi^{-1})'(\varphi(\eta))(\varphi\circ\chi^{-1})'(\chi(\eta))(\chi\circ g)'(0)
\end{equation}
%---(2.11)---
  for any $F\in\dot{\mathscr{D}}^{1\!s}_{\eta}$. Now let $\widetilde{\mathfrak{M}}$ be a shell of $(\mathfrak{M},\mathfrak{M}_0)$.
  Let $\widetilde{X}$ be the base Banach space of $\widetilde{\mathfrak{M}}$. For any $h\in\widetilde{X}^*$, let $F=h\circ\varphi$.
  It is easy to see that $F$ is fully strongly continuously differentiable at $\eta$, i.e., $F\in\dot{\mathscr{D}}^{1\!s}_{\eta}$.
  By applying (2.11) to $F=h\circ\varphi$, we get the relation
$$
  h[(\varphi\circ\psi^{-1})'(\psi(\eta))(\psi\circ f)'(0)]
  =h[(\varphi\circ\chi^{-1})'(\chi(\eta))(\chi\circ g)'(0)], \quad \forall h\in \widetilde{X}^*.
$$
  By arbitrariness of $h$ in $\widetilde{X}^*$, we obtain
$$
  (\varphi\circ\psi^{-1})'(\psi(\eta))(\psi\circ f)'(0)
  =(\varphi\circ\chi^{-1})'(\chi(\eta))(\chi\circ g)'(0),
$$
  which means $\varphi'(\eta)f'(0)=\varphi'(\eta)g'(0)$. This proves the assertion (1).

  (2)\ Let $C_1=\|(\psi\circ\varphi^{-1})'(\varphi(\eta))\|_{L(X)}$. Given $\upsilon\in\mathcal{T}_{\eta}^{00}(\mathfrak{M})$,
  let $f:(-\varepsilon,\varepsilon)\to\mathfrak{M}_0$ be such that $f(0)=\eta$, $f(t)$ is fully strongly differentiable
  at $t=0$ and $f'(0)=\upsilon$. Then we have
$$
\begin{array}{rl}
  \psi'(\eta)\upsilon=&(\psi\circ f)'(0)=(\psi\circ\varphi^{-1}\circ\varphi\circ f)'(0)\\
  =&(\psi\circ\varphi^{-1})'(\varphi(\eta))(\varphi\circ f)'(0)
  =(\psi\circ\varphi^{-1})'(\varphi(\eta))\varphi'(\eta)\upsilon,
\end{array}
$$
  which implies that
$$
  \|\upsilon\|_{\psi}=\|\psi'(\eta)\upsilon\|_{X}\leqslant\|(\psi\circ\varphi^{-1})'(\varphi(\eta))\|_{L(X)}
  \|\varphi'(\eta)\upsilon\|_{X}=C_1\|\upsilon\|_{\varphi}.
$$
  This proves the first inequality in (2.10). Similarly we see that the second inequality in (2.10) holds for
  $C_2=\|(\varphi\circ\psi^{-1})'(\varphi(\eta))\|_{L(X)}$.

  (3)\ We first prove that the mapping $\varphi'(\eta):\mathcal{T}_{\eta}(\mathfrak{M})\to X$ is an injection. Indeed,
  let $\upsilon_1,\upsilon_2\in\mathcal{T}_{\eta}(\mathfrak{M})$ be such that $\varphi'(\eta)\upsilon_1=
  \varphi'(\eta)\upsilon_2$. Let $f_1,f_2:(-\varepsilon,\varepsilon)\to\mathfrak{M}$ be two continuously differentiable
  functions such that $f_1(0)=f_2(0)=\eta$ and $f_1'(0)=\upsilon_1$, $f_2'(0)=\upsilon_2$. Let $(\mathcal{V},\psi)$ and
  $(\mathcal{W},\chi)$ be two $\mathfrak{M}_0$-regular local charts of $\mathfrak{M}$ at $\eta$ such that $[t\mapsto
  (\psi\circ f_1)(t)]$ and $[t\mapsto(\chi\circ f_2)(t)]$ are continuously differentiable at $t=0$, both in the topology
  of $X$. The condition $\varphi'(\eta)\upsilon_1=\varphi'(\eta)\upsilon_2$ implies
$$
  (\varphi\circ\psi^{-1})'(\psi(\eta))(\psi\circ f_1)'(0)=(\varphi\circ\chi^{-1})'(\chi(\eta))(\chi\circ f_2)'(0).
$$
  Let $F\in\dot{\mathscr{D}}^{1\!s}_{\eta}$. Then $F\circ\varphi^{-1}$ is continuously differentiable at $\varphi(\eta)$
  in the topology of $X$. From the above equality it follows that
$$
  (F\circ\varphi^{-1})'(\varphi(\eta))(\varphi\circ\psi^{-1})'(\psi(\eta))(\psi\circ f_1)'(0)
  =(F\circ\varphi^{-1})'(\varphi(\eta))(\varphi\circ\chi^{-1})'(\chi(\eta))(\chi\circ f_2)'(0).
$$
  By the assertion (3) of Lemma 2.5, this yields
$$
  (F\circ\psi^{-1})'(\psi(\eta))(\psi\circ f_1)'(0)=(F\circ\chi^{-1})'(\chi(\eta))(\chi\circ f_2)'(0).
$$
  This further implies
$$
  f_1'(0)F=(F\circ f_1)'(0)=f_2'(0)F,
$$
  i.e., $\upsilon_1F=\upsilon_2F$. By arbitrariness of $F\in\dot{\mathscr{D}}^{1\!s}_{\eta}$, we conclude $\upsilon_1=
  \upsilon_2$. This proves $\varphi'(\eta):\mathcal{T}_{\eta}(\mathfrak{M})\to X$ is an injection. Next we prove this
  map is also a surjection. To this end, for an arbitrary $v\in X$ we define a function $f:(-\varepsilon,\varepsilon)
  \to\mathfrak{M}$ (for sufficiently small $\varepsilon>0$) as follows:
$$
  f(t)=\varphi^{-1}(\varphi(\eta)+tv), \quad \forall t\in(-\varepsilon,\varepsilon).
$$
  Since $f(0)=\eta$ and clearly the function $\varphi\circ f:(-\varepsilon,\varepsilon)\to X$ is continuously differentiable
  at $t=0$ in the topology of $X$ with derivative $(\varphi\circ f)'(0)=v$, we see that $v=\varphi'(\eta)f'(0)\in\varphi'(\eta)
  \mathcal{T}_{\eta}(\mathfrak{M})$. Hence $\varphi'(\eta)\mathcal{T}_{\eta}(\mathfrak{M})\supseteq X$. This proves the desired
  assertion and completes the proof of Lemma 2.11. $\quad\Box$
\medskip

  {\bf Definition 2.11}\ \ {\em Let $\mathfrak{M}$ be a quasi-differentiable Banach manifold with $C^1$-kernel
  $\mathfrak{M}_0$. We define linear operations in $\mathcal{T}_{\eta}(\mathfrak{M})$ as follows:

  $(1)$\ \ For $\upsilon\in\mathcal{T}_{\eta}(\mathfrak{M})$ and $\lambda\in\mathbb{R}$, let $f:(-\varepsilon,\varepsilon)\to
  \mathfrak{M}$ $(\varepsilon>0)$ be a function such that $f(0)=\eta$, $f(t)$ is continuously differentiable at $t=0$,
  and $f'(0)=\upsilon$. Define $\lambda\upsilon=f_{\lambda}'(0)$, where $f_{\lambda}:(-\varepsilon/\lambda,\varepsilon/\lambda)
  \to\mathfrak{M}$ be the function $f_{\lambda}(t)=f(\lambda t)$.

  $(2)$\ \ For $\upsilon_1,\upsilon_2\in\mathcal{T}_{\eta}(\mathfrak{M})$, arbitrarily choose a $\mathfrak{M}_0$-regular local
  chart $(\mathcal{U},\varphi)$ of $\mathfrak{M}$ at $\eta$. Let $f:(-\varepsilon,\varepsilon)\to\mathfrak{M}$ $(\varepsilon>0$
  sufficiently small$)$ be the function $f(t)=\varphi^{-1}(\varphi(\eta)+t[\varphi'(\eta)\upsilon_1+\varphi'(\eta)\upsilon_2])$
  $($for $|t|<\varepsilon)$. Define $\upsilon_1+\upsilon_2=f'(0)$.}
\medskip

  {\bf Theorem 2.12}\ \ {\em Let $\mathfrak{M}$ be a quasi-differentiable Banach manifold with $C^1$-kernel $\mathfrak{M}_0$.
  Let $X$ be the base space of $\mathfrak{M}$. Then the above definition makes sense, and for any $\eta\in\mathfrak{M}_0$,
  $\mathcal{T}_{\eta}(\mathfrak{M})$ is a Banach space isomorphic to $X$.}
\medskip

  {\em Proof},\ \ The proof is immediate. $\quad\Box$
\medskip

  Let
$$
  \mathcal{T}_{\mathfrak{M}_0}(\mathfrak{M})=\bigcup_{\eta\in\mathfrak{M}_0}\mathcal{T}_{\eta}(\mathfrak{M})\;(\mbox{disjoint union})\;
  :=\{(\eta,\upsilon):\eta\in\mathfrak{M}_0,\upsilon\in\mathcal{T}_{\eta}(\mathfrak{M})\}.
$$
  We can use the standard method to make $\mathcal{T}_{\mathfrak{M}_0}(\mathfrak{M})$ into a (topological) Banach manifold built
  on the Banach space $X_0\times X$. More precisely, for any $(\eta_0,\upsilon_0)\in\mathcal{T}_{\mathfrak{M}_0}(\mathfrak{M})$,
  let $(\mathcal{U},\varphi)$ be a $\mathfrak{M}_0$-regular local chart of $\mathfrak{M}$ at $\eta_0$. Let $\mathcal{U}_0=\mathcal{U}
  \cap\mathfrak{M}_0$. We denote
$$
  \mathcal{O}=\{(\eta,\upsilon):\eta\in\mathcal{U}_0,\upsilon\in\mathcal{T}_{\eta}(\mathfrak{M})\},
$$
  and define $\Phi:\mathcal{O}\to X_0\times X$ as follows:
$$
  \Phi(\eta,\upsilon)=(\varphi(\eta),\varphi'(\eta)\upsilon), \quad \forall (\eta,\upsilon)\in\mathcal{O}.
$$
  We use $(\mathcal{O},\Phi)$ as a local chart of $\mathcal{T}_{\mathfrak{M}_0}(\mathfrak{M})$ at the point $(\eta_0,\upsilon_0)$.
  It is not hard to check that this indeed makes $\mathcal{T}_{\mathfrak{M}_0}(\mathfrak{M})$ into a (topological) Banach manifold.
  We call $\mathcal{T}_{\mathfrak{M}_0}(\mathfrak{M})$ endowed with this topological structure the {\em tangent bundle} of
  $\mathfrak{M}_0$ in $\mathfrak{M}$. Later on we often regard $\mathcal{T}_{\mathfrak{M}_0}(\mathfrak{M})$ as the set
  $\{\upsilon:\upsilon\in\mathcal{T}_{\eta}(\mathfrak{M}),\eta\in\mathcal{U}_0\}$, i.e., we identify the point $(\eta,\upsilon)$
  with $\upsilon$, because for different $\eta,\xi\in\mathfrak{M}_0$, the sets $\mathcal{T}_{\eta}(\mathfrak{M})$,
  $\mathcal{T}_{\xi}(\mathfrak{M})$ do not have common point.
\medskip

\section{Basic concepts on differential equations in quasi-differentiable Banach manifolds}
\setcounter{equation}{0}

\hskip 2em
  Let $\mathfrak{M}$ be a quasi-differentiable Banach manifold and $\mathfrak{M}_0$ be its $C^1$-kernel. Given an open interval
  $I\subseteq {\mathbb{R}}$, we use the notation $C^1\!(I,\mathfrak{M}_0)$ to denote the set of mappings $f:I\to\mathfrak{M}_0$
  such that for any $t_0\in I$, there exists a $\mathfrak{M}_0$-regular local chart $(\mathcal{U},\varphi)$ of $\mathfrak{M}$
  at $f(t_0)$ and $\varepsilon>0$ sufficiently small such that $\varphi\circ f\in C((t_0-\varepsilon,t_0+\varepsilon),X_0)\cap
  C^1((t_0-\varepsilon,t_0+\varepsilon),X)$. Note that this in particular implies $f\in C(I,\mathfrak{M}_0)\subseteq
  C(I,\mathfrak{M})$ and $f(t)$ is continuously differentiable at every point $t\in I$. Let $\mathscr{F}$ be a mapping from
  an open subset $\mathcal{O}$ of $\mathfrak{M}_0$ to $\mathcal{T}_{\mathfrak{M}_0}(\mathfrak{M})$, such that for any $\eta\in
  \mathfrak{M}_0$, $\mathscr{F}(\eta)\in\mathcal{T}_{\eta}(\mathfrak{M})$. We call $\mathscr{F}$ a {\em vector field in
  $\mathfrak{M}$ with domain $\mathcal{O}$}. The purpose of this and the next sections is to study the following differential
  equation in the quasi-differentiable Banach manifold $\mathfrak{M}$:
\begin{equation}
   \eta'=\mathscr{F}(\eta),
\end{equation}
%---(3.1)---
  or more precisely, the initial value problem of the above equation:
\begin{equation}
\left\{
\begin{array}{rlll}
   \eta'(t)&=&\mathscr{F}(\eta(t)), &\quad t>0,\\
   \eta(0)&=&\eta_0, &
\end{array}
\right.
\end{equation}
%---(3.2)---
  where $\eta_0$ is a given point in $\overline{\mathcal{O}}$=the closure of $\mathcal{O}$ in $\mathfrak{M}$. By a {\em solution}
  of the equation (3.1) in an open interval $I\subseteq{\mathbb{R}}$ we mean a function $\eta\in C^1\!(I,\mathfrak{M}_0)$ such
  that $\eta(t)\in\mathcal{O}$ and $\eta'(t)=\mathscr{F}(\eta(t))$ for $t\in I$, and by a {\em solution} of the problem (3.2)
  in a interval $[0,T)$ ($0<T\leqslant\infty$) we mean a function $\eta\in C([0,T),\mathfrak{M})\cap C^1((0,T),\mathfrak{M}_0)$
  such that $\eta(t)\in\mathcal{O}$ and $\eta'(t)=\mathscr{F}(\eta(t))$ for $t\in (0,T)$, and $\eta(0)=\eta_0$.

  Recall that a differential equation $x'=F(x)$ in  a Banach space $X$, where $F\in C^1(O,X)$ with $O$ being an open subset
  of an embedded Banach subspace $X_0$ of $X$, is said to be of {\em parabolic} type if for any $x\in O$, $F'(x)$ is a
  sectorial operator in $X$ with domain $X_0$, and the graph norm of $X_0={\rm Dom}\,F'(x)$ is equivalent to the norm of $X_0$,
  i.e., there exist positive constants $C_1,C_2$ such that
$$
  C_1\|y\|_{X_0}\leqslant\|y\|_{X}+\|F'(x)y\|_{X}\leqslant C_2\|y\|_{X_0}, \quad \forall y\in X_0.
$$
  Let $\mathscr{F}$ be as before. For $\eta_0\in\mathcal{O}\cap\mathfrak{M}_0$, let $(\mathcal{U},\varphi)$ be a
  $\mathfrak{M}_0$-regular local chart of $\mathfrak{M}$ at $\eta_0$ and let $\mathcal{U}_0=\mathcal{U}\cap\mathfrak{M}_0$. Then
  $(\mathcal{U}_0,\varphi|_{\mathcal{U}_0})$ is a local chart of $\mathfrak{M}_0$ at $\eta_0$. The {\em representation} of the
  vector field $\mathscr{F}$ in the local chart $(\mathcal{U},\varphi)$ is the mapping $F:\varphi(\mathcal{O}\cap\mathcal{U}_0)
  \subseteq X_0\to X$ defined as follows:
$$
  F(x)=\varphi'(\varphi^{-1}(x))\mathscr{F}(\varphi^{-1}(x)), \quad x\in\varphi(\mathcal{O}\cap\mathcal{U}_0).
$$
  The representation of the differential equation (3.1) in the local chart $(\mathcal{U},\varphi)$ is the following
  differential equation in $X$:
\begin{equation}
  x'=F(x).
\end{equation}
%---(2.7)---
  We introduce the following concept:

\medskip
  {\bf Definition 3.1}\ \ {\em Let $(\mathfrak{M},\mathscr{A})$ be a quasi-differentiable Banach manifold with a $C^1$-kernel
  $\mathfrak{M}_0$. Let $\mathscr{F}$ be a vector field in $\mathfrak{M}$ with domain $\mathcal{O}\subseteq\mathfrak{M}_0$,
  where $\mathcal{O}$ is an open subset of $\mathfrak{M}_0$. We have the following concepts:

  $(1)$\ We say the differential equation $(3.1)$ is of \mbox{\small\boldmath 4parabolic4} type at a point $\eta\in\mathcal{O}$
  if there exists a local chart $(\mathcal{U},\varphi)\in\mathscr{A}$ with $\eta\in\mathcal{U}$ such that the representation of
  the equation $(3.1)$ in this local chart is of parabolic type. We say $(3.1)$ is of parabolic type in $\mathcal{O}$ if for
  any $\eta\in\mathcal{O}$, $(3.1)$ is of parabolic at $\eta$.

  $(2)$\ Let $k\in\mathbb{N}$. We say the vector field $\mathscr{F}$ is of \mbox{\small\boldmath $C^k$-$class$} $($resp.
  \mbox{\small\boldmath $C^{k-0}$-$class$}$)$ at a point $\eta\in\mathcal{O}$ if there exists a local chart
  $(\mathcal{U},\varphi)\in\mathscr{A}$ with $\eta\in\mathcal{U}$ such that the representation of $\mathscr{F}$ in this local
  chart is of $C^k$-class $($resp. $C^{k-0}$-class$)$. We say $\mathscr{F}$ is of $C^k$-class $($resp. $C^{k-0}$-class$)$ in
  $\mathcal{O}$ if for any $\eta\in\mathcal{O}$, $\mathscr{F}$ is of $C^k$-class $($resp. $C^{k-0}$-class$)$ at $\eta$.

  $(3)$\ Let $k\in\mathbb{N}$. We say the differential equation $(3.1)$ is of \mbox{\small\boldmath $C^k$-$class$} $($resp.
  \mbox{\small\boldmath $C^{k-0}$-$class$}$)$ \mbox{\small\boldmath $parabolic$} type at a point $\eta\in\mathcal{O}$ if
  there exists a local chart $(\mathcal{U},\varphi)\in\mathscr{A}$ with $\eta\in\mathcal{U}$ such that the representation of
  the equation $(3.1)$ in this local chart is of parabolic type and also the representation of $\mathscr{F}$ in this local
  chart is of $C^k$-class $($resp. $C^{k-0}$-class$)$. We say $(3.1)$ is of $C^k$-class $($resp. $C^{k-0}$-class$)$ parabolic
  type in $\mathcal{O}$ if for any $\eta\in\mathcal{O}$, $(3.1)$ is of $C^k$-class $($resp. $C^{k-0}$-class$)$ parabolic at
  $\eta$.

  $(4)$\ A point $\eta\in\mathcal{O}$ is called a \mbox{\small\boldmath $stationary$ $point$} of the differential equation
  $(3.1)$ if $\mathscr{F}(\eta)=0$. In this case $\eta$ is also called a \mbox{\small\boldmath $zero$} of the vector field
  $\mathscr{F}$.}

\medskip
  Note that we do not exclude the possibility that the differential equation $(3.1)$ is of parabolic type
  in one local chart $(\mathcal{U},\varphi)\in\mathscr{A}$ with $\eta\in\mathcal{U}$ but is not of parabolic type in another
  local chart $(\mathcal{V},\psi)\in\mathscr{A}$ with $\eta\in\mathcal{V}$. Similarly, we do not exclude the possibility that
  the vector field $\mathscr{F}$ is of $C^k$-class $($resp. $C^{k-0}$-class$)$ in one local chart $(\mathcal{U},\varphi)\in
  \mathfrak{M}$ with $\eta\in\mathcal{U}$ but is not of $C^k$-class $($resp. $C^{k-0}$-class$)$ in another local chart
  $(\mathcal{V},\psi)\in\mathscr{A}$ with $\eta\in\mathcal{V}$. These drawbacks of analysis in quasi-differentiable
  Banach manifolds are not important because when we deal with a differential equation in a such manifold we usually always
  choose a local chart in which the representation of this equation has the best properties, and the other local charts in
  which these properties do not hold will not be considered. For instance, by choosing a local chart in which the equation (3.1)
  is of $C^{2-0}$-class parabolic type, transforming the problem (3.2) into corresponding initial value problem in a Banach
  space, and then applying local theory of parabolic differential equations in Banach spaces, we have the following basic result:

\medskip
  {\bf Theorem 3.2}\ \  {\em Assume that the differential equation $(3.1)$ is of $C^{2-0}$-class parabolic type. Then for any
  $\eta_0\in\overline{\mathcal{O}}$ there exists a corresponding $\delta>0$ such that the problem $(3.2)$ has a unique solution
  $\eta\in C([0,\delta),\mathfrak{M})\cap C^1((0,\delta),\mathfrak{M}_0)$.} $\quad\Box$

\medskip
  To extend the theory of asymptotic stability of stationary solution of parabolic differential equations in Banach spaces
  to such equations in quasi-differentiable Banach manifolds, we need to make some more preparations.

\medskip
  {\bf Lemma 3.3}\ \ {\em Let $(\mathfrak{M},\mathscr{A})$ be a quasi-differentiable Banach manifold with a $C^1$-kernel
  $\mathfrak{M}_0$ and an inner $C^2$-kernel $\mathfrak{M}_1$. Let $\mathscr{F}$ be a vector field in $\mathfrak{M}$ with
  domain $\mathcal{O}\subseteq\mathfrak{M}_0$. Let $(\mathcal{U},\varphi),(\mathcal{V},\psi)\in\mathscr{A}$ be two local
  charts of $\mathfrak{M}$ such that $\mathcal{O}\cap\mathcal{U}\cap\mathcal{V}\neq\emptyset$. Let $\Phi=\varphi\circ\psi^{-1}:
  \psi(\mathcal{U}_0\cap\mathcal{V}_0)\to\varphi(\mathcal{U}_0\cap\mathcal{V}_0)$ and $\Psi=(\psi\circ\varphi^{-1})'\circ\Phi$.
  Let $F_U$ and $F_V$ be representations of $\mathscr{F}$ in the local charts $(\mathcal{U},\varphi)$ and $(\mathcal{V},\psi)$,
  respectively, and assume that they are both differentiable as maps from open subsets of $X_0$ to $X$. Then the following
  relation holds:}
\begin{equation}
   F_V'(y)z=\Psi(y)F_U'(\Phi(y))\Psi(y)^{-1}z+\{\Psi'(y)z\}F_U(\Phi(y)),  \quad
   y\in\psi(\mathcal{O}\cap\mathcal{U}_1\cap\mathcal{V}_1), \;\; z\in X_0.
\end{equation}
%---(3.4)---

  {\em Proof}.\ \ Let $\mathcal{U}_0=\mathcal{U}\cap\mathfrak{M}_0$, $\mathcal{V}_0=\mathcal{V}\cap\mathfrak{M}_0$. By
  definition, we have
$$
  F_U(x)=\varphi'(\varphi^{-1}(x))\mathscr{F}(\varphi^{-1}(x)), \quad x\in\varphi(\mathcal{O}\cap\mathcal{U}_0),
$$
$$
  F_V(y)=\psi'(\psi^{-1}(y))\mathscr{F}(\psi^{-1}(y)), \quad y\in\psi(\mathcal{O}\cap\mathcal{V}_0).
$$
  From the definition of $\varphi'(\eta)$ and $\psi'(\eta)$ and the relation (2.2) we see that the following relation holds:
$$
  \psi'(\eta)=(\psi\circ\varphi^{-1})'(\varphi(\eta))\varphi'(\eta), \quad \forall\eta\in\mathcal{U}_0\cap\mathcal{V}_0.
$$
  Hence
\begin{equation}
  F_V(y)=\Psi(y)F_U(\Phi(y)), \quad y\in\psi(\mathcal{O}\cap\mathcal{U}_0\cap\mathcal{V}_0).
\end{equation}
%---(3.5)---
  The condition of this lemma implies that $[y\mapsto\Psi(y)],[y\mapsto\Psi(y)^{-1}]\in C(\psi(\mathcal{U}_0\cap
  \mathcal{V}_0),L(X))\cap C(\psi(\mathcal{U}_1\cap\mathcal{V}_1),L(X_0))$ and $\Phi\in
  \mathfrak{C}^1(\psi(\mathcal{U}_0\cap\mathcal{V}_0);X,X))\cap\mathfrak{C}^2(\psi(\mathcal{U}_1\cap\mathcal{V}_1);X,X))
  \cap\mathfrak{C}^1(\psi(\mathcal{U}_1\cap\mathcal{V}_1);X_0,X_0))$. Since $\Phi'(y)=\Psi(y)^{-1}$ (by Lemma 2.3) and
$$
  \Psi'(y)=(\psi\circ\varphi^{-1})''(\Phi(y))\Phi'(y)=(\psi\circ\varphi^{-1})''(\Phi(y))\Psi(y)^{-1},
$$
  we see $\Psi'\in C(\psi(\mathcal{U}_1\cap\mathcal{V}_1),L(X,L(X)))$. Hence, by differentiating (3.5) in $X_1$, we see that
  (3.4) holds for $z\in X_1$. By density of $X_1$ in $X_0$, we concludes that (3.4) also holds for $z\in X_1$. This proves
  the lemma. $\quad\Box$

\medskip
  {\bf Corollary 3.4}\ \ {\em Let $(\mathfrak{M},\mathscr{A})$ be a quasi-differentiable Banach manifold with a $C^1$-kernel
  $\mathfrak{M}_0$ and an inner $C^2$-kernel $\mathfrak{M}_1$. Let $\mathscr{F}$ be a vector field in $\mathfrak{M}$ with
  domain $\mathcal{O}\subseteq\mathfrak{M}_0$. Let $\eta\in\mathcal{O}\cap\mathfrak{M}_1$. If the representation of the
  equation $(3.1)$ in one local chart $(\mathcal{U},\varphi)\in\mathscr{A}$ with $\eta\in\mathcal{U}$ is of parabolic type at
  $\eta$, then in any other local chart $(\mathcal{V},\psi)\in\mathscr{A}$ with $\eta\in\mathcal{V}$, the representation of the
  equation $(3.1)$ is also of parabolic type at $\eta$.}

\medskip
  {\em Proof}.\ \ It is easy to see that if $F_U'(x)$ is a sectorial operator in $X$ for some $x\in\varphi(\mathcal{O}\cap\mathcal{U}_0)$ then
  $\Psi(y)F_U'(\Phi(y))\Psi(y)^{-1}$ is also a sectorial operator in $X$ for $y=\Phi^{-1}(x)\in\psi(\mathcal{O}\cap\mathcal{U}_0
  \cap\mathcal{V}_0)$, as $[y\mapsto\Psi(y)]\in C(\psi(\mathcal{U}_0\cap\mathcal{V}_0),L(X))$ and $[y\mapsto\Psi(y)^{-1}]\in
  C(\psi(\mathcal{U}_1\cap\mathcal{V}_1),L(X_0))$. Moreover, it is clear that the operator $z\mapsto
  \{\Psi'(y)z\}F_U(\Phi(y))$ is a bounded linear operator in $X$ for fixed $y\in\psi(\mathcal{O}\cap\mathcal{U}_1\cap\mathcal{V}_1)$.
  Hence, by a standard perturbation theorem for sectorial operators, it follows that $F_V'(y)$ is also a sectorial operator in $X$.
  This proves the lemma. $\quad\Box$

\medskip
  {\bf Corollary 3.5}\ \ {\em Let $(\mathfrak{M},\mathscr{A})$ be a quasi-differentiable Banach manifold with a $C^1$-kernel
  $\mathfrak{M}_0$ and an inner $C^2$-kernel $\mathfrak{M}_1$. Let $\mathscr{F}$ be a vector field in $\mathfrak{M}$ with domain
  $\mathcal{O}\subseteq\mathfrak{M}_0$. Let $\eta^*\in\mathcal{O}\cap\mathfrak{M}_1$ be a zero of $\mathscr{F}$. Then for any
  two local charts $(\mathcal{U},\varphi),(\mathcal{V},\psi)\in\mathscr{A}$ with $\eta^*\in\mathcal{U}\cap\mathcal{V}$, by letting
  $F_U(x)$ and $F_V(y)$ be the representations of $\mathscr{F}$ in these local charts, respectively, the following relation holds:
\begin{equation}
   F_V'(y^*)=\Psi(y^*)F_U'(x^*)\Psi(y^*)^{-1},
\end{equation}
%---(3.6)---
  where $\Psi(y)=(\psi\circ\varphi^{-1})'(\Phi(y))$ with $\Phi(y)=\varphi(\psi^{-1}(y))$, $x^*=\varphi(\eta^*)$, and $y^*=
  \psi(\eta^*)$.}

\medskip
  {\em Proof}.\ \ This is an immediate consequence of (3.4) and the fact that $F_U(x^*)=0$. $\quad\Box$

\medskip
  {\bf Lemma 3.6}\ \ {\em Let $(\mathfrak{M},\mathscr{A})$ be a quasi-differentiable Banach manifold with a $C^1$-kernel
  $\mathfrak{M}_0$ and an inner $C^1$-kernel $\mathfrak{M}_1$. Let $X,X_0,X_1$ be the base spaces of $\mathfrak{M}$,
  $\mathfrak{M}_0$ and $\mathfrak{M}_1$, respectively. Let $\eta\in\mathfrak{M}_1$. For any $(\mathcal{U},\varphi)\in
  \mathscr{A}$ such that $\eta\in\mathcal{U}$, we define a nonnegative-valued function $\|\cdot\|_{\varphi,X_0}$ in
  $\mathcal{T}_{\eta}^{00}(\mathfrak{M})$ as follows:
$$
  \|\upsilon\|_{\varphi,X_0}=\|\varphi'(\eta)\upsilon\|_{X_0}, \quad \forall\upsilon\in\mathcal{T}_{\eta}^{00}(\mathfrak{M}),
$$
  Then for any two local charts $(\mathcal{U},\varphi),(\mathcal{V},\psi)\in\mathscr{A}$ such that $\eta\in\mathcal{U}\cap
  \mathcal{V}$, $\|\cdot\|_{\varphi,X_0}$ and $\|\cdot\|_{\psi,X_0}$ are equivalent, i.e., there exists constants $C_1,C_2>0$
  such that}
\begin{equation}
  \|\upsilon\|_{\psi,X_0}\leqslant C_1\|\upsilon\|_{\varphi,X_0}, \quad
  \|\upsilon\|_{\varphi,X_0}\leqslant C_2\|\upsilon\|_{\psi,X_0}, \quad
   \forall\upsilon\in\mathcal{T}_{\eta}^{00}(\mathfrak{M}).
\end{equation}
%---(3.7)---

  {\em Proof}.\ \ From the proof of the assertion (2) of Lemma 2.10 we see that the following relation holds:
$$
  \psi'(\eta)\upsilon=(\psi\circ\varphi^{-1})'(\varphi(\eta))\varphi'(\eta)\upsilon, \quad
   \forall\upsilon\in\mathcal{T}_{\eta}^{00}(\mathfrak{M}).
$$
  Since $\eta\in\mathfrak{M}_1$ implies that $(\psi\circ\varphi^{-1})'(\varphi(\eta))\in L(X_0)$, it follows that by
  letting $C_1=\|(\psi\circ\varphi^{-1})'(\varphi(\eta))\|_{L(X_0)}$, we immediately obtain the first inequality in
  (3.7). Similarly, by letting $C_2=\|(\varphi\circ\psi^{-1})'(\varphi(\eta))\|_{L(X_0)}$, we obtain the second
  inequality in (3.7). $\quad\Box$
\medskip

  Let $\mathfrak{M}$ be a quasi-differentiable Banach manifold with $C^1$-kernel $\mathfrak{M}_0$. Let $X$ and $X_0$ be the base
  spaces of $\mathfrak{M}$ and $\mathfrak{M}_0$, respectively. For $\eta\in\mathfrak{M}_0$ and a $\mathfrak{M}_0$-regular
  local chart $(\mathcal{U},\varphi)$ of $\mathfrak{M}$ at $\eta$, we denote
$$
  \mathcal{T}_{\eta}^0(\mathfrak{M})=\{f'(0):f:(-\varepsilon,\varepsilon)\to\mathfrak{M}_0,\; f(0)=\eta,\;
  \mbox{$f(t)$ is strongly continuously differentiable at $t=0$}\}.
$$
$$
  \mathcal{T}_{\eta,\varphi}^0(\mathfrak{M})=\{f'(0):f:(-\varepsilon,\varepsilon)\to\mathfrak{M}_0,\; f(0)=\eta,\;
  \mbox{$f(t)$ is $\varphi$-strongly continuously differentiable at $t=0$}\}.
$$
  It is clear that the following relations hold:
$$
  \mathcal{T}_{\eta}^0(\mathfrak{M})=\displaystyle\bigcup_{(\mathcal{U},\varphi)}\mathcal{T}_{\eta,\varphi}^0(\mathfrak{M}),
  \qquad
  \mathcal{T}_{\eta}^{00}(\mathfrak{M})\subseteq\bigcap_{(\mathcal{U},\varphi)}\mathcal{T}_{\eta,\varphi}^0(\mathfrak{M}),
  \qquad
  \varphi'(\eta)\mathcal{T}_{\eta,\varphi}^0(\mathfrak{M})=X_0.
%  \qquad
%  X_1\subseteq\displaystyle\bigcap_{(\mathcal{U},\varphi)}\varphi'(\eta)\mathcal{T}_{\eta}^{00}(\mathfrak{M}),
$$

  {\bf Lemma 3.7}\ \ {\em Let $(\mathfrak{M},\mathscr{A})$ be a quasi-differentiable Banach manifold with a $C^1$-kernel
  $\mathfrak{M}_0$ and an inner $C^1$-kernel $\mathfrak{M}_1$. Let $X,X_0,X_1$ be the base spaces of $\mathfrak{M}$,
  $\mathfrak{M}_0$ and $\mathfrak{M}_1$, respectively. Let $\eta\in\mathfrak{M}_1$. Then for any two local charts
  $(\mathcal{U},\varphi),(\mathcal{V},\psi)\in\mathscr{A}$ such that $\eta\in\mathcal{U}\cap\mathcal{V}$, we have
  $\mathcal{T}_{\eta,\varphi}^0(\mathfrak{M})=\mathcal{T}_{\eta,\psi}^0(\mathfrak{M})$, so that for any local chart
  $(\mathcal{U},\varphi)\in\mathscr{A}$ such that $\eta\in\mathcal{U}$, we have}
\begin{equation}
  \mathcal{T}_{\eta}^0(\mathfrak{M})=\mathcal{T}_{\eta,\varphi}^0(\mathfrak{M}) \quad \mbox{and} \quad
  \varphi'(\eta)\mathcal{T}_{\eta}^0(\mathfrak{M})=X_0.
\end{equation}
%---(3.8)---

  {\em Proof}.\ \ Let $\upsilon\in\mathcal{T}_{\eta,\varphi}^0(\mathfrak{M})$. Then there exists a function $f:
  (-\varepsilon,\varepsilon)\to\mathfrak{M}_0$ such that $f(0)=\eta$, $f(t)$ is $\varphi$-strongly continuously differentiable
  at $t=0$, such that $\upsilon=f'(0)$. We define a function $g:(-\varepsilon,\varepsilon)\to\mathfrak{M}_0$ as follows:
$$
  g(t)=\psi^{-1}(\psi(\eta)+t(\psi\circ\varphi^{-1})'(\varphi(\eta))\varphi'(\eta)\upsilon), \quad
  \forall t\in(-\varepsilon,\varepsilon).
$$
  Since the condition $\eta\in\mathfrak{M}_1$ implies that $(\psi\circ\varphi^{-1})'(\varphi(\eta))\in L(X_0)$, we see that
  the above definition makes sense, $g(0)=\eta$, $g(t)$ is $\psi$-strongly continuously differentiable at $t=0$, and
  $(\psi\circ g)'(0)=(\psi\circ\varphi^{-1})'(\varphi(\eta))\varphi'(\eta)\upsilon$. This implies that $g'(0)\in
  \mathcal{T}_{\eta,\psi}^0(\mathfrak{M})$ and $g'(0)=\upsilon$. Hence $\upsilon\in\mathcal{T}_{\eta,\psi}^0(\mathfrak{M})$.
  This proves $\mathcal{T}_{\eta,\varphi}^0(\mathfrak{M})\subseteq\mathcal{T}_{\eta,\psi}^0(\mathfrak{M})$. Similarly we can
  prove $\mathcal{T}_{\eta,\psi}^0(\mathfrak{M})\subseteq\mathcal{T}_{\eta,\varphi}^0(\mathfrak{M})$. Hence
  $\mathcal{T}_{\eta,\varphi}^0(\mathfrak{M})=\mathcal{T}_{\eta,\psi}^0(\mathfrak{M})$. The other relation is an immediate
  consequence of this one. $\quad\Box$
\medskip

  Since $\mathcal{T}_{\eta,\varphi}^0(\mathfrak{M})$ is clearly a linear subspace of $\mathcal{T}_{\eta}(\mathfrak{M})$, it
  follows that if $\eta\in\mathfrak{M}_1$ then $\mathcal{T}_{\eta}^0(\mathfrak{M})$ is a linear subspace of
  $\mathcal{T}_{\eta}(\mathfrak{M})$. As a consequence of the above lemma, we have

\medskip
  {\bf Corollary 3.8}\ \ {\em Let $(\mathfrak{M},\mathscr{A})$ be a quasi-differentiable Banach manifold with a $C^1$-kernel
  $\mathfrak{M}_0$ and an inner $C^2$-kernel $\mathfrak{M}_1$. Let $X,X_0,X_1$ be the base spaces of $\mathfrak{M}$,
  $\mathfrak{M}_0$ and $\mathfrak{M}_1$, respectively. Let $\eta\in\mathfrak{M}_1$. Then for any local chart
  $(\mathcal{U},\varphi)\in\mathscr{A}$ such that $\eta\in\mathcal{U}$, when endowed with the norm
  $\|\cdot\|_{\varphi,X_0}$, $\mathcal{T}_{\eta}^0(\mathfrak{M})=\mathcal{T}_{\eta,\varphi}^0(\mathfrak{M})$ is an embedded
  Banach subspace of $\mathcal{T}_{\eta}(\mathfrak{M})$ isomorphic to $X_0$.} $\quad\Box$

\medskip
  The above discussion ensures that the following definition makes sense:

\medskip
  {\bf Definition 3.9}\ \ {\em Let notations and assumptions be as in Corollary 3.5; in particular, let $\eta\in\mathcal{O}
  \cap\mathfrak{M}_1$ be a zero of the vector field $\mathscr{F}$. Assume further that the equation $(3.1)$ is of $C^2$-class
  parabolic type at $\eta$. We have the following concepts:

  $(1)$\ We define $\mathscr{F}'(\eta)$ to be the following linear operator from $\mathcal{T}_{\eta}^0(\mathfrak{M})$
  to $\mathcal{T}_{\eta}(\mathfrak{M})$: For any local chart $(\mathcal{U},\varphi)\in\mathscr{A}$ with $\eta\in\mathcal{U}$,
$$
  \mathscr{F}'(\eta)\upsilon=\varphi'(\eta)^{-1}F_U'(\varphi(\eta))\varphi'(\eta)\upsilon, \quad
  \forall\upsilon\in\mathcal{T}_{\eta}^0(\mathfrak{M})=\mathcal{T}_{\eta,\varphi}^0(\mathfrak{M}),
$$
  where $F_U$ is the representation of $\mathscr{F}$ in the local chart $(\mathcal{U},\varphi)$.

  $(2)$\ We define the \mbox{\small\boldmath $spectrum$} of $\mathscr{F}'(\eta)$ as follows:
$$
  \sigma(\mathscr{F}'(\eta))=\sigma(F_U'(\varphi(\eta))).
$$
  It is the same with that by regarding $\mathscr{F}'(\eta)$ as an unbounded linear operator in $\mathcal{T}_{\eta}(\mathfrak{M})$
  with domain $\mathcal{T}_{\eta}^0(\mathfrak{M})$.

  $(3)$\ We say $\mathscr{F}'(\eta)$ is a \mbox{\small\boldmath $standard$ $Fredholm$ $operator$} if
  ${\rm dim}\,{\rm Ker}\mathscr{F}'(\eta)<\infty$, ${\rm Range}\mathscr{F}'(\eta))$ is closed, and}
$$
  \mathcal{T}_{\eta}(\mathfrak{M})={\rm Ker}\mathscr{F}'(\eta)\oplus\,{\rm Range}\mathscr{F}'(\eta).
$$
  Note that this particularly implies that $\mathscr{F}'(\eta)$ is a Fredholm operator of index zero.

\medskip
  {\bf Definition 3.10}\ \ {\em Let $(\mathfrak{M},\mathscr{A})$ be a quasi-differentiable Banach manifold with a $C^1$-kernel
  $\mathfrak{M}_0$. Let $\mathscr{F}$ be a vector field in $\mathfrak{M}$ with domain $\mathcal{O}\subseteq\mathfrak{M}_0$. Let
  $\eta\in\mathcal{O}\cap\mathfrak{M}_0$ be a zero of $\mathscr{F}$, so that $\eta$ is a stationary point of the equation
  $(3.1)$. We say $\eta$ is \mbox{\small\boldmath $stable$} if for any neighborhood $\mathscr{B}$ of $\eta$ in $\mathfrak{M}_0$,
  there exists corresponding neighborhood $\mathscr{B}'\subseteq\mathcal{O}$ of $\eta$, such that for any $\eta_0\in
  \mathscr{B}'$, the solution $\eta(t)$ of the problem $(3.2)$ exists for all $t\geqslant 0$ and $\eta(t)\in\mathscr{B}$ for
  all $t\geqslant 0$. We say $\eta$ is \mbox{\small\boldmath $asymptotically$ $stable$} if it is stable and, furthermore, for
  any $\eta_0\in\mathscr{B}'$ the solution $\eta(t)$ of the problem $(3.2)$ satisfies $\displaystyle\lim_{t\to\infty}\eta(t)=\eta$.}

\medskip
  {\bf Theorem 3.11}\ \ {\em Let notations and assumptions be as in Corollary 3.5; in particular, let $\eta\in\mathcal{O}
  \cap\mathfrak{M}_1$ be a zero of the vector field $\mathscr{F}$, so that $\eta$ a stationary point of the equation $(3.1)$.
  Assume further that the equation $(3.1)$ is of $C^2$-class parabolic type at $\eta$. We have the following assertions:

  $(1)$\ If $\sup\{{\rm Re}\lambda:\lambda\in\sigma(\mathscr{F}'(\eta))\}<0$ then the stationary point $\eta$ is asymptotically
  stable.

  $(2)$\ If $\sup\{{\rm Re}\lambda:\lambda\in\sigma(\mathscr{F}'(\eta))\}>0$ then the stationary point $\eta$ is unstable.} $\quad\Box$

\section{Invariant and quasi-invariant differential equations in Banach manifolds}
\setcounter{equation}{0}

\hskip 2em
  In this section we study the structure of local phase diagram of an invariant or quasi-invariant parabolic differential
  equation in a Banach manifold near its center manifold.

  Let $\mathfrak{M}$ be a (topological or $C^0$) Banach manifold and $G$ a Lie group of dimension $n$. An {\em action} of
  $G$ to $\mathfrak{M}$ is a mapping $p:G\times\mathfrak{M}\to\mathfrak{M}$ satisfying the following two conditions:
\begin{enumerate}
\item[]$(L1)$\ \ $p\in C(G\times\mathfrak{M},\mathfrak{M})$.\vspace*{-0.2cm}
\item[]$(L2)$\ \ $p(e,\eta)=\eta$, $\forall\eta\in\mathfrak{M}$, where $e$ denotes the unit element of $G$, and
$$
  p(a,p(b,\eta))=p(ab,\eta), \quad \forall a,b\in G, \;\; \forall\eta\in\mathfrak{M}.
$$
%\item[]$(L3)$\ \ If $a,b\in G$ are such that $p(a,\eta)=p(b,\eta)$ for some $\eta\in\mathfrak{M}$, then $a=b$.\vspace*{-0.2cm}
\end{enumerate}
  Assume further that $\mathfrak{M}$ has a $C^k$-embedded Banach submanifold $\mathfrak{M}_0$, $k\geqslant 1$. We say that
  the Lie group action $(G,p)$ to $\mathfrak{M}$ is {\em $\mathfrak{M}_0$-regular} if the following three additional conditions
  are also satisfied:
\begin{enumerate}
\item[]$(L3)$\ \ $p(a,\mathfrak{M}_0)\subseteq\mathfrak{M}_0$, $\forall a\in G$, and $p\in
  C(G\times\mathfrak{M}_0,\mathfrak{M}_0)$.\vspace*{-0.2cm}
\item[]$(L4)$\ \ For any $a\in G$, the mapping $\eta\mapsto p(a,\eta)$ is differentiable at every point $\eta\in\mathfrak{M}_0$
  ($\Rightarrow\partial_{\eta}p(a,\eta)\in L(T_{\eta}(\mathfrak{M}),T_{p(a,\eta)}(\mathfrak{M}))$, $\forall\eta\in\mathfrak{M}_0$),
  and $[(a,(\eta,\upsilon))\mapsto (p(a,\eta),\partial_{\eta}p(a,\eta)\upsilon)]\in C(G\times
  \mathcal{T}_{\mathfrak{M}_0}(\mathfrak{M}),\mathcal{T}_{\mathfrak{M}_0}(\mathfrak{M}))$.\vspace*{-0.2cm}
\item[]$(L5)$\ \ For any $\eta\in\mathfrak{M}_0$, the mapping $a\mapsto p(a,\eta)$ is differentiable at every point $a\in G$
  ($\Rightarrow\partial_ap(a,\eta)\in L(T_{a}(G),T_{p(a,\eta)}(\mathfrak{M}))$, $\forall a\in G$),
  $[((a,z),\eta)\mapsto(p(a,\eta),\partial_ap(a,\eta)z)]\in C(T(G)\times\mathfrak{M}_0,
  \mathcal{T}_{\mathfrak{M}_0}(\mathfrak{M}))$, and
$$
  {\rm rank}\,\partial_{a}p(a,\eta)=n, \quad \forall a\in G, \;\; \forall\eta\in\mathfrak{M}_0.
$$
\end{enumerate}

  {\bf Definition 4.1}\ \ {\em Let $(\mathfrak{M},\mathscr{A})$ be a Banach manifold and $\mathfrak{M}_0$ a $C^k$-embedded Banach
  submanifold of $\mathfrak{M}$, $k\geqslant 1$. Let $\mathscr{F}$ be a vector field in $\mathfrak{M}$ with domain $\mathfrak{M}_0$.
  Let $(G,p)$ be a $\mathfrak{M}_0$-regular Lie group action to $\mathfrak{M}$. We say $\mathscr{F}$ is \mbox{\small\boldmath
  $quasi$-$invariant$} under the Lie group action $(G,p)$ if there exists a positive-valued function $\theta$ defined in $G$
  such that the following relation holds:
$$
  \mathscr{F}(p(a,\eta))=\theta(a)\partial_{\eta}p(a,\eta)\mathscr{F}(\eta), \quad \forall a\in G, \;\; \forall\eta\in\mathfrak{M}_0.
$$
  In this case we also say $\mathscr{F}$ is \mbox{\small\boldmath $\theta$-$quasi$-$invariant$} and call $\theta$ \mbox{\small\boldmath
  $quasi$-$invariance$ $factor$}, and also say the differential equation $(2.2)$ is $\theta$-quasi-invariant under the Lie group
  action $(G,p)$. If in particular $\theta(a)=1$, $\forall a\in G$, then we simply say the vector field $\mathscr{F}$ and the
  differential equation $(2.2)$ are \mbox{\small\boldmath $invariant$} under the Lie group action $(G,p)$.}
\medskip

  Clearly, $\theta$ is a group homomorphism from $G$ to the multiplicative group ${\mathbb{R}}_+$: For any $a,b\in G$, $\theta(ab)
  =\theta(a)\theta(b)$. The following result is obvious:
\medskip

  {\bf Lemma 4.2}\ \ {\em Let $\mathscr{F}$ be a vector field in $\mathfrak{M}$ with domain $\mathfrak{M}_0$. Assume that
  $\mathscr{F}$ is $\theta$-quasi-invariant under the Lie group action $(G,p)$. Then we have the following assertion: If
  $t\mapsto\eta(t)$ is a solution of the differential equation $(3.1)$, then for any $a\in G$, $t\mapsto p(a,\eta(t\theta(a)))$
  is also a solution of this equation, and if $\eta^*\in\mathfrak{M}_0$ is a stationary point of $(3.1)$, then for any $a\in G$,
  $p(a,\eta^*)$ is also a stationary point of $(3.1)$.} $\quad\Box$

\medskip
  As a consequence, no stationary point of a quasi-invariant differential equation is isolated, and if $\eta^*\in\mathfrak{M}_0$
  is a stationary point of (3.1), then all points in the $n$-dimensional manifold $\{p(a,\eta^*):a\in G\}$ are stationary points
  of (3.1).

  Inspired by potential application to various evolutionary type free boundary problems, we consider a general situation where
  the vector field $\mathscr{F}$ is quasi-invariant under a finite number of Lie group actions $(G_i,p_i)$, $i=1,2,\cdots,N$,
  with possibly different quasi-invariance factors $\theta_i$, $i=1,2,\cdots,N$, respectively. We assume that the combined
  action of these Lie group actions satisfies the following additional conditions:
\begin{enumerate}
\item[]$(L6)$\ \ For any $1\leqslant i,j\leqslant N$ with $i\neq j$ there exists corresponding smooth function $f_{ij}:G_i\times G_j\to G_i$ such
  that
$$
   p_j(b,p_i(a,\eta))=p_i(f_{ij}(a,b),p_j(b,\eta)), \quad \forall\eta\in\mathfrak{M},\;\; \forall a\in G_i,\;\; \forall b\in G_j.
$$
  and for every fixed $b\in G_j$, the mapping $a\mapsto f_{ij}(a,b)$ is bijectiove with a smooth inverse.
\item[]$(L7)$\ \ The Lie group actions $(G_1,p_1)$, $(G_2,p_2)$, $\cdots$, $(G_N,p_N)$ are fully ranked, i.e., by letting
  $g:G\times\mathfrak{M}:=(G_1\times G_2\times\cdots\times G_N)\times\mathfrak{M}\to\mathfrak{M}$ be the function
$$
   g(a,\eta)=p_1(a_1,p_2(a_2,\cdots,p_N(a_N,\eta)))
$$
  for $\eta\in\mathfrak{M}$ and $a=(a_1,a_2,\cdots,a_N)\in G:=G_1\times G_2\times\cdots\times G_N$, there holds
$$
   {\rm rank}\,\partial_ag(a,\eta)=n:=n_1+n_2+\cdots+n_N,
$$
  $\forall\eta\in\mathfrak{M}$, $\forall a\in G$, where $n_i=\dim G_i$, $i=1,2,\cdots,N$.
\end{enumerate}

  All concepts and notations appearing in Theorem 1.1 have been explained. In what follows we give the proof of Theorem 1.1.
  We need a preliminary lemma. Recall that for a given Banach space $X$ and given numbers $0<\alpha<1$
  and $T>0$, the notation $C^{\alpha}_{\alpha}((0,T],X)$ denotes the Banach space of bounded vector functions $u:(0,T]\to X$
  such that the vector function $t\mapsto t^{\alpha}u(t)$ is uniformly $\alpha$-H\"{o}lder continuous in $(0,T]$, with norm
$$
  \|u\|_{C^{\alpha}_{\alpha}((0,T],X)}=\sup_{0<t\leqslant T}\|u(t)\|_X+\sup_{0<s<t\leqslant T}
  \frac{\|t^{\alpha}u(t)-s^{\alpha}u(s)\|_X}{(t-s)^{\alpha}}
$$
  (cf. the introduction of Chapter 4 of \cite{Lun2}), and for given $\omega>0$, the notation $C^{\alpha}([T,\infty),X,-\omega)$
  denotes the Banach space of vector functions $u:[T,\infty)\to X$ such that the vector function $t\mapsto e^{\omega t}u(t)$
  is bounded and uniformly $\alpha$-H\"{o}lder continuous in $([T,\infty)$, with norm
$$
  \|u\|_{C^{\alpha}([T,\infty),X,-\omega)}=\sup_{t\geqslant T}\|e^{\omega t}u(t)\|_X+\sup_{t>s\geqslant T}
  \frac{\|e^{\omega t}u(t)-e^{\omega s}u(s)\|_X}{(t-s)^{\alpha}}
$$
  (cf. Section 4.4 of \cite{Lun2}).
\medskip

  {\bf Lemma 4.3}\ \ {\em Let $X$ be a Banach space and $X_0$ an embedded Banach subspace of $X$. Let $A$ be a sectorial
  operator in $X$ with domain $X_0$. Assume $\omega_-=-\sup\{{\rm Re}\lambda:\lambda\in\sigma(A)\}>0$ and $f\in
  C^{\alpha}_{\alpha}((0,1],X)\cap C^{\alpha}([1,\infty),X,-\omega)$, where  $0<\alpha<1$ and $\omega\in (0,\omega_-)$.
  Given $u_0\in X_0$, let $u(t)=e^{tA}u_0+\int_0^t e^{(t-s)A}f(s)ds$. Then $u\in C^{\alpha}_{\alpha}((0,1],X_0)\cap
  C^{\alpha}([1,\infty),X_0,-\omega)$, and there exists a constant $C=C(\alpha,\omega)>0$ such that}
\begin{equation}
   \|u\|_{C^{\alpha}_{\alpha}((0,T],X_0)}+\|u\|_{C^{\alpha}([T,\infty),X_0,-\omega)}\leqslant
   C(\|u_0\|_{X_0}+\|f\|_{C^{\alpha}_{\alpha}((0,T],X)}+\|f\|_{C^{\alpha}([T,\infty),X,-\omega)}).
\end{equation}
%---(4.1)---

  {\em Proof}.\ \ See Lemma 2.2 of \cite{Cui2}. $\quad\Box$
\medskip

  {\em Proof of Theorem 1.1}.\ \ The assertion (1) is an immediate consequence of the properties $(L3)$, $(L5)$ and $(L7)$
  of the Lie group actions $(G_i,p_i)$, $i=1,2,\cdots,N$. In what follows we prove the assertions (2) $\sim$ (4). For
  simplicity of notations we only consider the case $N=2$; for general $N$ the proof is similar.

  Since the equation (3.1) is of $C^2$-class parabolic type, there exists a local chart $(\mathcal{U},\varphi)\in\mathscr{A}$
  with $\eta^*\in\mathcal{U}$, such that the representation of this equation (3.3) in this local chart is of parabolic type
  and the representation $F$ of the vector field $\mathscr{F}$ is of $C^2$-class. Without loss of generality we assume that
  $\varphi(\eta^*)=0$, so that $F(0)=0$. The condition $(G3)$ and $(G4)$ imply that
$$
  {\rm dim}\,{\rm Ker}F'(0)=n, \quad \mbox{${\rm Range}F'(0)$ is closed in $X$},  \quad \mbox{and}
$$
$$
  X={\rm Ker}F'(0)\oplus{\rm Range}F'(0).
$$

  Firstly, since the equation (3.3) is of parabolic type and $F$ is of $C^2$-class, by Theorem 8.1.1 of \cite{Lun2}, there
  exists a neighborhood $U_0$ of the origin of $X_0$ such that for any $x_0\in U_0$, the initial value problem
\begin{equation}
\left\{
\begin{array}{ll}
  x'(t)=F(x(t)), &\quad t>0,\\
  x(0)=x_0
\end{array}
\right.
\end{equation}
%---(4.2)---
  has a unique local solution $x\in C([0,T],X_0)\cap C^1([0,T],X)\cap C^{\alpha}_{\alpha}((0,T],X_0)$, where $T=T(x_0)$
  depends on $x_0$ and $\alpha$ is an arbitrary number in $(0,1)$. Furthermore, denoting by $T^*(x_0)$ the supreme of all
  such $T$, then by Proposition 9.1.1 of \cite{Lun2} there exists $\varepsilon>0$ independent of $x_0$ such that if
  $\|x(t)\|_{X_0}<\varepsilon$ for all $t\in [0,T^*(x_0))$, then $T^*(x_0)=\infty$.

  In what follows we denote $M_c=\{\varphi(\eta):\eta\in \mathcal{M}_c\cap\mathcal{U}_0\cap\varphi^{-1}(U_0)\}$. Note that
  $0\in M_c$.

  Let $G_i^0$ be a small neighborhood of the unit element $e$ in $G_i$, $i=1,2$, and $\mathcal{U}'$, $\mathcal{U}'_0$ small
  neighborhoods of $\eta^*$ in $\mathfrak{M}$ and $\mathfrak{M}_0$, respectively, such that $p_i(a_i,\eta)\in
  \mathcal{U},\mathcal{U}_0$, respectively, for all $a_i\in G_i^0$, $i=1,2$, and $\eta\in\mathcal{U}',\mathcal{U}'_0$,
  respectively. The Lie group actions $(G_1,p_1)$ and $(G_2,p_2)$ in the Banach manifold $\mathfrak{M}$ induce local Lie
  group actions $(G_1^0,\tilde{p}_1)$ and $(G_2^0,\tilde{p}_2)$, respectively, in the open subset $U'=\varphi(\mathcal{U}')$
  of the Banach spaces $X$ as follows: For any $a_i\in G_i^0$ and $x\in U'$, define
$$
  \tilde{p}_i(a_i,x)=\varphi(p_i(a_i,\varphi^{-1}(x))), \quad i=1,2.
$$
  When restricted to $U_0'=\varphi(\mathcal{U}_0')$, they are local Lie group actions in this open subset of $X_0$.
  A simple computation easily shows that the condition $(G2)$ implies that for each $i=1,2$, $F$ is $\theta_i$-quasi-invariant
  with respect to the local Lie group action $(G_i^0,\tilde{p}_i)$, i.e.
\begin{equation}
  F(\tilde{p}_i(a_i,x))=\theta_i(a_i)\partial_{x}\tilde{p}_i(a_i,x)F(x), \quad
  \forall a_i\in G_i^0, \;\; \forall x\in U_0'.
\end{equation}
%---(4.3)---
  Note that since $\eta\in\mathcal{M}_c$ if and only if there exists $a=(a_1,a_2)\in G_1\times G_2$ such that
  $\eta=g(a,\eta^*)=p_1(a_1,p_2(a_2,\eta^*))$ and $\eta\in\mathfrak{M}_0$ in a small neighborhood of $\eta^*$
  has a such expression if and only if $\varphi(\eta)=\tilde{p}_1(a_1,\tilde{p}_2(a_2,0))$ for some $a_i\in G_i^0$,
  $i=1,2$, it follows that by suitably choosing the neighborhood $U_0$ of $0$ in $X_0$ and the neighborhoods
  $G_1^0$, $G_2^0$ of the unit elements in the groups $G_1$, $G_2$, respectively, we have $x\in M_c$ if and only
  there exists $a_i\in G_i^0$, $i=1,2$, such that $x=\tilde{p}_1(a_1,\tilde{p}_2(a_2,0))$. From this fact we see
$$
  M_c\subseteq S:=\{x\in U_0:F(x)=0\}.
$$

  Let $A=F'(0)$ and $N(x)=F(x)-Ax$. The $N\in C^2(U_0,X)$ and $N(0)=0$, $N'(0)=0$, so that by shrinking $U_0$ when necessary,
  there exists a constant $C>0$ such that
\begin{equation}
  \|N(x)\|_{X}\leqslant C\|x\|_{X_0}^2, \quad \|N(x)-N(y)\|_{X}\leqslant C(\|x\|_{X_0}+\|y\|_{X_0})\|x-y\|_{X_0},
  \quad \forall x,y\in U_0.
\end{equation}
%---(4.4)---
  The equation $(4.2)_1$ can be rewritten as follows:
\begin{equation}
  x'=Ax+N(x).
\end{equation}
%---(4.5)---
  Let $\sigma_-(A)=\sigma(A)\backslash\{0\}$. Choose a closed smooth curve in the complex plane such that it encloses the
  origin and separates it from $\sigma_-(A)$. We denote by $\Gamma$ this curve with anticlockwise orientation. Let $P\in L(X)$
  be the following operator:
$$
  P=\frac{1}{2\pi i}\int_{\Gamma}R(\lambda,A){\rm d}\lambda.
$$
  We know that $P$ is a bounded projection: $P^2=P$. Since $0$ is the unique element in $\sigma(A)$ enclosed by $\Gamma$
  and $X={\rm Ker}A\oplus {\rm Range}A$ with ${\rm Range}A$ closed, we have (cf. Proposition A.2.2 of \cite{Lun2})
$$
  PX=PX_0={\rm Ker}A, \quad (I-P)X={\rm Range}A.
$$
  Note that these relations imply that $APX=0$ and $A(I-P)X_0=(I-P)X$. Let $A_-=A|_{(I-P)X_0}$. Then $A_-:(I-P)X_0\to (I-P)X$
  is an isomorphism and $\sigma(A_-)=\sigma_-(A)$, so that
\begin{equation}
  \sup\{{\rm Re}\lambda:\lambda\in\sigma(A_-)\}=-\omega_-<0.
\end{equation}
%---(4.6)---
  Given $\delta,\delta'>0$, we denote by $B_1(0,\delta)$ and $B_2(0,\delta')$ spheres in $PX$ and $(I-P)X_0$, respectively,
  both centered at the origin but with radius $\delta$ and $\delta'$, respectively, i.e.,
$$
  B_1(0,\delta)=\{u\in PX=PX_0:\|u\|_{X_0}<\delta\}, \quad B_2(0,\delta')=\{v\in (I-P)X_0:\|v\|_{X_0}<\delta'\},
$$
  and let $O=\{u+v:u\in B_1(0,\delta),v\in B_2(0,\delta')\}\subseteq X_0$, $D=B_1(0,\delta)\times B_2(0,\delta')
  \subseteq PX\times(I-P)X_0$. We prove that if $\delta,\delta'>0$ are sufficiently small then there exists a mapping
  $h\in C^{2-0}(B_1(0,\delta),B_2(0,\delta'))$, $h(0)=0$, $h'(0)=0$, such that
\begin{equation}
  M_c\cap O=\{u+h(u):u\in B_1(0,\delta)\}.
\end{equation}
%---(4.7)---
  We first choose $\delta,\delta'>0$ small enough such that $O\subseteq U'$, and define $F_1:D\to PX$ and $F_2:D\to (I-P)X$
  respectively by
$$
  F_1(u,v)=PF(u+v), \quad F_2(u,v)=(I-P)F(u+v), \quad \forall(u,v)\in D.
$$
  Clearly $F_2\in C^{2-0}(D,(I-P)X)$ (surely also $F_1\in C^{2-0}(D,PX)$) and $F_2(0,0)=0$. Since $\partial_vF_2(0,0)
  =A_-:(I-P)X_0\to (I-P)X$ is an isomorphism, by the implicit function theorem we see that if $\delta,\delta'$ are
  sufficiently small then there exists a unique mapping $h\in C^{2-0}(B_1(0,\delta),B_2(0,\delta'))$ such that
  $h(0)=0$, $F_2(u,h(u))=0$, $\forall u\in B_1(0,\delta)$, and, furthermore, for any $u\in B_1(0,\delta)$, $v=h(u)$
  is the unique solution of the equation $F_2(u,v)=0$ in $B_2(0,\delta')$. For $x\in S\cap O$ let $u=Px$ and
  $v=(I-P)x$. Since $F(x)=0$, which implies $F_2(u,v)=(I-P)F(x)=0$, we infer that $v=h(u)$. Hence $S\cap
  O\subseteq\{u+h(u):u\in B_1(0,\delta)\}$, and, consequently,
$$
  M_c\cap O\subseteq S\cap O\subseteq\{u+h(u):u\in B_1(0,\delta)\}
$$
  Since $\mbox{dim}M_c=n=\mbox{dim}\{u+h(u):u\in B_1(0,\delta)\}$ (the second equality follows from the condition
  $(G_3$): $\mbox{dim Ker} A=n$), it follows that $M_c\cap O=\{u+h(u):u\in B_1(0,\delta)\}$. This proves (4.7).
  Note that this further implies that
$$
  F_1(u,h(u))=PF(u+h(u))=0, \quad \forall u\in B_1(0,\delta).
$$
  Note also that since $\partial_uF_2(0,0)PX=(I-P)APX_0=0$, we have $h'(0)=-[\partial_vF_2(0,0)]^{-1}\partial_uF_2(0,0)=0$.

  Let
$$
  N_1(u,v)=PN(u+v), \quad N_2(u,v)=(I-P)N(u+v), \quad \forall(u,v)\in D.
$$
  Since $AX_0=(I-P)X$, we have $PA=0$. Using this fact we see the differential equation (4.5) reduces into the following
  system of differential equations:
\begin{equation}
\left\{
\begin{array}{l}
  u'=N_1(u,v),\\
  v'=A_-v+N_2(u,v).
\end{array}
\right.
\end{equation}
%---(4.8)---
  Given $x_0\in O$, set $u_0=Px_0$, $v_0=(I-P)x_0$, and let $(u,v)=(u(t),v(t))$ be the solution of (4.8) subject to the
  initial condition $(u(0),v(0))=(u_0,v_0)$ defined in a maximal interval $[0,T^*)$ such that $(u(t),v(t))\in D$ for
  all $t\in[0,T^*)$. Since $(u,v)\equiv(0,0)$ is a solution of (4.8) defined for all $t\geqslant 0$, by continuous
  dependence of the solution on initial data we infer that by replacing $\delta,\delta'$ with smaller numbers when necessary,
  we may assume that $T^*=T^*(x_0)>1$ for all $x_0\in D$. Let $\phi(t)=v(t)-h(u(t))$, $t\in [0,T^*)$. Since
  $A_-h(u)+N_2(u,h(u))=F_2(u,h(u))=0$ and $N_1(u,h(u))=F_1(u,h(u))=0$ for all $u\in B_1(0,\delta)$, we see that $\phi'(t)=
  A_-\phi(t)+\sigma(t)$, $\forall t\in [0,T^*)$, where
$$
  \sigma(t)=[N_2(u(t),v(t))-N_2(u(t),h(u(t)))]-h'(u(t))[N_1(u(t),v(t))-N_1(u(t),h(u(t)))].
$$
  By the variation of constant formula, it follows that
$$
  \phi(t)={\rm e}^{tA_-}\phi(0)+\int_0^t{\rm e}^{(t-s)A_-}\sigma(s){\rm d}s, \quad t\in [0,T^*).
$$
  Using this expression and applying Lemma 4.3, we see that for any $0<\alpha<1$ and $\omega\in (0,\omega_-)$,
\begin{equation}
   \|\phi\|_{C^{\alpha}_{\alpha}((0,1],X_0)}+\|\phi\|_{C^{\alpha}([1,T^*),X_0,-\omega)}\leqslant
   C(\|\phi(0)\|_{X_0}+\|\sigma\|_{C^{\alpha}_{\alpha}((0,1],X)}+\|\sigma\|_{C^{\alpha}([1,T^*),X,-\omega)}).
\end{equation}
%---(4.9)---
  From (4.4) it is not hard to deduce that (cf. the proof of Theorem 9.1.2 of \cite{Lun2})
$$
   \|\sigma\|_{C^{\alpha}_{\alpha}((0,1],X)}+\|\sigma\|_{C^{\alpha}([1,T^*),X,-\omega)}\leqslant
   C\max\{\delta,\delta'\}(\|\phi\|_{C^{\alpha}_{\alpha}((0,1],X_0)}+\|\phi\|_{C^{\alpha}([1,T^*),X_0,-\omega)}).
$$
  Substituting this estimate into (4.9) we see that if $\delta,\delta'$ are chosen sufficiently small then
$$
  \|\phi\|_{C^{\alpha}_{\alpha}((0,1],X_0)}+\|\phi\|_{C^{\alpha}([1,T^*),X_0,-\omega)}\leqslant C\|\phi(0)\|_{X_0},
$$
  which implies, in particular, that
\begin{equation}
  \|\phi(t)\|_{X_0}\leqslant C{\rm e}^{-\omega t}\|\phi(0)\|_{X_0}, \quad \forall t\in [0,T^*).
\end{equation}
%---(4.10)---
  Next, since $N_1(u,h(u))=F_1(u,h(u))=0$, $\forall u\in B_1(0,\delta)$, we have
$$
  u'(t)=N_1(u(t),v(t))-N_1(u(t),h(u(t)))=:\sigma_1(t), \quad \forall t\in [0,T^*).
$$
  Again from (4.4) we have
$$
  \|\sigma_1(t)\|_{X}\leqslant C(\|u(t)\|_{X_0}+\|v(t)\|_{X_0})\|\phi(t)\|_{X_0}\leqslant
   C\|\phi(t)\|_{X_0}, \quad \forall t\in [0,T^*).
$$
   Hence
\begin{eqnarray}
  \|u(t)\|_{X_0}&\leqslant& C\|u(t)\|_{X}\leqslant C(\|u_0\|_{X}+\int_0^t\|\sigma_1(s)\|_{X}{\rm d}s)
  \leqslant C(\|u_0\|_{X}+\int_0^t\|\phi(s)\|_{X_0}{\rm d}s) \nonumber \\
  &\leqslant& C(\|u_0\|_{X}+\|\phi(0)\|_{X_0}), \quad \forall t\in [0,T^*).
\end{eqnarray}
%---(4.11)---
  In getting the first inequality we used the fact that $PX_0=PX={\rm Ker}A$ is a finite-dimensional space so that
  all norms in it are mutually equivalent. From (4.10) and (4.11) we see that if $\delta,\delta'$ are chosen sufficiently
  small then $\|x(t)\|_{X_0}<\varepsilon$ for all $t\in [0,T^*)$, so that $T^*=\infty$.

  Next, similarly as in the proof of (4.11) we see that for any $s>t\geqslant 0$,
\begin{equation}
  \|u(s)-u(t)\|_{X_0}\leqslant C\int_t^s\|\sigma_1(\tau)\|_{X}{\rm d}\tau\leqslant C\int_t^s\|\phi(\tau)\|_{X_0}{\rm d}\tau
  \leqslant C({\rm e}^{-\omega t}-{\rm e}^{-\omega s})\|\phi(0)\|_{X_0}.
\end{equation}
%---(4.12)---
  Hence $\displaystyle\lim_{t\to\infty}u(t)$ exists in $X_0$, which we denote as $\bar{u}$. Let $\bar{x}=\bar{u}+h(\bar{u})$.
  Then $\bar{x}\in M_c$ and
\begin{equation}
  \lim_{t\to\infty}x(t)=\lim_{t\to\infty}u(t)+\lim_{t\to\infty}h(u(t))+\lim_{t\to\infty}\phi(t)=\bar{u}+h(\bar{u})=\bar{x}.
\end{equation}
%---(4.13)---
  This shows that the solution of the equation (4.5) starting from an initial point located in the neighborhood $O$ of the
  origin converges to a point lying in $M_c$ as $t\to\infty$.

  Note that by letting $s\to\infty$ in (4.12), we have
\begin{equation}
  \|u(t)-\bar{u}\|_{X_0}\leqslant C{\rm e}^{-\omega t}\|\phi(0)\|_{X_0}, \quad \forall t\geqslant 0.
\end{equation}
%---(4.14)---
  From (4.10) and (4.14) we easily obtain
\begin{equation}
  \|x(t)-\bar{x}\|_{X_0}\leqslant C{\rm e}^{-\omega t}\|\phi(0)\|_{X_0}, \quad \forall t\geqslant 0.
\end{equation}
%---(4.15)---
  Hence, the solution of the equation (4.5) converges to its limit as $t\to\infty$ in exponential rate.

  We now prove that there exists a $C^{2-0}$-Banach manifold $M_s\subseteq O$ of codimension $n_1+n_2$, such that
  the solution $x=x(t)$ of the equation (4.5) converges to the origin as $t\to\infty$ if and only if $x(0)\in M_s$.
  To this end, for any $(u_0,v_0)\in D$ we denote by $u=u(t;u_0,v_0)$, $v=v(t;u_0,v_0)$ the unique solution of (4.8)
  with initial data $u(0)=u_0$, $v(0)=v_0$, and set $x(t;u_0,v_0)=u(t;u_0,v_0)+v(t;u_0,v_0)$, $\forall t\geqslant 0$.
  In what follows we first prove that $\displaystyle\lim_{t\to\infty}x(t;u_0,v_0)=0$ if and only if the improper
  integral $\displaystyle\int_0^{\infty}PN(x(\tau;u_0,v_0)){\rm d}\tau$ is convergent and the following relation holds:
\begin{equation}
  u_0+\int_0^{\infty}PN(x(\tau;u_0,v_0)){\rm d}\tau=0.
\end{equation}
%---(4.16)---
  Indeed, if $\displaystyle\lim_{t\to\infty}x(t;u_0,v_0)=0$ then from (4.4) and (4.15) (with $\bar{x}=0$) we easily
  see that the improper integral $\displaystyle\int_0^{\infty}PN(x(\tau;u_0,v_0)){\rm d}\tau$ is convergent. Integration
  the first equation in (4.8) we get
\begin{equation}
  u(t;u_0,v_0)=u_0+\int_0^tPN(x(\tau;u_0,v_0)){\rm d}\tau, \quad t\geqslant 0.
\end{equation}
%---(4.17)---
  Since $\displaystyle\lim_{t\to\infty}x(t;u_0,v_0)=0$ implies both $\displaystyle\lim_{t\to\infty}u(t;u_0,v_0)=0$ and
  $\displaystyle\lim_{t\to\infty}v(t;u_0,v_0)=0$, by letting $t\to\infty$ in both sides of the above equation we see
  that (4.16) follows. Conversely, if the improper integral $\displaystyle\int_0^{\infty}PN(x(\tau;u_0,v_0)){\rm d}\tau$
  is convergent and (4.16) holds true, then by letting $t\to\infty$ in (4.17) we see that $\displaystyle\lim_{t\to\infty}
  u(t;u_0,v_0)=0$. It then follows from (4.13) that $\displaystyle\lim_{t\to\infty}x(t;u_0,v_0)=0$. This proves the
  desired assertion. Now, by regarding (4.16) as an implicit function equation and applying the implicit function theorem,
  we can easily show that if $\delta,\delta'$ are sufficiently small then (4.16) defines an implicit function
  $u_0=q(v_0)$, where $q\in C^{2-0}(B_2(0,\delta'),B_1(0,\delta))$ and $q(0)=0$, cf. the proof of Theorem 2.1 of
  \cite{Cui2} for details. Hence, the equation (4.16) defines a $C^{2-0}$-Banach manifold in $D$. We denote by $M_s$ the
  corresponding $C^{2-0}$-Banach manifold in $O$. Note that the definition of $M_s$ ensures that $\displaystyle\lim_{t\to\infty}
  x(t;u_0,v_0)=0$ if and only if $x_0=u_0+v_0\in M_s$.

  We now prove that for any $x_0\in O$ there exist unique $a\in G_1$, $b\in G_2$ and $y_0\in M_s$ such that $x_0=
  \tilde{p}_1(a,\tilde{p}_2(b,y_0))$ and
$$
  \lim_{t\to\infty}x(t)=\tilde{p}_1(a,\tilde{p}_2(b,0)),
$$
  where $x=x(t)$ is the solution of (4.5) with initial data $x(0)=x_0$. Indeed, let $\bar{x}$ be as in (4.13). Since
  $\bar{x}\in M_c$, there exist unique $a\in G_1^0$, $b\in G_2^0$ such that $\tilde{p}_1(a,\tilde{p}_2(b,0))=\bar{x}$.
  Let $y_0=\tilde{p}_2(b^{-1},\tilde{p}_1(a^{-1},x_0))$ and
$$
  y(t)=\tilde{p}_2(b^{-1},\tilde{p}_1(a^{-1},x(t\theta_1(a^{-1})\theta_2(b^{-1})))) \quad  \mbox{for}\;\, t\geqslant 0.
$$
  $y(t)$ is the solution of (4.5) with initial data $y(0)=y_0$. From the condition $\displaystyle\lim_{t\to\infty}x(t)
  =\bar{x}$ we have
$$
  \lim_{t\to\infty}y(t)=\tilde{p}_2(b^{-1},\tilde{p}_1(a^{-1},\bar{x}))=0.
$$
  Hence $y_0\in M_s$. This proves existence. Uniqueness is obvious.

  Now let $\mathcal{O}_{\eta^*}=\varphi^{-1}(O)$ and $\mathcal{M}_s=\varphi^{-1}(M_s)$. Then $\mathcal{O}_{\eta^*}$ is a
  neighborhood of $\eta^*$ and $\mathcal{M}_s$ is a submanifold of $\mathfrak{M}_0$ satisfying the condition (3). Moreover,
  for any $\eta_0\in\mathcal{O}_{\eta^*}$ the initial value problem (3.2) has a unique solution $\eta\in
  C([0,\infty),\mathfrak{M})\cap C^1((0,\infty),\mathfrak{M}_0)$ and there exist unique $a\in G_1$, $b\in G_2$ and
  $\xi_0\in\mathcal{M}_s$ such that $\eta_0=g(a,b,\xi_0)$ and for the solution $\eta=\eta(t)$ of (3.2) with initial data
  $\eta(0)=\eta_0$,
$$
   \lim_{t\to\infty}\eta(t)=g(a,b,\eta^*).
$$
  To finish the proof of Theorem 4.3 we now only need to repeat the above argument for every point $\eta\in\mathcal{M}_c$
  and then glue all the open sets $\mathcal{O}_{\eta}$ together to form the neighborhood $\mathcal{O}$ of $\mathcal{M}_c$.
  This completes the proof. $\quad\Box$
\medskip

  {\em Remark}.\ \ $\mathcal{M}_c$ is called the {\em center manifold} of the equation (3.1), and $\mathcal{M}_s$ is
  called the {\em stable manifold} of the equation (3.1) corresponding to the stationary point $\eta^*$.

\section{Proof of Theorem 1.2}
\setcounter{equation}{0}

\hskip 2em
  Before giving the proof of Theorem 1.2, let us first point out that application of Theorem 1.1 to the Neumann initial
  boundary value problem of the heat equation gives a new explanation of the dynamical
  behavior of the solution of that problem. Let $\Omega$ be a given bounded domain in ${\mathbb{R}}^n$ with a $C^2$-boundary.
  Recall that the Neumann initial-boundary value problem is as follows:
\begin{equation}
\left\{
\begin{array}{ll}
  \partial_tu(x,t)=\Delta u(x,t), &\quad x\in\Omega, \;t>0,\\
  \partial_{\bfn}u(x,t)=0, &\quad x\in\partial\Omega, \;t>0,\\
  u(x,0)=u_0(x), &\quad x\in\Omega,
\end{array}
\right.
\end{equation}
%---(5.1)---
  where $\bfn$ denotes the outward-pointing unit normal field of $\partial\Omega$ and $u_0\in L^2(\Omega)$. Clearly,
  stationary solutions of this problem make up an one-dimensional manifold $\mathcal{M}_c={\mathbb{R}}{\mathbf{1}}$,
  where ${\mathbf{1}}$ denotes the function in $\overline{\Omega}$ with values identically $1$. The solution of
  the above problem has the following asymptotic behavior as $t\to\infty$:
\begin{equation}
   \lim_{t\to\infty}u(x,t)=\frac{1}{|\Omega|}\int_{\Omega}u_0(x){\rm d}x \quad \mbox{in $L^2(\Omega)$ norm}.
\end{equation}
%---(5.2)---
  Let $\mathfrak{M}=X=L^2(\Omega)$ and $\mathfrak{M}_0=X_0=\{u\in H^2(\Omega):
  \partial_{\bfn}u|_{\partial\Omega}=0\}$. Let $G={\mathbb{R}}$ be as in Section 1. We introduce an action $p$ of $G$ to
  $X=L^2(\Omega)$ as follows:
$$
  p(a,u)=a{\mathbf{1}}+u, \quad \forall u\in L^2(\Omega), \quad \forall a\in{\mathbb{R}}.
$$
  Note that the restriction of this group action to $X_0$ is also a group action to $X_0$. It is easy to see that
  the differential equation in the Banach space $X=L^2(\Omega)$ corresponding to the problem (5.1) is invariant under
  the Lie group action $(G,p)$ defined here, and Theorem 1.1 applies to it. The stable manifold of this equation
  corresponding to the stationary point $u^*=0$ is
$$
   \mathcal{M}_s=\Big\{u\in X_0:\int_{\Omega}u(x){\rm d}x=0\Big\}.
$$
  For any $u_0\in X_0$, let $v_0=u_0-\displaystyle\Big(\frac{1}{|\Omega|}\int_{\Omega}u_0(x){\rm d}x\Big){\mathbf{1}}$
  and $a=\displaystyle\frac{1}{|\Omega|}\int_{\Omega}u_0(x){\rm d}x$. Clearly $v_0\in\mathcal{M}_c$, $u_0=p(a,v_0)$,
  and $u(t)=p(a,v(t))$, $\forall t\geqslant 0$, where $u(t)$ and $v(t)$ are solutions of (4.3) with respect to initial
  data $u_0$ and $v_0$, respectively. The last relation implies that
$$
   \lim_{t\to\infty}u(t)=p(a,0)=\Big(\frac{1}{|\Omega|}\int_{\Omega}u_0(x){\rm d}x\Big){\mathbf{1}}
   \quad \mbox{in $L^2(\Omega)$ norm},
$$
  which recovers the formula (5.2).

  Let us now consider the problem (1.7). Let $m$, $\mu$, $\mathfrak{M}$ and $\mathfrak{M}_0$ be as in Section 1, i.e.,
  $m$ is a positive integer $\geqslant 2$, $0<\mu<1$, $\mathfrak{M}:=\dot{\mathfrak{S}}^{m+\mu}({\mathbb{R}}^n)$ and
  $\mathfrak{M}_0:=\dot{\mathfrak{S}}^{m+3+\mu}({\mathbb{R}}^n)$. We define a quasi-differentiable structure in $\mathfrak{M}$
  with kernel $\mathfrak{M}_0$ as follows: Let $\Sigma$ be an arbitrary smooth closed hypersurface in $\mathbb{R}^n$
  homeomorphic to $\mathbb{S}^{n-1}$. Choose a sufficiently small number $\delta>0$ such that the $\delta$-neighborhood of
  $\Sigma$ is homeomorphic to the $n$-dimensional submanifold $\Sigma\times(-\delta,\delta)$ of $\mathbb{R}^{n+1}$. Let
  ${\mathbf{n}}$ be the normal field of $\Sigma$, outward-pointing with respect to the domain enclosed by $\Sigma$. For any
  $S\in\mathfrak{M}$ contained in the $\delta$-neighborhood of $\Sigma$, let $\rho\in\dot{C}^{m+\mu}(\mathbb{S}^{n-1})$ be
  the unique function in $\mathbb{S}^{n-1}$ such that the following relation holds:
$$
   S=\{x+\rho(x){\mathbf{n}}(x):\; x\in\mathbb{S}^{n-1}\}.
$$
  We define a map $\varphi:\mathfrak{M}\to\dot{C}^{m+\mu}(\mathbb{S}^{n-1})$ by defining $\varphi(S)=\rho$. Let $\mathcal{U}$
  be the subset of $\mathfrak{M}$ consisting of all $S\in\mathfrak{M}$ contained in the $\delta$-neighborhood of $\Sigma$.
  Then $(\mathcal{U},\varphi)$ is a local chart of $\mathfrak{M}$. Let $\mathscr{A}$ be the set of all such pairs
  $(\mathcal{U},\varphi)$. It follows that $(\mathfrak{M},\mathscr{A})$ is a quasi-differentiable Banach manifold with kernel
  $\mathfrak{M}_0$, inner kernel $\mathfrak{M}_1=\dot{\mathfrak{S}}^{m+4+\mu}({\mathbb{R}}^n)$, and shell
  $\widetilde{\mathfrak{M}}=\dot{\mathfrak{S}}^{m-1+\mu}({\mathbb{R}}^n)$. The base Banach spaces of $\mathfrak{M}$,
  $\mathfrak{M}_0$, $\mathfrak{M}_1$, and $\widetilde{\mathfrak{M}}$ are respectively $X=\dot{C}^{m+\mu}(\mathbb{S}^{n-1})$,
  $X_0=\dot{C}^{m+3+\mu}(\mathbb{S}^{n-1})$, $X_1=\dot{C}^{m+4+\mu}(\mathbb{S}^{n-1})$ and $\widetilde{X}=
  \dot{C}^{m-1+\mu}(\mathbb{S}^{n-1})$.

  In Section 1 we stated that the problem (1.7) can be reduced into a differential equation in the Banach manifold
  $\mathfrak{M}$. In what follows we give the proof of this statement. For this purpose, we need first make some preparations.
  Let $K,R$ be two positive numbers such that the condition (1.8) is satisfied. Given $\rho,\eta\in C^{m+\mu}({\mathbb{S}}^{n-1})$
  with $\|\rho\|_{C^{m+\mu}({\mathbb{S}}^{n-1})}$ and $\|\eta\|_{C^{m+\mu}({\mathbb{S}}^{n-1})}$ sufficiently small, we denote
$$
%  \Omega_{\rho}=\{x\in{\mathbb{R}}^n: r<R[1+\rho(\omega)]\}, \quad
  D_{\rho,\eta}=\{x\in{\mathbb{R}}^n: K[1+\eta(\omega)]<r<R[1+\rho(\omega)]\},
$$
$$
   S_{\rho}=\{x\in{\mathbb{R}}^n: r=R[1+\rho(\omega)]\}, \quad \mbox{and} \quad
   \Gamma_{\eta}=\{x\in{\mathbb{R}}^n: r= K[1+\eta(\omega)]\}.
$$
  Here and hereafter we use $r=r(x)$ and $\omega=\omega(x)$ to denote the polar and spherical coordinates of the
  variable $x$, i.e., $r(x)=|x|$ and $\omega(x)=x/|x|$ for $x\neq 0$. Later on we shall also use the following abbreviations:
$$
  D=D_{0,0}=B(0,R)\backslash\overline{B(0,K)}, \quad S_0=\partial B(0,R), \quad \mbox{and} \quad  \Gamma_0=\partial B(0,K).
$$
  Consider the following free-boundary problem: Given $\rho\in\dot{C}^{m+\mu}({\mathbb{S}}^{n-1})$ with
  $\|\rho\|_{C^{m+\mu}({\mathbb{S}}^{n-1})}$ sufficiently small, find a triple $(\eta,u,\bfc)$ with $\eta\in
  \dot{C}^{m+\mu}({\mathbb{S}}^{n-1})$, $u\in\dot{C}^{m-2+\mu}(\overline{D}_{\rho,\eta})$ and $\bfc\in\mathbb{R}^n$ such that
  the following equations are satisfied:
\begin{equation}
\left\{
\begin{array}{ll}
   \Delta u(x)=0, &\quad x\in D_{\rho,\eta},\\
   u(x)=\gamma\kappa(x),   &\quad x\in S_{\rho},\\
   u(x)=\mu\kappa(x),   &\quad x\in\Gamma_{\eta},\\
   \partial_{\mbox{\footnotesize$\bfn$}}u(x)=\bfc\cdot\bfn,   &\quad x\in\Gamma_{\eta},\\
   \displaystyle\frac{1}{|S_{\rho}|}\oint_{S_{\rho}}x{\rm d}\sigma
   =&\displaystyle\frac{1}{|\Gamma_{\eta}|}\oint_{\Gamma_{\eta}}x{\rm d}\sigma.
\end{array}
\right.
\end{equation}
%---(5.3)---

  {\bf Lemma 5.1}\ \ {\em Given $\rho\in\dot{C}^{m+\mu}({\mathbb{S}}^{n-1})$ with $\|\rho\|_{C^{m+\mu}({\mathbb{S}}^{n-1})}$
  sufficiently small, the problem $(5.3)$ has a unique solution $(\eta,u,\bfc)$ with $\eta\in\dot{C}^{m+\mu}({\mathbb{S}}^{n-1})$,
  $u\in\dot{C}^{m-2+\mu}(\overline{D}_{\rho,\eta})$ and $\bfc\in\mathbb{R}^n$. Moreover, the map $\rho\mapsto(\eta,u,\bfc)$ from
  a small neighborhood of the origin of $\dot{C}^{m+\mu}({\mathbb{S}}^{n-1})$ to $\dot{C}^{m+\mu}({\mathbb{S}}^{n-1})\times
  \dot{C}^{m-2+\mu}(\overline{D}_{\rho,\eta})\times\mathbb{R}^n$ is smooth.}
\medskip

  {\em Remark}.\ \ Note that the map $\rho\mapsto u$ from a small neighborhood of the origin of
  $\dot{C}^{m+\mu}({\mathbb{S}}^{n-1})$ to $\dot{C}^{m-2+\mu}(\overline{D}_{\rho,\eta})$ is smooth means the following: By
  letting $\Psi_{\rho,\eta}:\overline{D}_{\rho,\eta}\to\overline{D}$ be the Hanzawa transformation (see (5.4) below), the
  map $\rho\mapsto u\circ\Psi_{\rho,\eta}^{-1}$ from a small neighborhood of the origin of $\dot{C}^{m+\mu}({\mathbb{S}}^{n-1})$
  to $\dot{C}^{m-2+\mu}(\overline{D})$ is smooth. Note that $\Psi_{\rho,\eta}$ depends smoothly on $\rho$ and $\eta$.
\medskip

  {\em Proof of Lemma 5.1}.\ \ We use the inverse function theorem to prove this lemma. For sufficiently small $\delta,\delta'>0$
  we denote
$$
  O_{\delta}=\{\rho\in\dot{C}^{m+\mu}({\mathbb{S}}^{n-1}):\|\rho\|_{C^{m+\mu}({\mathbb{S}}^{n-1})}<\delta\}, \quad
  O_{\delta'}'=\{\eta\in\dot{C}^{m+\mu}({\mathbb{S}}^{n-1}):\|\eta\|_{C^{m+\mu}({\mathbb{S}}^{n-1})}<\delta'\};
$$
  they are open subsets of $\dot{C}^{m+\mu}({\mathbb{S}}^{n-1})$. Choose a function $\phi\in C^{\infty}[ K,R]$ such that it
  satisfies the following conditions:
$$
  0\leqslant\phi\leqslant 1; \quad \phi(R)=\phi( K)=1; \quad
  \phi(t)=0 \;\; \mbox{for}\;\,\frac{3}{4} K+\frac{1}{4}R\leqslant t\leqslant\frac{1}{4} K+\frac{3}{4}R;
$$
$$
  \phi'(t)\leqslant 0 \;\; \mbox{for}\;\, K\leqslant t\leqslant\frac{3}{4} K+\frac{1}{4}R; \quad
  \phi'(t)\geqslant 0 \;\; \mbox{for}\;\,\frac{1}{4} K+\frac{3}{4}R\leqslant t\leqslant R.
$$
  Let $M_0=\displaystyle\max_{ K\leqslant t\leqslant R}|\phi'(t)|$ and assume $\delta,\delta'$ are small enough such that
  $\delta<(1+M_0R)^{-1}$, $\delta'<(1+M_0 K)^{-1}$ and $\max\{\delta,\delta'\}<\displaystyle\frac{1}{3}\frac{R- K}{R+ K}$.
  Consider the variable transformation $y=\Psi_{\rho,\eta}(x)$ from $\overline{D}_{\rho,\eta}$ to $\overline{D}$, where for
  $x\in\overline{D}_{\rho,\eta}$,
\begin{equation}
  \Psi_{\rho,\eta}(x)=\left\{
\begin{array}{ll}
   \displaystyle x-R\rho(\omega)\phi\Big(\frac{r}{1+\rho(\omega)}\Big)\omega &\quad\;\;
   \mbox{if} \;\; r\geqslant\frac{1}{2}( K+R),\\ [0.3cm]
   \displaystyle x- K\eta(\omega)\phi\Big(\frac{r}{1+\eta(\omega)}\Big)\omega &\quad\;\;
   \mbox{if} \;\; r<\frac{1}{2}( K+R).
\end{array}
\right.
\end{equation}
%---(5.4)---
  It is easy to see that $\Psi_{\rho,\eta}$ is a $\dot{C}^{m+\mu}$ diffeomorphism from $\overline{D}_{\rho,\eta}$ onto $\overline{D}$.
  Moreover, denoting
$$
  E^1_{\rho}=\{x\in{\mathbb{R}}^n:  \frac{1}{2}( K+R)<r<R[1+\rho(\omega)]\}, \quad
  E^1=\{x\in{\mathbb{R}}^n:  \frac{1}{2}( K+R)<r<R\},
$$
$$
  E^2_{\eta}=\{x\in{\mathbb{R}}^n:  K[1+\eta(\omega)]<r<\frac{1}{2}( K+R)\}, \quad
  E^2=\{x\in{\mathbb{R}}^n:  K<r<\frac{1}{2}( K+R)\},
$$
  we see that the restriction of $\Psi_{\rho,\eta}$ on $\overline{E}^1_{\rho}$ is independent of $\eta$ and the restriction of
  $\Psi_{\rho,\eta}$ on $\overline{E}^2_{\eta}$ is independent of $\rho$. Hence we re-denote the restrictions of $\Psi_{\rho,\eta}$
  on $\overline{E}^1_{\rho}$ and $\overline{E}^2_{\eta}$ as $\Psi_{1\rho}$ and $\Psi_{2\eta}$, respectively, and denote by
  $\psi_{1\rho}$ and $\psi_{2\eta}$ the restrictions of $\Psi_{1\rho}$ and $\Psi_{2\eta}$ to $S_{\rho}$ and $\Gamma_{\eta}$,
  respectively, which are clearly $\dot{C}^{m+\mu}$-diffeomorphisms from $S_{\rho}$ onto $S_0$ and from $\Gamma_{\eta}$ onto
  $\Gamma_0$, respectively. Now define operators $\mathcal{A}(\rho,\eta):\dot{C}^{m+\mu}(\overline{D})\to
  \dot{C}^{m-2+\mu}(\overline{D})$ and $\mathcal{N}(\eta,\bfc):\dot{C}^{m-2+\mu}(\overline{E^2})\to\dot{C}^{m-3+\mu}(\Gamma_0)$
  respectively as follows:
$$
  \mathcal{A}(\rho,\eta)v=[\Delta(v\circ\Psi_{\rho,\eta})]\circ\Psi_{\rho,\eta}^{-1} \quad\;\;
   \mbox{for}\;\,v\in\dot{C}^{m+\mu}(\overline{D}),
$$
$$
  \mathcal{N}(\eta,\bfc)v=[\partial_{{\mathbf{n}}}(v\circ\Psi_{2\eta})|_{\Gamma_{\eta}}-\bfc\cdot\bfn]
  \circ\psi_{2\eta}^{-1} \quad\;\; \mbox{for}\;\,v\in\dot{C}^{m-2+\mu}(\overline{E^2}).
$$
  We further introduce two operators $\mathcal{K}_1:\dot{C}^{m+\mu}({\mathbb{S}}^{n-1})\to\dot{C}^{m-2+\mu}(S_0)$ and
  $\mathcal{K}_2:\dot{C}^{m+\mu}({\mathbb{S}}^{n-1})\to\dot{C}^{m-2+\mu}(\Gamma_0)$ as follows: For arbitrary
  $\rho,\eta\in\dot{C}^{m+\mu}({\mathbb{S}}^{n-1})$,
$$
  \mathcal{K}_1(\rho)(x)=\mbox{the mean curvature of $S_{\rho}$ at $\psi_{1\rho}^{-1}(x)$}, \;\; \forall x\in S_0,
$$
$$
  \mathcal{K}_2(\eta)(x)=\mbox{the mean curvature of $\Gamma_{\eta}$ at $\psi_{2\eta}^{-1}(x)$}, \;\; \forall x\in\Gamma_0.
$$
  It follows that by letting $v=u\circ\Psi_{\rho,\eta}^{-1}$, the problem (5.3) transforms into the following problem:
\begin{equation}
\left\{
\begin{array}{ll}
   \mathcal{A}(\rho,\eta)v=0, &\quad \mbox{in}\; D,\\
   v=\gamma\mathcal{K}_1(\rho),   &\quad  \mbox{on}\; S_0,\\
   v=\mu\mathcal{K}_2(\eta),   &\quad  \mbox{on}\; \Gamma_0,\\
   \mathcal{N}(\eta,\bfc)v=0,   &\quad  \mbox{on}\; \Gamma_0,\\
   \displaystyle\frac{1}{|S_{\rho}|}\oint_{S_{\rho}}x{\rm d}\sigma=&\displaystyle\frac{1}{|\Gamma_{\eta}|}\oint_{\Gamma_{\eta}}x{\rm d}\sigma.
\end{array}
\right.
\end{equation}
%---(5.5)---
  Given $(\rho,\eta)\in O_{\delta}\times O_{\delta'}'\subseteq C^{m+\mu}({\mathbb{S}}^{n-1})\times C^{m+\mu}({\mathbb{S}}^{n-1})$,
  the equation $(5.5)_1$ subject to the boundary value conditions $(5.5)_2$ and $(5.5)_3$ has clearly a unique solution
  $v\in C^{m-2+\mu}(\overline{D})$ which we denote as $v_{\rho\eta}$. We put $\mathcal{Q}(\rho,\eta)=v_{\rho\eta}$. In addition,
  we denote
$$
  \mathcal{J}(\rho,\eta)=\frac{1}{|\Gamma_{\eta}|}\oint_{\Gamma_{\eta}}x{\rm d}\sigma
  -\frac{1}{|S_{\rho}|}\oint_{S_{\rho}}x{\rm d}\sigma.
$$
  Note that $\mathcal{J}(\rho,\eta)\in\mathbb{R}^n$. Now we define a mapping
$$
  \mathcal{P}:O_{\delta}\times O_{\delta'}'\times\mathbb{R}^n\to C^{m-3+\mu}(\Gamma_0)\times\mathbb{R}^n
$$
  as follows: For $\rho\in O_{\delta}$, $\eta\in O_{\delta'}'$ and $\bfc\in\mathbb{R}^n$,
$$
  \mathcal{P}(\rho,\eta,\bfc)=\Big(\mathcal{N}(\eta,\bfc)\mathcal{Q}(\rho,\eta),%+\mathcal{L}(\bfc),
  \mathcal{J}(\rho,\eta)\Big)
$$
  Clearly, $\mathcal{P}$ is a smooth map and
\begin{equation}
  \mathcal{P}(0,0,0)=(0,0).
\end{equation}
%---(5.6)---
  The above relation follows from the fact $v_{00}=\gamma/R=\mu/K$ (the constant function). In what follows we compute
  $\partial_{\eta,\mbox{\small $\bfc$}}\mathcal{P}(0,0,0)$ to show that it is an isomorphism of $C^{m+\mu}({\mathbb{S}}^{n-1})\times\mathbb{R}^n$
  onto $C^{m-3+\mu}(\overline{D})\times\mathbb{R}^n$.

  Given $\zeta\in C^{m+\mu}({\mathbb{S}}^{n-1})$, we denote $w=\partial_{\eta}\mathcal{Q}(0,0)\zeta$. By some similar computation
  as those made in the references \cite{Cui2, Cui4, FriHu}, we see that $w$ is the solution of the following problem:
\begin{equation}
\left\{
\begin{array}{ll}
   \Delta w=0 &\quad \mbox{in}\; D,\\
   w=0   &\quad  \mbox{on}\; S_0,\\
   w=-\frac{\mu}{K}(\zeta+\frac{1}{n\!-\!1}\Delta_{\omega}\zeta)  &\quad  \mbox{on}\; \Gamma_0,
\end{array}
\right.
\end{equation}
%---(5.7)---
  where $\Delta_{\omega}$ is the Laplace-Beltrami operator on ${\mathbb{S}}^{n-1}$. We now compute $w$ by using spherical
  harmonic expansion of the function on ${\mathbb{S}}^{n-1}$.

  Let $\{\lambda_k\}_{k=0}^{\infty}$ be the increasing sequence of all distinct eigenvalues of the operator $-\Delta_{\omega}$,
  and $\mathscr{H}^k(\mathbb{S}^{n-1})$ be the eigenspace of $-\Delta_{\omega}$ corresponding to the eigenvalue $\lambda_k$,
  $k=0,1,2,\cdots$. Let $d_k=\mbox{\rm dim}\,\mathscr{H}^k(\mathbb{S}^{n-1})$, $k=0,1,2,\cdots$. Recall that
$$
  \lambda_0=0, \quad \lambda_1=n\!-\!1, \quad \lambda_k=k^2+(n\!-\!2)k, \quad k=2,3,\cdots,
$$
$$
  d_0=1, \quad d_1=n \quad\mbox{and} \quad d_k={n\!+\!k\!-\!1\choose k}-{n\!+\!k\!-\!3\choose k\!-\!2}\quad\mbox{for}\;\,k\geqslant 2.
$$
  For each integer $k\geqslant 0$, let $\{Y_{kl}(\omega)\}_{l=1}^{d_k}$ be a normalized orthogonal basis of
  $\mathscr{H}^k(\mathbb{S}^{n-1})$ as a subspace of $L^2({\mathbb{S}}^{n-1})$. Note that $Y_{00}(\omega)=1/\sqrt{\sigma_n}$,
  where $\sigma_n=\frac{2\pi^{n/2}}{\Gamma(n/2)}$ is the surface measure of $\mathbb{S}^{n-1}$. For $k=1$ we particularly put
$$
  Y_{1l}(\omega)=\sqrt{n}\sigma_n^{-\frac{1}{2}}\omega_l, \qquad l=1,2,\cdots,n.
$$
  Note that the above definitions imply that the following relations hold:
$$
  \Delta_{\omega}Y_{kl}(\omega)=-\lambda_kY_{kl}(\omega), \quad k=0,1,2,\cdots,\;\; l=1,2,\cdots,d_k.
$$
  It is well-known that the following relation holds:
$$
  \Delta w=\frac{\partial^2w}{\partial r^2}+\frac{n\!-\!1}{r}\frac{\partial w}{\partial r}
  +\frac{1}{r^2}\Delta_{\omega}w.
$$
  Using these facts, we easily see that if a given function $\zeta\in C^{\infty}({\mathbb{S}}^{n-1})$ has a spherical
  harmonics expansion
$$
  \zeta(\omega)=\sum_{k=0}^{\infty}\sum_{l=1}^{d_k}a_{kl}Y_{kl}(\omega),
$$
  then the solution $w$ of the problem (5.7) is given by
$$
  w(r,\omega)=\frac{\mu}{K}\sum_{k=0}^{\infty}\sum_{l=1}^{d_k}\frac{\lambda_k\!-\!n\!+\!1}{n\!-\!1}
  \frac{K^{2k+n-2}}{R^{2k+n-2}\!-\!K^{2k+n-2}}\Big(\frac{r}{K}\Big)^k\Big[\Big(\frac{R}{r}\Big)^{2k+n-2}-1\Big]a_{kl}Y_{kl}(\omega)
  \quad  \mbox{for}\;\;  K\leqslant r\leqslant R.
$$
  Since $\mathcal{Q}(0,0)=\mbox{const.}$ so that $[\partial_{\eta}\mathcal{N}(0,0)\zeta]\mathcal{Q}(0,0)=0$, and $\mathcal{N}(0,0)
  =\partial_r|_{r=K}$, it follows that the first component of $\partial_{\eta}\mathcal{P}(0,0,0)\zeta$ is given by
\begin{align*}
  \mathcal{B}_1(\zeta)=&\mathcal{N}(0,0)\partial_{\eta}\mathcal{Q}(0,0)\zeta+[\partial_{\eta}\mathcal{N}(0,0)\zeta]\mathcal{Q}(0,0)
  =(\partial_rw)(K,\omega)
\nonumber\\
  =&-\frac{\mu}{K^2}\sum_{k=0}^{\infty}\sum_{l=1}^{d_k}\frac{\lambda_k\!-\!n\!+\!1}{n\!-\!1}
  \frac{kK^{2k+n-2}\!+\!(k\!+\!n\!-\!2)R^{2k+n-2}}{R^{2k+n-2}\!-\!K^{2k+n-2}}a_{kl}Y_{kl}(\omega).
\end{align*}

  To compute the second component of $\partial_{\eta}\mathcal{P}(0,0,0)\zeta$, we need to use the following formula: If a closed
  hypersurface $\Sigma$ in $\mathbb{R}^n$ is given by the equation $r=r(\omega)$, $\omega\in\mathbb{S}^{n-1}$, where
  $r(\omega)$ is a $C^1$-function on $\mathbb{S}^{n-1}$, then for any continuous function $f$ on $\Sigma$ we have
$$
  \int_{\Sigma}\!f(x)\mbox{\rm d}\sigma=\int_{\mathbb{S}^{n-1}}\!\!f(r(\omega)\omega)
  \sqrt{|r(\omega)|^2+|\nabla_{\omega}r(\omega)|^2}|r(\omega)|^{n-2}\mbox{\rm d}\omega,
$$
  where $\mbox{\rm d}\sigma$ and $\mbox{\rm d}\omega$ denote the measure elements on $\Sigma$ and $\mathbb{S}^{n-1}$, respectively,
  induces by the Lebesgue measure on $\mathbb{R}^n$, and $\nabla_{\omega}$ denotes gradient of functions on $\mathbb{S}^{n-1}$.
  Using this formula, we can easily verify that the second component of $D_{\eta}\mathcal{P}(0,0,0)\zeta$ is as follows:
$$
  \mathcal{B}_2(\zeta)=\sqrt{n}\sigma_n^{-\frac{1}{2}}K\mathcal{I}(\zeta), \qquad
  \mathcal{I}(\zeta)=\sqrt{n}\sigma_n^{-\frac{1}{2}}\int_{\mathbb{S}^{n-1}}\!\!\omega\zeta(\omega)\mbox{\rm d}\omega.
$$
  Note that $\mathcal{I}$ is an isomorphism of $\mathscr{H}^1(\mathbb{S}^{n-1})$ onto $\mathbb{R}^n$ when restricted to
  $\mathscr{H}^1(\mathbb{S}^{n-1})$. Note also that both $\mathcal{B}_1$ and $\mathcal{B}_2$ are Fourier multipliers,
  with orders $3$ and $-\infty$, respectively. As a result, both $\mathcal{B}_1$ and $\mathcal{B}_2$ can be extended
  into continuous linear operators in the distribution space $\mathscr{D}'({\mathbb{S}}^{n-1})=
  \mathscr{S}'({\mathbb{S}}^{n-1})$ on $\mathbb{S}^{n-1}$.

  Let $\mathcal{L}$ be the following map from $\mathbb{R}^n$ to $\mathscr{H}^1(\mathbb{S}^{n-1})$: Given $\bfc\in\mathbb{R}^n$,
  $\mathcal{L}(\bfc)$ is the function in $\mathbb{S}^{n-1}$ defined by $\mathcal{L}(\bfc)(\omega)=\bfc\cdot\omega$ for
  $\omega\in\mathbb{S}^{n-1}$. A simple computation shows that $\partial_{\mbox{\small $\bfc$}}\mathcal{P}(0,0,0)\bfc=
  (-\mathcal{L}(\bfc),0)$ for all $\bfc\in\mathbb{R}^n$. Hence we have
\begin{equation}
  \partial_{\eta,\bfc}\mathcal{P}(0,0,0)(\zeta,\bfc)=(\mathcal{B}_1(\zeta)-\mathcal{L}(\bfc),\mathcal{B}_2(\zeta)), \quad
  \forall(\zeta,\bfc)\in C^{\infty}({\mathbb{S}}^{n-1})\times\mathbb{R}^n.
\end{equation}
%---(5.8)---
  From this expression we can easily show that $\partial_{\eta,\mbox{\small $\bfc$}}\mathcal{P}(0,0,0)$ is an isomorphism of
  $C^{m+\mu}({\mathbb{S}}^{n-1})\times\mathbb{R}^n$ onto $C^{m-3+\mu}({\mathbb{S}}^{n-1})\times\mathbb{R}^n$, by using some
  standard argument. Briefly speaking, this argument is as follows: Since $\mathcal{B}_1$ is a bijection of
  $\displaystyle\bigoplus_{k=0\atop k\neq 1}^{\infty}\mathscr{H}^k(\mathbb{S}^{n-1})$ (in the topology of
  $\mathscr{D}'({\mathbb{S}}^{n-1})$) onto itself, $\mathcal{B}_2$ is an isomorphism of $\mathscr{H}^1(\mathbb{S}^{n-1})$
  onto $\mathbb{R}^n$ when restricted to $\mathscr{H}^1(\mathbb{S}^{n-1})$, and $\mathcal{L}$ is an isomorphism of
  $\mathbb{R}^n$ onto $\mathscr{H}^1(\mathbb{S}^{n-1})$, it follows immediately that $\partial_{\eta,\mbox{\small $\bfc$}}
  \mathcal{P}(0,0,0): C^{m+\mu}({\mathbb{S}}^{n-1})\times\mathbb{R}^n\to C^{m-3+\mu}({\mathbb{S}}^{n-1})\times\mathbb{R}^n$
  has a trivial kernel. Since it is clear that $\mathcal{B}_1$ is also a third-order elliptic pseudo-differential operator
  on $\mathbb{S}^{n-1}$ (it is the composition of the linearization of the mean curvature operator which is a second-order
  elliptic partial differential operator with the Dirichlet-Newmann operator on ${\mathbb{S}}^{n-1}=(1/K)\Gamma_0$ which is
  a first-order elliptic pseudo-differential operator), it follows that the operator $(\zeta,\bfc)\mapsto(\mathcal{B}_1(\zeta)
  -\mathcal{L}(\bfc),\mathcal{B}_2(\zeta))$ is a Fredholm operator of index zero, and, consequently, the assertion that this
  operator is an injection of $C^{m+\mu}({\mathbb{S}}^{n-1})\times\mathbb{R}^n$ onto $C^{m-3+\mu}(\overline{D})\times\mathbb{R}^n$
  ensures that it is also an isomorphism as we have desired. An alternative argument to prove that
  $\partial_{\eta,\mbox{\small $\bfc$}}\mathcal{P}(0,0,0): C^{m+\mu}({\mathbb{S}}^{n-1})\times\mathbb{R}^n\to
  C^{m-3+\mu}({\mathbb{S}}^{n-1})\times\mathbb{R}^n$ is surjective is as follows: From the Fourier expansion expression of
  $\mathcal{B}_1$ we easily see that $\partial_{\eta,\mbox{\small $\bfc$}}\mathcal{P}(0,0,0): \mathscr{D}'({\mathbb{S}}^{n-1})
  \times\mathbb{R}^n\to\mathscr{D}'({\mathbb{S}}^{n-1})\times\mathbb{R}^n$ is a surjection, so that for any $(\upsilon,\bfb)
  \in C^{m-3+\mu}({\mathbb{S}}^{n-1})\times\mathbb{R}^n$ the equation $\partial_{\eta,\mbox{\small $\bfc$}}\mathcal{P}(0,0,0)(\zeta,\bfc)
  =(\upsilon,\bfb)$ has a unique solution $(\zeta,\bfc)\in\mathscr{D}'({\mathbb{S}}^{n-1})\times\mathbb{R}^n$. Since
  $\mathcal{B}_1$ is a third-order elliptic pseudo-differential operator, from the property $\upsilon\in C^{m-3+\mu}({\mathbb{S}}^{n-1})$
  it follows that $\zeta\in C^{m+\mu}({\mathbb{S}}^{n-1})$. Hence the desired assertion follows.

  Having proved that $D_{\eta,\mbox{\small $\bfc$}}\mathcal{P}(0,0,0)$ is an isomorphism of $C^{m+\mu}({\mathbb{S}}^{n-1})
  \times\mathbb{R}^n$ onto $C^{m-3+\mu}(\overline{D})\times\mathbb{R}^n$, it follows from (5.6) and the implicit function
  theorem that by choosing $\delta$ and $\delta'$ smaller when necessary, there exists a unique mapping $\mathcal{G}=
  (\mathcal{G}_1,\mathcal{G}_2):O_{\delta}\to O_{\delta'}'\times\mathbb{R}^n$ which is smooth such that $\mathcal{G}_1(0)=0$,
  $\mathcal{G}_2(0)=0$, and for any $\rho\in O_{\delta}$, by letting $\eta=\mathcal{G}_1(\rho)$ and $\bfc=\mathcal{G}_2(\rho)$,
  $(\eta,\bfc)$ is the unique solution of the following equation in $O_{\delta'}'\times\mathbb{R}^n$:
$$
  \mathcal{P}(\rho,\eta,\bfc)=0, \quad \mbox{i.e.} \quad \mathcal{N}(\eta,\bfc)\mathcal{Q}(\rho,\eta)=0%+\mathcal{L}(\bfc),
  \quad \mbox{and} \quad  \mathcal{J}(\rho,\eta)=0.
$$
  Now let $v=\mathcal{Q}(\rho,\eta)$. Then $(\eta,v,\bfc)$ is clearly a solution of the free-boundary problem (5.5) for given
  $\rho\in O_{\delta}$. Consequently, by letting $u=v\circ\Psi_{\rho,\eta}$, $(\eta,u,\bfc)$ is a solution of the free-boundary
  problem (5.3) for given $\rho\in O_{\delta}$. Since $\mathcal{G}_1$ and $\mathcal{G}_2$ are smooth, it follows that the
  mapping $\rho\mapsto(\eta,u,\bfc)$ is smooth. This completes the proof of Lemma 5.1. $\quad\Box$
\medskip

  {\em Remark}.\ \ Later we shall need the expression of $\partial_{\eta,\mbox{\small $\bfc$}}\mathcal{P}(0,0,0)^{-1}$. From
  (5.8) we easily see that for an arbitrary $(\upsilon,\bfb)\in C^{\infty}({\mathbb{S}}^{n-1})\times\mathbb{R}^n$, if we
  denote $(\zeta,\bfc)=\partial_{\eta,\mbox{\small $\bfc$}}\mathcal{P}(0,0,0)^{-1}(\upsilon,\bfb)$ and assume that
$$
  \upsilon(\omega)=\sum_{k=0}^{\infty}\sum_{l=1}^{d_k}b_{kl}Y_{kl}(\omega),
$$
  then $\bfc=-\sqrt{n}\sigma_n^{-\frac{1}{2}}(b_{11},b_{12},\cdots,b_{1n})$, and
\begin{align}
  \zeta(\omega)=&\frac{K^2}{\mu}\frac{R^{n-2}\!-\!K^{n-2}}{(n\!-\!2)R^{n-2}}b_{00}Y_{00}(\omega)+K^{-1}\bfb\cdot\omega
\nonumber\\
  &-\frac{K^2}{\mu}\sum_{k=2}^{\infty}\sum_{l=1}^{d_k}\frac{n\!-\!1}{\lambda_k\!-\!n\!+\!1}
  \frac{R^{2k+n-2}\!-\!K^{2k+n-2}}{kK^{2k+n-2}\!+\!(k\!+\!n\!-\!2)R^{2k+n-2}}b_{kl}Y_{kl}(\omega).
\end{align}
%---(5.9)---

  Let $\mathfrak{M}_0$ and $\mathfrak{M}$ be as before, i.e., $\mathfrak{M}=\dot{\mathfrak{S}}^{m+\mu}({\mathbb{R}}^n)$
  and $\mathfrak{M}_0=\dot{\mathfrak{S}}^{m+3+\mu}({\mathbb{R}}^n)$, where $m$ is a positive integer $\geqslant 2$ and
  $0<\mu<1$. Let $\Sigma\in\mathfrak{M}_0$. The {\em standard local chart} $(\mathcal{U}_{\Sigma},\varphi_{\Sigma})$
  of $\mathfrak{S}^{m+\mu}({\mathbf{R}}^n)$ at $\Sigma$ is defined as follows: Since $\Sigma\in\dot{\mathfrak{S}}^{m+3+\mu}({\mathbf{R}}^n)$,
  we see that $\Sigma$ is a $\dot{C}^{m+3+\mu}$-hypersurface and its normal field $\bfn$ is of $C^{m+2+\mu}$-class: $\bfn\in
  C^{m+2+\mu}(\Sigma,{\mathbf{R}}^n)$. Choose a positive number $\delta$ sufficiently small such that the map
$$
  (x,t)\mapsto x+t\bfn(x), \quad x\in\Sigma,\;\; |t|<\delta
$$
  is a $\dot{C}^{m+2+\mu}$-diffeomorphism of $\Sigma\times(-\delta,\delta)$ onto the $\delta$-neighborhood of $\Sigma$. Let
  $\mathcal{U}_{\Sigma}=\{S\in\mathfrak{M}:d(S,\Sigma)<\frac{1}{2}\delta\}$. For every $S\in\mathcal{U}_{\Sigma}$
  there exists a unique function $\rho\in\dot{C}^{m+\mu}(\Sigma)$ such that the following relation holds:
$$
  S=\{x+\rho(x)\bfn(x):x\in\Sigma\}.
$$
  We define $\varphi_{\Sigma}(S)=\rho$. This defines the map $\varphi_{\Sigma}:\mathcal{U}_{\Sigma}\to\dot{C}^{m+\mu}(\Sigma)$
  and therefore defines the local chart $(\mathcal{U}_{\Sigma},\varphi_{\Sigma})$ of $\mathfrak{S}^{m+\mu}({\mathbf{R}}^n)$
  at $\Sigma$. In this local chart, the tangent space $\mathcal{T}_{\Sigma}(\mathfrak{M})$ of $\mathfrak{M}$ at $\Sigma$ can
  be naturally identified with the Banach space $\dot{C}^{m+\mu}(\Sigma)$.
  It follows that if $I\subseteq{\mathbb{R}}$ is an open interval and $S:I\to\mathfrak{M}_0$ is a $C^1$ curve, then
$$
   S'(t)=V_n(\cdot,t) \quad \mbox{for}\;\; t\in I;
$$
  see \cite{Cui3} for details. (Note that in \cite{Cui3} instead of $\dot{\mathfrak{S}}^{m+\mu}({\mathbb{R}}^n)$ the discussion
  is made for $\dot{\mathfrak{D}}^{m+\mu}({\mathbb{R}}^n)$. Since $\dot{\mathfrak{S}}^{m+\mu}({\mathbb{R}}^n)$ can be naturally
  identified with $\dot{\mathfrak{D}}^{m+\mu}({\mathbb{R}}^n)$, all results for $\dot{\mathfrak{D}}^{m+\mu}({\mathbb{R}}^n)$
  obtained in \cite{Cui3} work for $\dot{\mathfrak{S}}^{m+\mu}({\mathbb{R}}^n)$). This shows that $V_n(\cdot,t)$ appearing on
  the left-hand side of $(1.7)_4$ can be explained as $S'(t)$.

  Let $\mathcal{M}_c$ be the $C^{\infty}$-submanifold of $\mathfrak{M}_0$ and $\mathfrak{M}$ consisting of all surface spheres
  in $\mathbb{R}^n$. For each $S_0\in\mathcal{M}_c$ we let $R$ be its radius and $K$ be the corresponding number such that the
  relation (1.8) holds. Let $\Gamma_0$ be the surface sphere of radius $K$ concentric to $S_0$, and $D$ be the annular domain
  in $\mathbb{R}^n$ enclosed by $S_0$ and $\Gamma_0$. Choose a small number $\delta>0$ such that the assertion of Lemma 5.1
  holds when $m$ there is replaced by $m\!+\!3$, and denote by $\mathcal{U}(S_0,\delta)$ the neighborhood of $S_0$ in
  $\mathfrak{M}_0$ consisting of all hypersurfaces in $\mathbb{R}^n$ of the form $r=R[1+\rho(\omega)]$, where $\rho\in
  \dot{C}^{m+3+\mu}({\mathbb{S}}^{n-1})$ and $\|\rho\|_{\dot{C}^{m+3+\mu}({\mathbb{S}}^{n-1})}<\delta$. We denote
$$
  \mathcal{O}=\displaystyle\bigcup_{S_0\in\mathcal{M}_c}\mathcal{U}(S_0,\delta).
$$
  We now define a vector field $\mathscr{F}$ in $\mathfrak{M}$ with domain $\mathcal{O}$ as follows: Given $S\in
  \mathcal{O}$, let $S_0,\Gamma_0,D,\delta$ be as above such that $S\in\mathcal{U}(S_0,\delta)$. Let $\rho\in
  \dot{C}^{m+3+\mu}({\mathbb{S}}^{n-1})$ with $\|\rho\|_{\dot{C}^{m+3+\mu}({\mathbb{S}}^{n-1})}<\delta$ be the function
  on ${\mathbb{S}}^{n-1}$ such that $S$ is the hypersurface $r=R[1+\rho(\omega)]$. It follows that the problem
  (5.3) has a unique solution $(\eta,u,\bfc)$ with $\eta\in\dot{C}^{m+3+\mu}({\mathbb{S}}^{n-1})$,
  $u\in\dot{C}^{m+1+\mu}(\overline{D}_{\rho,\eta})$ and $\bfc\in\mathbb{R}^n$. Moreover, the map $\rho\mapsto(\eta,u,\bfc)$
  from the neighborhood $\|\rho\|_{\dot{C}^{m+3+\mu}({\mathbb{S}}^{n-1})}<\delta$ of the origin of
  $\dot{C}^{m+3+\mu}({\mathbb{S}}^{n-1})$ to $\dot{C}^{m+3+\mu}({\mathbb{S}}^{n-1})\times
  \dot{C}^{m+1+\mu}(\overline{D}_{\rho,\eta})\times\mathbb{R}^n$ is smooth. We define
$$
  \mathscr{F}(S)=\partial_{\mbox{\footnotesize$\bfn$}}u|_{S}
  \in\dot{C}^{m+\mu}(S)=\mathcal{T}_{S}(\mathfrak{M}).
$$
  It follows that the problem (1.7) reduces into the following initial value problem of a differential equation (1.11)
  in the Banach manifold $\mathfrak{M}$. Note that by putting $\varphi(S)=\rho$ and $\mathcal{U}=\mathcal{U}(S_0,\delta)$,
  $(\mathcal{U},\varphi)$ is a regular local chart of $\mathfrak{M}$.
\medskip

  {\bf Lemma 5.2}\ \ {\em If for each $S_0\in\mathcal{M}_c$ the number $\delta$ is chosen sufficiently small then the
  equation $(1.11)_7$ is of parabolic type in $\mathcal{O}$, and the representation of the vector field $\mathscr{F}$
  in the local chart $(\mathcal{U},\varphi)$ of $\mathfrak{M}$ as prescribed above is smooth.}
\medskip

  {\em Proof}.\ \ The second assertion immediately follows from the fact that the mapping $\rho\mapsto (\eta,u)$ is smooth
  ensured by Lemma 5.1. In what follows we prove the first assertion of this lemma.

  Let $F$ be the representation of $\mathscr{F}$ in the local chart $(\mathcal{U},\varphi)$, i.e., $F(\rho)=
  \varphi'(\varphi^{-1}(\rho))\mathscr{F}(\varphi^{-1}(\rho))$. Let $O_{\delta}=\{\rho\in\dot{C}^{m+3+\mu}({\mathbb{S}}^{n-1})
  :\|\rho\|_{\dot{C}^{m+3+\mu}({\mathbb{S}}^{n-1})}<\delta\}$. We regard $O_{\delta}$ as a neighborhood of the origin of
  $X_0=\dot{C}^{m+3+\mu}({\mathbb{S}}^{n-1})$. By using Lemma 5.1 (with $m$ there replaced by $m\!+\!3$), we easily see
  $F\in C^{\infty}(O_{\delta},X)$, where $X=\dot{C}^{m+\mu}({\mathbb{S}}^{n-1})$. To give the expression of $F$, we introduce
  an operator $\mathcal{D}:C^{m+1+\mu}(\overline{E^1})\to C^{m+\mu}(S_0)$ as follows:
$$
  \mathcal{D}(\rho)v=[\partial_{{\mathbf{n}}}(v\circ\Psi_{1\rho})|_{S_{\rho}}]\circ\psi_{1\rho}^{-1} \quad\;\;
  \mbox{for}\;\,v\in C^{m+1+\mu}(\overline{E^1}).
$$
  Here the notations $E^1$, $\Psi_{1\rho}$ and etc. are similar as in Lemma 5.1. We further introduce an operator $G:
  C^{m+3+\mu}({\mathbb{S}}^{n-1})\times C^{m+3+\mu}({\mathbb{S}}^{n-1})\to C^{m+\mu}({\mathbb{S}}^{n-1})$ as follows:
  For arbitrary $\rho,\eta\in C^{m+3+\mu}({\mathbb{S}}^{n-1})$,
$$
  G(\rho,\eta)=-\theta_{\ast}\circ\mathcal{D}(\rho)\mathcal{Q}(\rho,\eta),
$$
  where $\mathcal{Q}(\rho,\eta)$ is similar as in Lemma 5.1, $\theta$ represents natural projection of $S_0$ onto
  $\mathbb{S}^{n-1}$, and $\theta_{\ast}$ represents the push-forward operator induced by $\theta$, i.e.,
  $\theta_{\ast}(f)=f\circ\theta^{-1}\in C^{m+\mu}({\mathbb{S}}^{n-1})$ for $f\in C^{m+\mu}(S_0)$. It follows that
$$
  F(\rho)=G(\rho,\mathcal{G}_1(\rho)), \qquad \forall\rho\in O_{\delta}.
$$
  Hence for $\xi\in C^{m+3+\mu}({\mathbb{S}}^{n-1})$ we have
\begin{align}
   F'(0)\xi=&\partial_{\rho}G(0,0)\xi+\partial_{\eta}G(0,0)\mathcal{G}_1'(0)\xi
\nonumber\\
   =&-\theta_{\ast}\circ \partial_rw_1|_{S_0}-\theta_{\ast}\circ \partial_rw_2|_{S_0},
\end{align}
%---(5.10)---
  where $w_1,w_2$ are solutions of the following problems respectively:
$$
\left\{
\begin{array}{ll}
   \Delta w_1=0 &\quad \mbox{in}\; D,\\
   w_1=-\frac{\gamma}{R}(\xi+\frac{1}{n\!-\!1}\Delta_{\omega}\xi)   &\quad  \mbox{on}\; S_0,\\
   w_1=0  &\quad  \mbox{on}\; \Gamma_0,
\end{array}
\right.
$$
$$
\left\{
\begin{array}{ll}
   \Delta w_2=0 &\quad \mbox{in}\; D,\\
   w_2=0   &\quad  \mbox{on}\; S_0,\\
   w_2=-\frac{\mu}{K}(\zeta+\frac{1}{n\!-\!1}\Delta_{\omega}\zeta)  &\quad  \mbox{on}\; \Gamma_0,
\end{array}
\right.
$$
  where $\zeta=\mathcal{G}_1'(0)\xi$. Similarly as in the proof of Lemma 5.1, we see that if $\xi$, $\zeta$ has spherical
  harmonics expansions
\begin{equation}
  \xi(\omega)=\sum_{k=0}^{\infty}\sum_{l=1}^{d_k}c_{kl}Y_{kl}(\omega), \qquad
  \zeta(\omega)=\sum_{k=0}^{\infty}\sum_{l=1}^{d_k}c'_{kl}Y_{kl}(\omega),
\end{equation}
%---(5.11)---
  respectively, then
$$
  w_1(r,\omega)=\frac{\gamma}{R}\sum_{k=0}^{\infty}\sum_{l=1}^{d_k}\frac{\lambda_k\!-\!n\!+\!1}{n\!-\!1}
  \Big(\frac{r}{R}\Big)^k\frac{R^{2k+n-2}}{R^{2k+n-2}\!-\!K^{2k+n-2}}
  \Big[1-\Big(\frac{K}{r}\Big)^{2k+n-2}\Big]c_{kl}Y_{kl}(\omega)   \quad
   \mbox{for}\;\;  K\leqslant r\leqslant R,
$$
$$
  w_2(r,\omega)=\frac{\mu}{K}\sum_{k=0}^{\infty}\sum_{l=1}^{d_k}\frac{\lambda_k\!-\!n\!+\!1}{n\!-\!1}
  \Big(\frac{r}{K}\Big)^k\frac{K^{2k+n-2}}{R^{2k+n-2}\!-\!K^{2k+n-2}}
  \Big[\Big(\frac{R}{r}\Big)^{2k+n-2}-1\Big]c'_{kl}Y_{kl}(\omega)   \quad
   \mbox{for}\;\;  K\leqslant r\leqslant R.
$$
  Hence
\begin{equation}
  \theta_{\ast}\circ \partial_rw_1|_{S_0}=\frac{\gamma}{R^2}\sum_{k=0}^{\infty}\sum_{l=1}^{d_k}\frac{\lambda_k\!-\!n\!+\!1}{n\!-\!1}
  \Big(\frac{kR^{2k+n-2}\!+\!(k\!+\!n\!-\!2)K^{2k+n-2}}{R^{2k+n-2}\!-\!K^{2k+n-2}}\Big)c_{kl}Y_{kl}(\omega),
\end{equation}
%---(5.12)---
\begin{equation}
  \theta_{\ast}\circ \partial_rw_2|_{S_0}=-\frac{\mu}{K^2}\sum_{k=0}^{\infty}\sum_{l=1}^{d_k}\frac{\lambda_k\!-\!n\!+\!1}{n\!-\!1}
  \Big(\frac{(2k\!+\!n\!-\!2)R^{k-1}K^{k+n-1}}{R^{2k+n-2}\!-\!K^{2k+n-2}}\Big)c'_{kl}Y_{kl}(\omega),
\end{equation}
%---(5.13)---
  Note that the operators $\xi\mapsto\theta_{\ast}\circ\partial_rw_1|_{S_0}$ and $\zeta\mapsto\theta_{\ast}\circ
  \partial_rw_2|_{S_0}$ are Fourier multipliers of order $3$ and $-\infty$, respectively (the latter is ensured by the
  fact that $(K/R)^k$ converges to zero in exponential rate as $k\to\infty$). Note also that
$$
  \partial_{\rho}\mathcal{P}(0,0,0)\xi=(\theta_{\ast}\circ \partial_rw_1|_{\Gamma_0},-(R/K)\mathcal{B}_2(\xi)),
$$
  and $\zeta=\mathcal{G}_1'(0)\xi$ is the first component of $-\partial_{\eta,\bfc}\mathcal{P}(0,0,0)^{-1}
  \partial_{\rho}\mathcal{P}(0,0,0)\xi$. Clearly,
\begin{equation}
  \theta_{\ast}\circ \partial_rw_1|_{\Gamma_0}=\frac{\gamma}{R^2}\sum_{k=0}^{\infty}\sum_{l=1}^{d_k}\frac{\lambda_k\!-\!n\!+\!1}{n\!-\!1}
  \Big(\frac{(2k\!+\!n\!-\!2)K^{k-1}R^{k+n-1}}{R^{2k+n-2}\!-\!K^{2k+n-2}}\Big)c_{kl}Y_{kl}(\omega),
\end{equation}
%---(5.14)---
  From the above expression of $\partial_{\rho}\mathcal{P}(0,0,0)$ and the expression (5.9) of
  $\partial_{\eta,\mbox{\small $\bfc$}}\mathcal{P}(0,0,0)$
  we easily see that the operator $\xi\mapsto\zeta$ given by the relation $\zeta=\mathcal{G}_1(0)\xi$ is a zeroth order
  Fourier multiplier, so that the operator $\xi\mapsto\theta_{\ast}\circ\partial_rw_2|_{S_0}$ is a Fourier multiplier of
  order $-\infty$, and consequently a compact operator of $\dot{C}^{m+\mu}({\mathbb{S}}^{n-1})$ into itself. Now, since the
  operator $\xi\mapsto\theta_{\ast}\circ\partial_rw_1|_{S_0}$ is clearly a sectorial operator in
  $\dot{C}^{m+\mu}({\mathbb{S}}^{n-1})$ with domain $\dot{C}^{m+3+\mu}({\mathbb{S}}^{n-1})$, due to the fact that it is a
  third-order elliptic pseudo-differential operator with spectrum lying in the left-half part of the real axis, it
  follows immediately that $F'(0)$ is a sectorial operator. Now, by using a standard perturbation theorem, we see that
  $F'(\rho)$ is also a sectorial operator if $\|\rho\|_{\dot{C}^{m+3+\mu}({\mathbb{S}}^{n-1})}$ is sufficiently small,
  which proves the lemma.  $\quad\Box$
\medskip

  Let $G_{tl}$ and $G_{\!dl}$ be as in Section 1. We now prove that (1.12) and (1.13) hold.
\medskip

  {\bf Lemma 5.3}\ \ {\em Let notation be as in Lemma 5.2. The vector field $\mathscr{F}$ is invariant under the translation
  group action $(G_{tl},p)$, and quasi-invariant under the dilation group action $(G_{dl},q)$ with quasi-invariant factor
  $\theta(\lambda)=\lambda^{-3}$, $\lambda>0$, i.e., the relations $(1.12)$ and $(1.13)$ hold.}
\medskip

  {\em Proof}.\ \ Indeed, since the equations $(1.7)_1$--$(1.7)_6$ are invariant under translation of coordinate, it follows
  that given $S_0\in\mathfrak{M}_0$, if $S(t)$ is a solution of the equation $(1.11)_1$ with initial data $S(0)=S_0$, then
  for any $z\in G_{tl}$, $S_1(t)=p(z,S(t))$ is also a solution of $(5.11)_1$, but with initial data $S_1(0)=p(z,S_0)$. This
  implies that the following relation holds:
$$
  \frac{d}{dt}\big(p(z,S(t))\big)=\mathscr{F}(p(z,S(t))) \quad \mbox{for}\;\; t>0.
$$
  Since
$$
  \frac{d}{dt}\big(p(z,S(t))\big)=\partial_{S}p(z,S(t))S'(t)=\partial_{S}p(z,S(t))\mathscr{F}(S(t)),
$$
  we get
$$
  \partial_{S}p(z,S(t))\mathscr{F}(S(t))=\mathscr{F}(p(z,S(t))) \quad \mbox{for}\;\; t>0.
$$
  Letting $t\to 0^+$, we see that (1.12) follows. The proof of (1.13) is similar. $\quad\Box$
\medskip

  From \cite{Cui3} we know that the Lie group actions $(G_{tl},p)$ and $(G_{dl},q)$ satisfy the conditions
  $(L1)$--$(L7)$. Concerning the other conditions in Theorem 1.1, we have the following preliminary result:
\medskip

  {\bf Lemma 5.4}\ \ {\em Let notation be as in Lemma 5.2 and its proof. We have the following assertions:

  $(1)$\ \ $F'(0)$ has the following expression: If $\xi\in C^{\infty}(\mathbb{S}^{n-1})$ has an expression as given in
  $(5.11)$, then
$$
  F'(0)\xi=\sum_{k=2}^{\infty}\sum_{l=1}^{d_k}\mu_kc_{kl}Y_{kl}(\omega),
$$
  where
\begin{align*}
   \mu_k=&\displaystyle-\frac{\gamma}{R^2}\frac{\lambda_k\!-\!n\!+\!1}{n\!-\!1}
   \frac{k(k\!+\!n\!-\!2)(R^{2k+n-2}\!-\!K^{2k+n-2})}{kK^{2k+n-2}\!+\!(k\!+\!n\!-\!2)R^{2k+n-2}}\\
   =&-\frac{\gamma}{R^2}\frac{\lambda_k\!-\!n\!+\!1}{n\!-\!1}
  \frac{(k\!+\!n\!-\!2)[1\!-\!(\mu/\gamma)^{2k+n-2}]}{1\!+\!k^{-1}(n\!-\!2)\!+\!(\mu/\gamma)^{2k+n-2}},\qquad k=2,3,\cdots.
\end{align*}

  $(2)$\ \ $\sigma(F'(0))=\{0\}\cup\{\mu_k:k=2,3,\cdots\}$.

  $(3)$\ \ ${\rm Ker}\,F'(0)=\mathscr{H}^0(\mathbb{S}^{n-1})\cup\mathscr{H}^1(\mathbb{S}^{n-1})$.
  Hence $\dim{\rm Ker}\,F'(0)=n+1$.

  $(4)$\ \ $F'(0)$ is a standard Fredholm operator.}
\medskip

  {\em Proof.}\ \ The assertion (1) follows easily from (5.10), (5.12), (5.13), (5.14) and (5.9). The assertions (2)$\sim$(4)
  are immediate consequences of the assertion (1). $\quad\Box$
\medskip

  {\em Proof of Theorem 1.2.}\ \ It is easy to see that $\{\mu_k\}_{k=2}^{\infty}$ is a strictly monotone decreasing
  sequence of negative numbers, so that $\sup\{{\rm Re}\lambda:\lambda\in\sigma(F'(0))\backslash\{0\}\}=\mu_2<0$. This shows
  that all conditions (with $N=2$) of Theorem 1.1 are satisfied. Hence, by applying Theorem 1.1 we obtain Theorem 1.2.$\quad\Box$
\medskip

%   {\bf Acknowledgement}.\hskip 1em

\end{document}